\documentclass[12pt]{article}

\usepackage{amsfonts}
\usepackage{amsmath}
\usepackage{amsthm}
\usepackage{graphics}
\usepackage{setspace}
\usepackage{makeidx}
\usepackage{setspace}

\makeindex

\theoremstyle{definition}

\theoremstyle{definition}

\theoremstyle{remark}

\theoremstyle{plain}

\newcommand{\href}[1]{#1} 

\usepackage{amsmath,amsthm, amsfonts, amssymb,amstext, latexsym} 
\usepackage{graphicx} 
\usepackage[bindingoffset=0.125in,centering,includeheadfoot,margin=1in]{geometry}
\usepackage[english]{babel}
\usepackage{graphics}
\usepackage[mathscr]{eucal}
\usepackage{tikz}
\usetikzlibrary{arrows,snakes,backgrounds}
\usepackage[labelfont=bf]{caption}
\usepackage{float}

\newcommand{\pr}[1]{\!\left(#1\right)}
\newcommand{\fr}[2]{\frac{#1}{#2}}

\newcommand{\lr}[1]{\left|\;\!#1\;\!\right|}

\newcommand{\mb}[1]{\mathbb{#1}}

\newcommand{\rt}{\rightarrow}

\newcommand{\tb}{\textbf}

\newcommand{\tx}[1]{\text{#1}}

\newcommand{\ep}{\epsilon}

\newcommand{\Dl}[1]{\Delta{#1}}

\newcommand{\sq}{\sqrt}

\newcommand{\mc}[1]{\mathcal{#1}}

\newcommand{\s}{\sigma}

\newcommand{\p}{\prod}

\newcommand{\eq}{\equiv}
\newcommand{\B}{O}

\newcommand{\ch}{\chi}

\newcommand{\z}{\zeta}

\newcommand{\gm}{\gamma}
\newcommand{\Gm}{\Gamma}
\newcommand{\Q}{\mb{Q}}
\newcommand{\Z}{\mb{Z}}
\newcommand{\C}{\mb{C}}

\newcommand{\kr}[2]{\pr{\fr{#1}{#2}}}

\newcommand{\msum}{\sum_{m\leq \fr{X}{\Dl{X}}}}
\newcommand{\asum}{\sum_{a \leq \sq{\fr{m\Dl{X}}{3}}}}
\newcommand{\bsum}{\sum_{0 \leq b \leq a}}
\newcommand{\pasum}{\sum_{a\leq x}}

\newcommand{\Digits}{\tx{Digits}}

\newcommand{\sumd}{\sum_{d\in D(\infty)}}

\title{Conjectures and experiments concerning the moments of $L(1/2,\chi_d)$.}
\author{Matthew W. Alderson and Michael O. Rubinstein
\footnote{Support for work on this paper was provided by the
National Science Foundation under awards DMS-0757627 (FRG grant),
an NSERC Discovery Grant, and an OGS Scholarship.}}
\date{}

\begin{document}

\maketitle

\begin{abstract}
    We report on some extensive computations and experiments
    concerning the moments of quadratic Dirichlet $L$-functions at the critical
    point. We computed the values of $L(1/2,\chi_d)$ for $- 5\times 10^{10} < d
    < 1.3 \times 10^{10}$ in order to numerically test conjectures concerning
    the moments $\sum_{|d|<X} L(1/2,\chi_d)^k$. Specifically, we tested
    the full asymptotics for the moments conjectured
    by Conrey, Farmer, Keating,
    Rubinstein, and Snaith, as well as the conjectures of Diaconu,
    Goldfeld, Hoffstein, and Zhang concerning additional
    lower terms in the moments. We also describe the algorithms used
    for this large scale computation.
\end{abstract}


\section{Introduction}

Let $D$ be a squarefree integer, $D\neq 0,1$, and let $K= \Q(\sqrt{D})$
be the corresponding quadratic field. The fundamental discriminant $d$ of $K$ equals
$D$ if $D = 1 \mod 4$, and $4D$ if $D = 2,3 \mod 4$.
Let $\chi_d(n)$ be the Kronecker symbol $\kr{d}{n}$,
and $L\pr{s,\ch_d}$ the quadratic Dirichlet $L$-function given by the
Dirichlet series
\begin{equation}
    L\pr{s,\ch_d} = \sum_{n=1}^\infty \fr{\chi_d(n)}{n^s}, \qquad \Re(s)>0,
\end{equation}
satisfying the functional equation
\begin{equation}
    L\pr{s,\ch_d} = \lr{d}^{\fr{1}{2} - s} X(s,a) L\pr{1-s, \ch_d},
\end{equation}
where
\begin{equation}
    \label{eq:X}
    X(s,a) =
    \pi^{s-\fr{1}{2}} \fr{\Gm\pr{\fr{1-s+ a}{2}}}{\Gm\pr{\fr{s+a}{2}}}, \qquad a =
    \begin{cases} 0 \quad &\tx{if $d> 0$,} \\ 1 \quad &\tx{if $d<0$.} \end{cases}
\end{equation}

In this paper, we describe some experiments concerning
the moments of quadratic Dirichlet $L$-functions at the central critical point:
\begin{equation}
  \sum_{d \in D(X)} L\pr{1/2,\ch_d}^k. \label{eq: integralmoments}
\end{equation}
Here, $k$ is a positive integer, and $D(X)$ denotes the set of fundamental discriminants
with $|d|< X$.

Several conjectures exist for these moments. For instance,
Keating and Snaith \cite{keat: snaith}, motivated by the fundamental work of
Katz and Sarnak \cite{katz: sarnak} and based on an analogous result in Random
Matrix Theory, conjectured a formula for the leading asymptotics of (\ref{eq: integralmoments}).
Specifically, they conjectured
that, as $X\rt\infty$,
\begin{equation}
  \fr{1}{\lr{D(X)}} \sum_{d\in D(X)} L(1/2,\ch_d)^k \sim
  a_k \p_{j=1}^k \fr{j!}{(2j)!} \log(X)^{\fr{k(k+1)}{2}}, \label{eq: KeatingSnaith}
\end{equation}
where $a_k$ is an arithmetic factor, described by Conrey and Farmer \cite{cf: meanLfunc},
of the form
$$
    a_k = \p_p \fr{\pr{1 - \fr{1}{p}}^{\fr{k(k+1)}{2}}}{1 + \fr{1}{p}} \pr{\fr{\pr{1 -
    \fr{1}{\sq{p}}}^{-k} + \pr{1 + \fr{1}{\sq{p}}}^{-k}}{2} + \fr{1}{p}}.
$$
In a few cases, Keating and Snaith's conjecture
agrees with known theorems, e.g. Jutila for
$k =1,2$ \cite{jutila: character}, and Soundararajan for $k=3$ \cite{sound: quadDir}.

Subsequently, a more precise asymptotic expansion for (\ref{eq: integralmoments}) was predicted by
Conrey, Farmer, Keating, Rubinstein, and Snaith \cite{cfkrs: intmom}. The lower terms are affected
by the form of the Gamma factors in the functional equation for $L(s,\chi_d)$, so naturally they
considered the subset of $d>0$ separately from $d<0$. Therefore, let
\begin{eqnarray}
    D_+(X) &=& \{ d \in D(X) : d>0\} \notag \\
    D_-(X) &=& \{ d \in D(X) : d<0\}.
\end{eqnarray}
The conjecture of CFKRS states that:
\begin{equation}
  \label{eq: moment asympt}
  \sum_{d \in D\pm(X)} L\pr{1/2,\ch_d}^k \sim
  \frac{3}{\pi^2} X  \mathcal Q_\pm(k,\log X),
\end{equation}
where $\mathcal Q_\pm(k,x)$ is a polynomial of degree $k(k+1)/2$ in $x$.
The fraction $3/\pi^2$ accounts for the density of fundamental discriminants amongst all the integers.
In their paper, CFKRS also conjectured a remainder term of size $O(X^{1/2+\ep})$ but the evidence,
both theoretical and numerical (discussed below), suggests the existence of lower order terms
for $k \geq 3$.

The beauty in their conjecture lies in the fact that it gives a formula for
the polynomial in question.
The polynomial $\mathcal Q_\pm(k,\log X)$ is expressed
in terms of a more fundamental polynomial $Q_\pm(k,x)$ of the same degree
that captures the moments locally:
\begin{equation}
  \mathcal Q_\pm(k,\log{X}) = \frac{1}{X} \int_1^X Q_\pm(k,\log{t}) dt.
\end{equation}
The polynomial $Q_\pm(k,x)$ is described below, and its leading coefficient
agrees with the conjecture of Keating and Snaith (\ref{eq: KeatingSnaith}).

The polynomial $Q_\pm(k,x)$ of CFKRS is given by them as the $k$-fold residue:
\begin{equation}\label{eq: main integral}
    Q_\pm(k,x)= \frac{(-1)^{k(k-1)/2} 2^k}{k!} \frac {1}{(2\pi i)^k}
    \oint \cdots \oint
    \frac{G_\pm(z_1,\ldots,z_k)\Delta(z_1^2,\ldots,z_k^2)^2}{\prod_{j=1}^k z_j^{2k-1}}
    e^{\frac x 2 \sum_{j=1}^k z_j} \,  dz_1\ldots dz_k
  \end{equation}
where
\begin{equation}
    G_\pm(z_1,\ldots,z_k) = A_k(z_1,\ldots,z_k) \prod_{j=1}^k X(\frac 1 2 +z_j,
    a)^{-1/2} \prod_{1\leq i\leq j \leq k} \zeta(1+z_i+z_j),
\end{equation}
and
\begin{equation}
    \Delta(z_1^2,\ldots,z_k^2)
    = \prod_{1\leq i < j \leq k} (z_j^2 -z_i^2)
\end{equation}
is a Vandermonde determinant.
Here, $a=0$ for $G_+$ and $a=1$ for $G_-$, $X(s,a)$ is given in~\eqref{eq:X},
and $A_k$ equals the Euler product, absolutely convergent for $|\Re z_j|
< 1/2$, defined by
\begin{multline}
    \index{$A_k$}A_k(z_1, \ldots, z_k) = \prod_p \prod_{1\leq i\leq j \leq k}
    \left(1-\frac{1}{p^{1+z_i+z_j}} \right)\\
    \times \left(\frac 1 2 \left( \prod_{j=1}^{k}\left(1-\frac 1
          {p^{\frac 1 2 +z_j}} \right)^{-1} +
        \prod_{j=1}^{k}\left(1+\frac 1{p^{\frac 1 2 +z_j}}
        \right)^{-1} \right) +\frac 1 p \right) \left(1+\frac 1 p
    \right)^{-1}.
  \end{multline}

More generally, CFKRS predicted that for suitable
weight functions $g$,
\begin{equation}
    \label{eq: cfkrs}
    \sum_{d \in D_\pm(\infty)}  L(1/2,\ch_d)^k g\pr{\lr{d}}
    \sim \frac{3}{\pi^2}
    \int_1^\infty Q_\pm(k,\log{t}) g(t) dt.
\end{equation}

The method that CFKRS used
to heuristically derive this formula relies on number theoretic techniques, specifically
the approximate functional equation, but was guided by analogous results in random matrix theory
to help determine the form of the conjecture.

An alternative approach for conjecturing moments exists. In their paper \cite{dgh: multDir},
Diaconu, Goldfeld, and Hoffstein used the double Dirichlet series
$$
     Z_k\pr{s,w} = \sumd \fr{L\pr{s,\ch_d}^k}{\lr{d}^w}
$$
to study the moments of $L(1/2,\chi_d)$.
In particular, they showed how one can derive a formula for the
cubic moments of $L\pr{1/2,\ch_d}$ by investigating the polar behavior of
$Z_3\pr{s,w}$.

The method of DGH produces a proof for the cubic moment of $L(1/2,\chi_d)$.
Specifically they show that the difference of both sides of \eqref{eq: moment asympt}
for $k=3$ is, for any $\epsilon>0$, of size $O_\epsilon(X^{\theta+\epsilon})$,
where $\theta=.85366\ldots$. They also gave a remainder term for a smoothed
cubic moment with $\theta=4/5$.
Recently, Young~\cite{young2} has obtained an improved
estimate for the remainder term of size $O(X^{3/4+\ep})$. The moments he
considered were smoothed, and, for simplicity, he considered the subset of
discriminants divisible by 8 and positive. The appearance of $3/4+\ep$ is very
interesting in light of the next paragraph.

The method of DGH predicts the existence of
a further lower order term of size $X^{3/4}$. In particular,
DGH conjectured that there exists a constant $b$ such that
\begin{equation}
     \sum_{d\in D(X)} L(1/2,\ch_d)^3
     = \frac{6}{\pi^2} X\mc{Q}\pr{3,\log X} + bX^{\fr{3}{4}} + \B\pr{X^{\fr{1}{2} + \ep}}.  \label{eq: 421}
\end{equation}
The existence of such a term comes from a pole of the double
Dirichlet series at $w=3/4$ and $s=1/2$,
the conjectured meromorphic continuation of
$Z_3\pr{1/2,w}$ to $\Re(w) < 3/4$, and assumes a growth condition on
$Z_3\pr{1/2,w}$.

DGH also suggest that additional lower order terms, infinitely many for each $k\geq 4$,
are expected to persist. The form of these terms is described in Zhang's
survey \cite{zhang: appDir}, along with an exposition of the approach
using double Dirichlet series. Their conjecture
involving lower terms is stated in the following form. For $k\geq 4$, and every $\ep>0$,
\begin{equation}
  \sum_{d \in D(X)} L\pr{1/2,\ch_d}^k= 
  \sum_{l=1}^\infty X^{(l+1)/(2l)} P_l(\log{x}) + O(X^{1/2+\ep}),
  \label{eq: integralmoments_k}
\end{equation}
where every $P_l$ is a polynomial depending on $k$. In particular, $P_1$ is a polynomial
of degree $k(k+1)/2$, presumably agreeing with the polynomial predicted by the CFKRS
conjecture, but to our knowledge this agreement has not been checked.

Zhang \cite{zhang: cubmom} further conjectured that $b \approx - .2154$,
and, in a private communication to one of the authors, reported that
he also computed the constants associated with the $X^{3/4}$ term
when one restricts to $d<0$, or to $d>0$, thus predicting:
\begin{eqnarray}
     \sum_{d\in D_\pm(X)} L(1/2,\ch_d)^3
     &=& \frac{3}{\pi^2} X\mathcal Q_\pm(3,\log X) +
     b_\pm X^{\fr{3}{4}} + \B\pr{X^{\fr{1}{2} + \ep}},
\end{eqnarray}
with $b_+ \approx -.14$ and $b_- \approx -.07$ (note that $b=b_+ + b_-$).
His evaluation of $b$, $b_+$, and $b_-$ involves an elaborate sieving process,
and also depends on unproven hypotheses regarding the meromorphic continuation
and rate of growth of $Z_3$.

While it might seem that a term as large as $X^{3/4}$ in the cubic moment
should be easily detected, two things make it very challenging in this context:
the small size of the constants involved, and also the fact that the remainder
term, conjecturally of size $O(X^{1/2+\ep})$, dominates even in our large
data set - presumably the $X^\ep$ can get as large as some power of $\log(x)$,
which can dominate $X^{1/4}$ even for values of $X$ as large as $10^{11}$.

For this reason, we embarked on a large scale
computation in order to see whether such a lower main term in the cubic moment
could be detected or not. We also carried out extensive
verification of the predictions of CFKRS for $k=1,\ldots,8$. While CFKRS provided
some modest data in \cite{cfkrs: intmom}, for $|d|<10^7$, we carried out tests
for $-5\times10^{10}<d<0$ and $0<d<1.3 \times 10^{10}$.
In order to dampen the effect of the noisy remainder term, we also considered
smoothed moments.

Our numerical results are described in Section 2. Interestingly, they
lend support to both the full asymptotic expansion conjectured by CFKRS,
and to the existence of lower terms predicted by DGH and Zhang.

In Section 3 we describe the two methods that we used to compute a large number
of  $L(1/2,\chi_d)$, for $d<0$ and, separately, for $d>0$. The first, for
$d<0$, is based on the theory of binary quadratic forms, and uses Chowla and Selberg's
$K$-Bessel expansion of the Epstein zeta function~\cite{cs: epzeta}. The
second, for $d>0$, uses a traditional smooth approximate functional equation.
Both methods have comparable runtime complexities, but the former has the
advantage of being faster by a constant factor. See the
end of the paper for a discussion that compares the two runtimes.

\section{Numerical Data}

In this section, we numerically examine the conjectures of CFKRS, DGH, and Zhang
for the moments of $L(1/2,\ch_d)$.

The collected data provides further evidence in favour of the
CFKRS conjecture concerning the full asymptotics of the moments of
$L(1/2,\chi_d)$. With respect to the remainder term, the numerics also seem
to suggest the presence of additional lower terms as predicted by DGH and
Zhang.

In Tables 1--2 and Figures 1--4 we depict the quantities
\begin{equation}
    R_\pm(k,X) := \fr{\displaystyle \sum_{d\in D_{\pm}(X)} L(1/2,\ch_d)^k}
    {\displaystyle \frac{3}{\pi^2}\int_1^X Q_\pm(k,\log{t}) dt} , \label{eq: quotient}
\end{equation}
and the related difference
\begin{equation}
    \Delta_\pm(k,X) := \sum_{d\in D_{\pm}(X)} L(1/2,\ch_d)^k\;\;
    - \frac{3}{\pi^2}\int_1^X Q_\pm(k,\log{t}) dt ,
    \label{eq: difference}
\end{equation}
for $k = 1,\ldots, 8$ and both positive and negative discriminants $d$.

The quantity \eqref{eq: quotient} measures the consistency of CFKRS prediction,
while~\eqref{eq: difference} allows one to see the associated remainder term.
The numerator of (\ref{eq: quotient}) was calculated by computing many values of
$L(1/2,\ch_d)$, using the methods described in the next two sections. The
denominator was obtained from numerically approximated values of the coefficients of
$Q_\pm(k,\log t)$, computed
in the same manner performed in \cite{cfkrs: intmom}, though to higher precision. Tables of the
coefficients of the polynomials $Q_\pm(k,x)$ can be found in \cite{cfkrs: intmom}.
These values were then also used in graphing the difference (\ref{eq: difference}).

Tables 1--2 provide strong numerical support in favor of the asymptotic
formula predicted by CFKRS, described in equations~\eqref{eq: moment asympt}
to~\eqref{eq: main integral}, for both $d<0$ and $d>0$, agreeing to 7--8 decimal places
for $k=1$, and 4--5 decimal places for $k=8$.

In the figures below we depict thousands of values of the
quantities~\eqref{eq: quotient} and~\eqref{eq: difference}
at multiples of $10^7$, i.e. $X=10^7,2\times10^7,\ldots$. We display data up to
$X=1.3\times10^{10}$ for $d>0$, and $X=5\times 10^{10}$ for $d<0$. The larger
amount of data for $d<0$ reflects the faster method that we used for computing
the corresponding $L$-values.

In Figures~\ref{fig: fig1} and~\ref{fig: fig2}, notice that each graph fluctuates tightly about
one, with the extent of fluctuation becoming progressively larger as
$k$ increases, as indicated by the varying vertical scales. The graphs show
excellent agreement with the full asymptotics as predicted by CFKRS across all
eight moments computed, for both $d<0$ and $d>0$. One does also notice a
slight downward shift from $1$ in the $k=3$ plots, as predicted by DGH and Zhang.

We also depict in Figures~\ref{fig: fig3} and~\ref{fig: fig4} the
differences~\eqref{eq: difference} as dots, as well as the running average of
the plotted differences as a solid curve. While we plot the average every
$10^7$, these running averages were computed by sampling the differences every
$10^6$. We chose to display $1/10$th of the computed values in order to make
our plots more readable given the limited resolution of computer displays and
printers.

Averaging has the effect of smoothing the moment and reducing the impact
of the noisy remainder term. It allows one to more clearly see if
there are any biases hiding within the noise. This running average gives a
discrete approximation to the smoothed difference:
\begin{equation}
    \Delta_\pm(k,X) := \sum_{d\in D_{\pm}(X)} L(1/2,\ch_d)^k (1-|d|/X)\;\;
    - \frac{3}{\pi^2}\int_1^X Q_\pm(k,\log{t}) (1-t/X) dt , \label{eq: smoothed difference}
\end{equation}

We make several observations concerning Figures~\ref{fig: fig3} and~\ref{fig:
fig4}. First, for $k=1$, the observed remainder term is seemingly of size
$X^{1/4+\ep}$, much smaller than Goldfeld and Hoffstein's proven bound of
$O(X^{19/32+\ep})$ for the first moment, and also smaller than the bound of
$O(X^{1/2+\ep})$ implied in their work for a smoothed first moment (for the
latter, see also Young~\cite{young1}). Furthermore, there appears to be a bias
in the average remainder term. This is especially apparent for $d<0$ and
further supported by our log log plots in Figure~\ref{fig: fig6}.

For $k=2$, the remainder term, even when averaged, fluctuates above and below
0. The largest average remainder for $k=2$ in our data set was of size roughly
$1.3\times10^4$, consistent with a conjectured remainder, for $k=2$, of size
$O(X^{1/2+\ep})$.

In Figures~\ref{fig: fig3} and~\ref{fig: fig4}, a bias of the kind predicted by
DGH can be seen. For $k\geq 3$ a noticeable bias is evident in the remainder
term, especially when averaged, and most prominently for $d>0$.

We elaborate on this last point further.
In the case of $k=3$, DGH predict a single main lower term of
the form $b X^{3/4}$ and, as described in the introduction, Zhang worked out
the value of $b$, separately for $d<0$ and $d>0$, as equal to $-.07$ and $-.14$
respectively.

One possible explanation for the fact that the bias appears more prominently
for $d>0$, even though we have less data in that case, is that the `noise'
appears numerically to be much larger for $d<0$. Comparing Figure~\ref{fig:
fig3} with Figure~\ref{fig: fig4}, we notice that the plotted remainder terms
seem to be about ten times larger for $d<0$ as compared to $d>0$, even when restricted
to $|d|<1.3\times 10^{10}$. A larger amount of `noise' in the remainder term
makes it harder to detect a lower order term hiding within the noise,
especially when the the lower terms are married to such small constant factors:
$-.07 X^{3/4}$ and $-.14 X^{3/4}$, as predicted by Zhang, before averaging, and
$4/7$ as large after averaging over $X$.

Thus, even though one expects, in the long run, to see an $X^{3/4}$ term
dominate over noise of size $X^{1/2+\ep}$, in the ranges examined it seems that
the noise has a large impact.

Another factor possibly affecting the poorer quality of the lower term detected
when $d<0$ is that the predicted constant factor is about one half as large:
$-.07$ for $d<0$, compared to $-.14$ for $d>0$. Combined with noise that is ten
times bigger, it is not surprising that the quality of the average remainder
term for $d<0$ seems to be more affected by the noise.

In Figure \ref{fig: fig5} we redisplay the $k=3$ plots from Figures~\ref{fig:
fig3} and~\ref{fig: fig4}, zoomed in to allow one to see the average
remainder term in greater detail. Here we also depict the prediction of Zhang
(dashed line). More precisely, the dashed line represents the average predicted
lower term:
\begin{equation}
    \frac{1}{X}\int_0^X b_\pm t^{3/4} dt = \frac{4}{7} b_\pm x^{3/4},
\end{equation}
with $b_-=-.07$ for $d<0$ and $b_+=-.14$ for $d>0$.
For $d>0$, the fit of the average remainder term against Zhang's prediction
is very nice. For $d<0$, the plot supports a bias in the sense that the average value
is mainly negative, but the fit against $-.07 \frac{4}{7} X^{3/4}$ is far
from conclusive.

Plots of the average reaminder term on a log log scale are shown in
Figure~\ref{fig: fig6}, for $1\leq k \eq 4$, and positive/negative $d$.
On a log log scale, a function of $X$ of the form $f(X) = B X^{3/4}$
is transformed into the function of $u=\log{X}$ given by
$\log{B} + 3/4 u$, i.e. a straight line with slope $3/4$. We compare, in the
$k=3$ plots, the log log plot of the average remainder against Zhang's
prediction. The fit is especially nice for $d>0$, but only in crude qualitative
terms for $d<0$, consistent with the observations made concerning the poorer
fit in the $k=3$, $d<0$, plots of Figures~\ref{fig: fig3}--\ref{fig: fig5}.

For $k \geq 4$, DGH predict infinitely many lower terms, with
the largest one of size $X^{3/4}$ times a power of $\log{X}$
which they did not make explicit.
Interestingly, a bias in support of this does appear evident,
especially for $k=4$ and $d>0$. In that log log plot, the remainder term does
seem reasonably straight suggesting a lower order term obeying a power law,
perhaps with some additional powers of log.

It is reasonable to contest that the observed biases here exist due to
a persistent small error in the calculation of the moment polynomials or in the
values of $L(1/2,\chi_d)$. In an effort to alleviate such concerns, the
computations yielding our numerics were executed again, in a limited way, using
higher precision. As anticipated, these higher precision results remained
consistent with the initial results, reducing the possibility of such a bias
existing. Furthermore, the overall excellent agreement of the computed moments
with the predicted asymptotic formula of CFKRS supports the correctness
of the computation.

\subsection{Tables and Figures}
\begin{table}[H]
\centerline{\fontsize{11pt}{12pt}\selectfont
\begin{tabular}{|c|c|c|c|c|}
\hline
$k$ & $\sum_{d\in D_-(X)} L(1/2,\chi_d)^k$ & $\frac{3}{\pi^2} \int_1^X Q_-(k,\log{t}) dt$& $R_-(k,x)$ & $\Delta_-(k,X)$ \\ \hline
1 & 25458527125.376 & 25458526443.085 & 1.0000000268001 & 682.291 \\ 
1 & 52401254983.398 & 52401252573.351 & 1.0000000459922 & 2410.047 \\ 
1 & 79904180421.746 & 79904180600.902 & .99999999775786 & -179.156 \\ 
1 & 107770905413.09 & 107770904521.07 & 1.0000000082770 & 892.02 \\ 
1 & 135908144579.9 & 135908144595.65 & .99999999988411 & -15.75 \\ 
\hline
2 & 695798091128.96 & 695797942880.62 & 1.0000002130623 & 148248.34 \\ 
2 & 1505736931971.7 & 1505736615082.0 & 1.0000002104549 & 316889.7 \\ 
2 & 2362905062077.2 & 2362905209666.9 & .99999993753888 & -147589.7 \\ 
2 & 3251727763805.6 & 3251727486319.2 & 1.0000000853351 & 277486.4 \\ 
2 & 4164586513531.5 & 4164586544704.8 & .99999999251467 & -31173.3 \\ 
\hline
3 & 35923488939396. & 35923434720074. & 1.0000015093023 & 54219322. \\ 
3 & 82792501873632. & 82792433101707. & 1.0000008306547 & 68771925. \\ 
3 & .13470723693602e15 & .13470723096090e15 & 1.0000000443563 & 5975116 \\ 
3 & .19013982678941e15 & .19013979175101e15 & 1.0000001842770 & 35038394 \\ 
3 & .24831500039182e15 & .24831501538879e15 & .99999993960505 & -14996973 \\ 
\hline
4 & .26221677201508e16 & .26221542614856e16 & 1.0000051326749 & 13458665240 \\ 
4 & .64846065425230e16 & .64845918799277e16 & 1.0000022611439 & 14662595290 \\ 
4 & .10987196470794e17 & .10987187884822e17 & 1.0000007814531 & .8585972e10 \\ 
4 & .15956123181403e17 & .15956125546013e17 & .99999985180550 & -.2364610e10 \\ 
4 & .21299535514803e17 & .21299540911015e17 & .99999974665125 & -.5396212e10 \\ 
\hline
5 & .23541937472178e18 & .23541622006477e18 & 1.0000134003384 & .315465701e13 \\ 
5 & .62771726711464e18 & .62771414322685e18 & 1.0000049766089 & .312388779e13 \\ 
5 & .11106890853615e19 & .11106862772711e19 & 1.0000025282480 & .28080904e13 \\ 
5 & .16628632428499e19 & .16628683849741e19 & .99999690767817 & -.51421242e13 \\ 
5 & .22724025077610e19 & .22724048423231e19 & .99999897264693 & -.23345621e13 \\ 
\hline
6 & .24225487162243e20 & .24224780818937e20 & 1.0000291578822 & .706343306e15 \\ 
6 & .69880224640908e20 & .69879554487455e20 & 1.0000095901220 & .670153453e15 \\ 
6 & .12937968210632e21 & .12937887586288e21 & 1.0000062316467 & .80624344e15 \\ 
6 & .19996752978479e21 & .19997013306315e21 & .99998698166411 & -.260327836e16 \\ 
6 & .28005925088677e21 & .28006019455853e21 & .99999663046810 & -.94367176e15 \\ 
\hline
7 & .274712571777e22 & .274697762672e22 & 1.00005391054 & .14809105e18 \\ 
7 & .859431066562e22 & .859415893116e22 & 1.00001765553 & .15173446e18 \\ 
7 & .166743403869e23 & .166740957095e23 & 1.00001467410 & .2446774e18 \\ 
7 & .266330275024e23 & .266339641978e23 & .999964830793 & -.9366954e18 \\ 
7 & .382588166641e23 & .382591322018e23 & .999991752617 & -.3155377e18 \\ 
\hline
8 & .3351697755e24 & .3351406841e24 & 1.000086804 & .290914e20 \\ 
8 & .1139465805e25 & .1139429048e25 & 1.000032259 & .36757e20 \\ 
8 & .2319359069e25 & .2319282301e25 & 1.000033100 & .76768e20 \\ 
8 & .3831454627e25 & .3831738559e25 & .9999259000 & -.283932e21 \\ 
8 & .5649093016e25 & .5649183210e25 & .9999840342 & -.90194e20 \\ 
\hline
\end{tabular}
}
\caption[Moments of $L(1/2,\ch_d)$ versus the CFKRS conjectured asymptotics,
$d<0$]{Moments $\sum_{d\in D_-(X)} L(1/2,\chi_d)^k$ versus CFKRS' $\frac{3}{\pi^2}
\int_1^X Q_-(k,\log{t}) dt$, for $k=1,\ldots,8$ and $d<0$. Five values for
each $k$ are shown, at $X = 10^{10},2\times 10^{10}, \ldots, 5\times
10^{10}$.}\label{tab:Lhalfchin}
\end{table}

\begin{table}[H]
\centerline{\fontsize{11pt}{12pt}\selectfont
\begin{tabular}{|c|c|c|c|c|}
\hline
$k$ & $\sum_{d\in D_+(X)} L(1/2,\chi_d)^k$ & $\frac{3}{\pi^2} \int_1^X Q_+(k,\log{t}) dt$& $R_+(k,x)$ & $\Delta_+(k,X)$  \\ \hline
1 & 4074391863.4447 & 4074392042.9388 & .99999995594580 & -179.4941 \\ 
1 & 8445624718.0243 & 8445624023.3138 & 1.0000000822569 & 694.7105 \\ 
1 & 12928896894.590 & 12928896383.146 & 1.0000000395582 & 511.444 \\ 
1 & 17484928279.579 & 17484927921.500 & 1.0000000204793 & 358.079 \\ 
1 & 22095062063.114 & 22095062690.738 & .99999997159438 & -627.624 \\ 
\hline
2 & 76310075816.466 & 76310057832.320 & 1.0000002356720 & 17984.146 \\ 
2 & 168051689378.93 & 168051603484.03 & 1.0000005111222 & 85894.90 \\ 
2 & 266303938917.29 & 266303916920.62 & 1.0000000825999 & 21996.67 \\ 
2 & 368948427173.22 & 368948308826.37 & 1.0000003207681 & 118346.85 \\ 
2 & 474942139636.16 & 474942177549.68 & .99999992017235 & -37913.52 \\ 
\hline
3 & 2478393690176.2 & 2478391641054.5 & 1.0000008267950 & 2049121.7 \\ 
3 & 5878735240405.9 & 5878729153410.4 & 1.0000010354271 & 6086995.5 \\ 
3 & 9720154390088.4 & 9720158187579.5 & .99999960931797 & -3797491.1 \\ 
3 & 13873264940982. & 13873252832529. & 1.0000008727912 & 12108453. \\ 
3 & 18271480140004. & 18271496263135. & .99999911758015 & -16123131. \\ 
\hline
4 & .10868425484737e15 & .10868409751016e15 & 1.0000014476562 & 157337203 \\ 
4 & .27974980520169e15 & .27974915668497e15 & 1.0000023182079 & 648516719 \\ 
4 & .48473276073219e15 & .48473329605692e15 & .99999889563038 & -535324726 \\ 
4 & .71493167429315e15 & .71492961664246e15 & 1.0000028781164 & 2057650687 \\ 
4 & .96564046289913e15 & .96564334647659e15 & .99999701382764 & -2883577466 \\ 
\hline
5 & .57022430562904e16 & .57022322406897e16 & 1.0000018967310 & 10815600670 \\ 
5 & .15999737676260e17 & .15999653478756e17 & 1.0000052624580 & .84197504e11 \\ 
5 & .29130430291967e17 & .29130495012249e17 & .99999777826357 & -.64720282e11 \\ 
5 & .44482716417300e17 & .44482376920928e17 & 1.0000076321545 & .339496372e12 \\ 
5 & .61707290890367e17 & .61707708869778e17 & .99999322646362 & -.417979411e12 \\ 
\hline
6 & .33658290814098e18 & .33658163201404e18 & 1.0000037914337 & .127612694e13 \\ 
6 & .10326933113376e19 & .10326816848898e19 & 1.0000112585010 & .116264478e14 \\ 
6 & .19792425806612e19 & .19792515491256e19 & .99999546875969 & -.89684644e13 \\ 
6 & .31332379844474e19 & .31331890401641e19 & 1.0000156212353 & .489442833e14 \\ 
6 & .44685941512069e19 & .44686487402482e19 & .99998778399367 & -.545890413e14 \\ 
\hline
7 & .215991539086e20 & .215989246213e20 & 1.00001061568 & .2292873e15 \\ 
7 & .726312167992e20 & .726295668031e20 & 1.00002271797 & .16499961e16 \\ 
7 & .146733199900e21 & .146734533114e21 & .999990914109 & -.1333214e16 \\ 
7 & .241042340834e21 & .241036160843e21 & 1.00002563927 & .6179991e16 \\ 
7 & .353694078736e21 & .353700808054e21 & .999980974547 & -.6729318e16 \\ 
\hline
8 & .1475899774e22 & .1475859642e22 & 1.000027192 & .40132e17 \\ 
8 & .5449090667e22 & .5448853612e22 & 1.000043505 & .237055e18 \\ 
8 & .1161602962e23 & .1161622793e23 & .9999829282 & -.19831e18 \\ 
8 & .1981618159e23 & .1981550523e23 & 1.000034133 & .67636e18 \\ 
8 & .2993403001e23 & .2993484649e23 & .9999727248 & -.81648e18 \\ 
\hline
\end{tabular}
}
\caption[Moments of $L(1/2,\ch_d)$ versus conjectured asymptotics,
$d<0$]{Moments $\sum_{d\in D_+(X)} L(1/2,\chi_d)^k$ versus $\frac{3}{\pi^2}
\int_1^X Q_+(k,\log{t}) dt$, for $k=1,\ldots,8$ and $d>0$. Five values for
each $k$ are shown, $X = 2\times 10^9, 4\times 10^9, \ldots, 10^{10}$.}
\label{tab:Lhalfchip}
\end{table}

\renewcommand{\thefigure}{\arabic{figure}}

\newpage
\thispagestyle{empty}
\begin{figure}[H]
    \centerline{
       \includegraphics[width=.48\textwidth,height=2in]{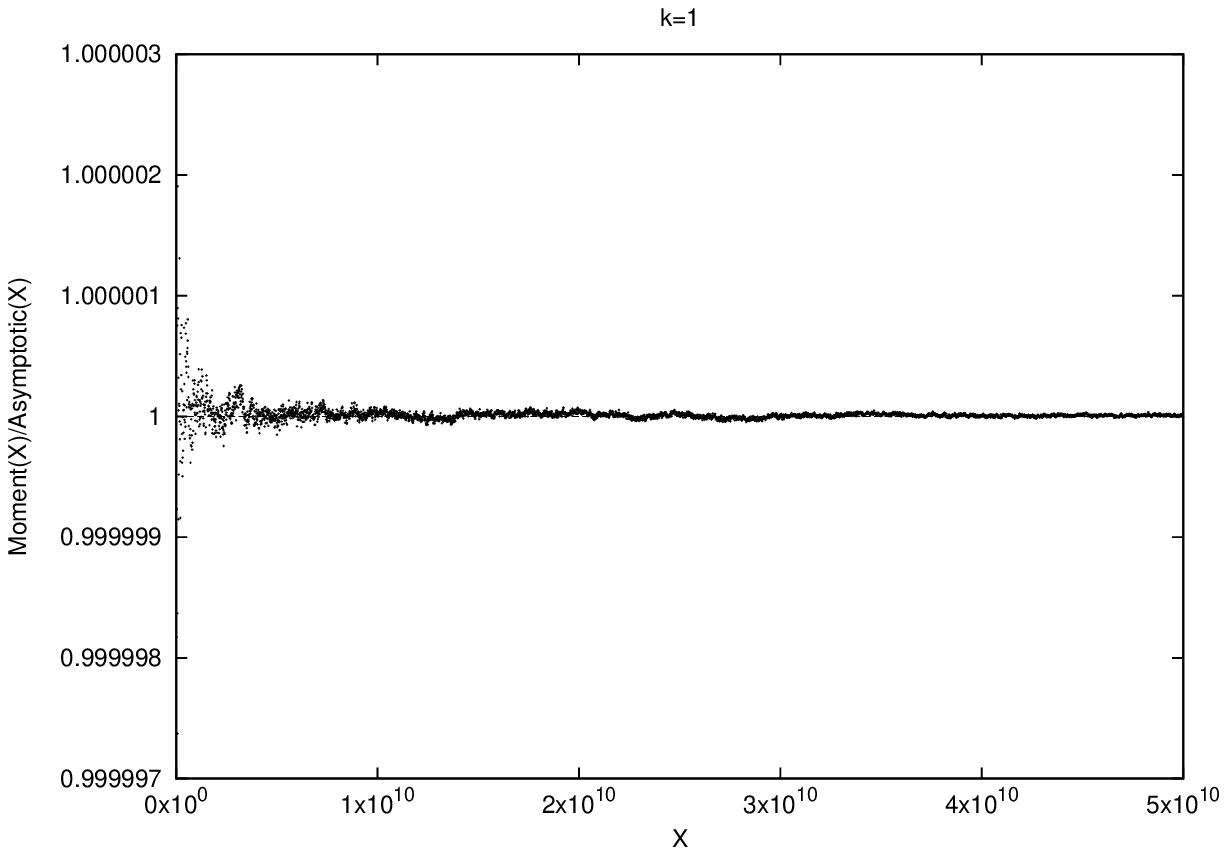}
       \includegraphics[width=.48\textwidth,height=2in]{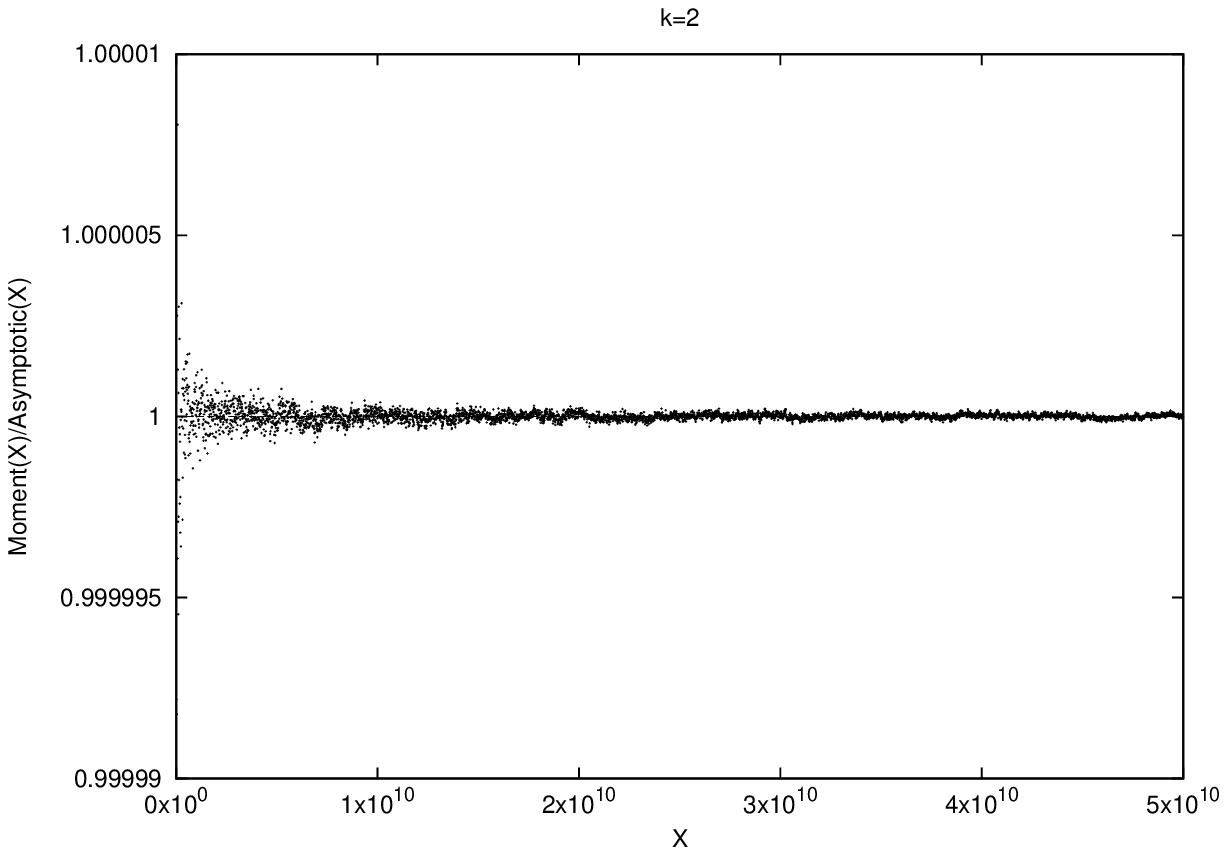}
     }
     \centerline{
        \includegraphics[width=.48\textwidth,height=2in]{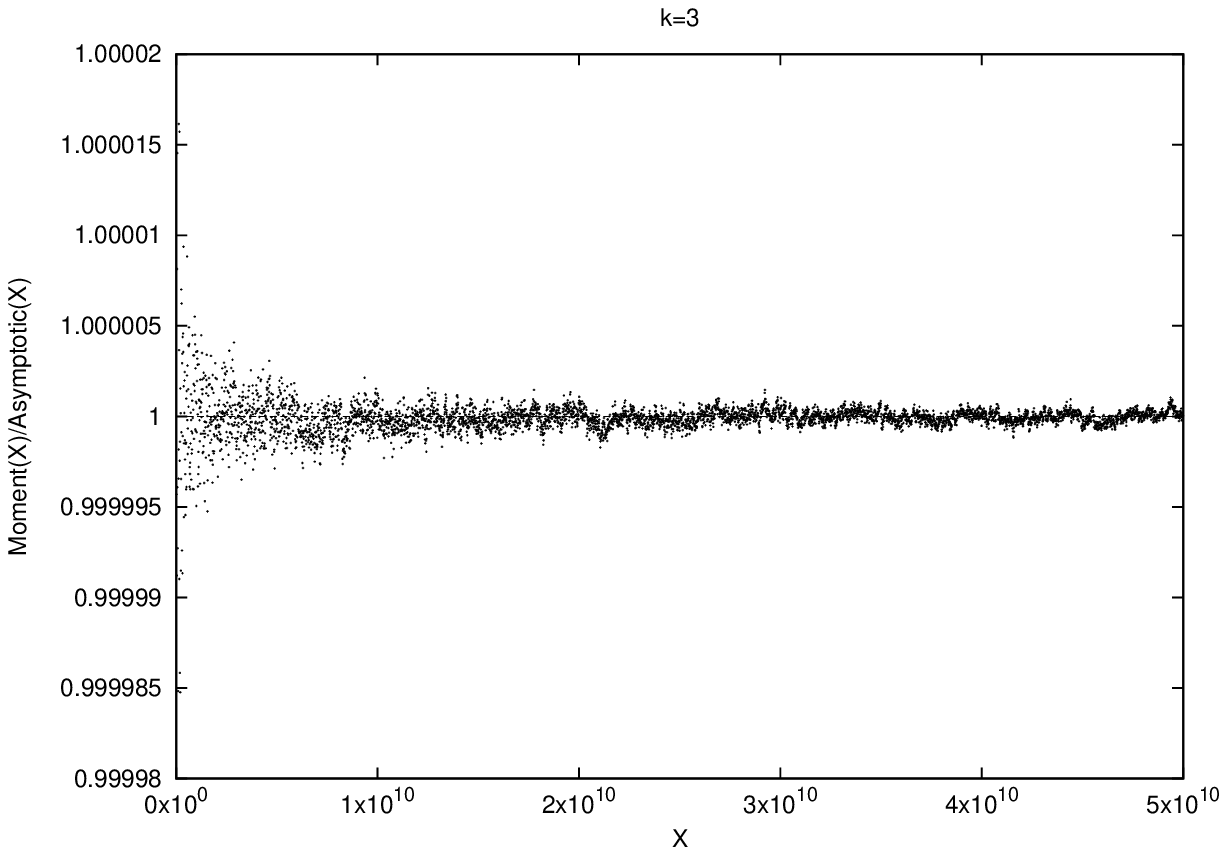}
        \includegraphics[width=.48\textwidth,height=2in]{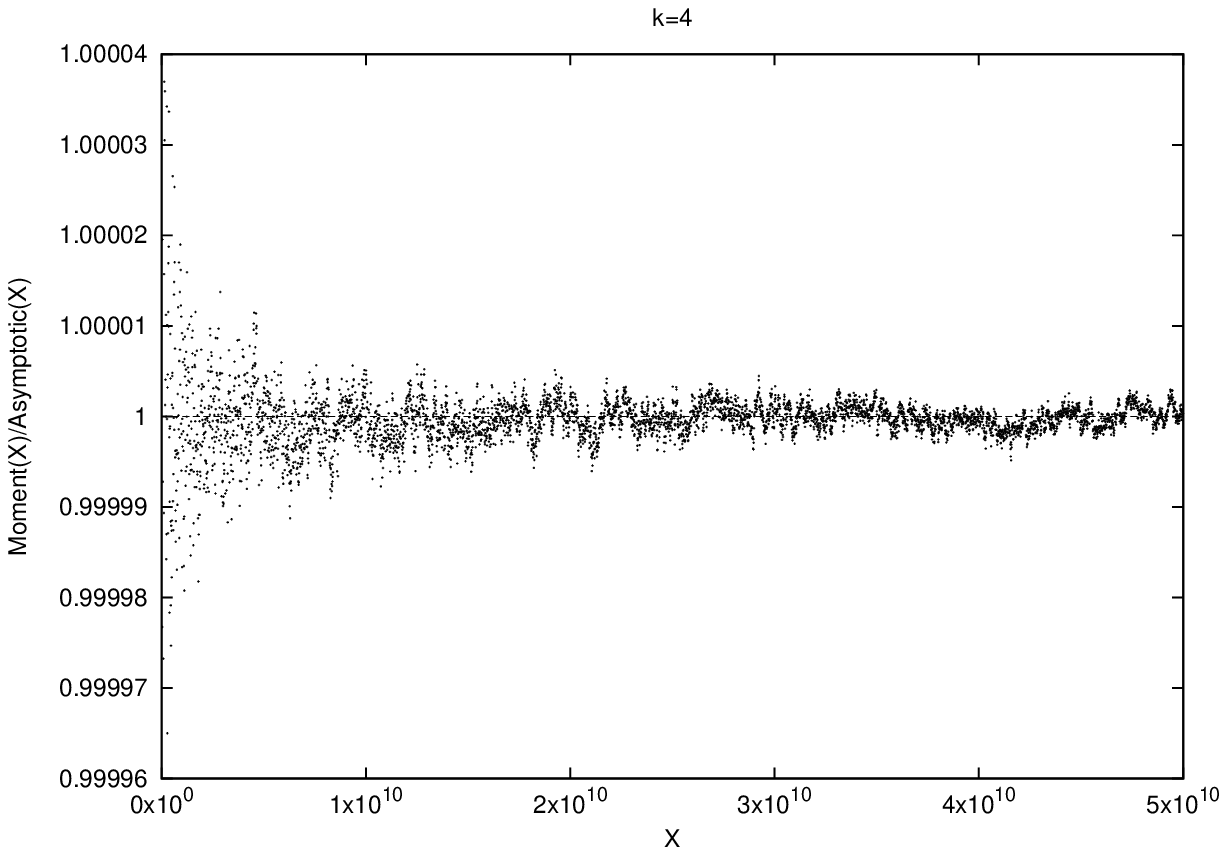}
      }
     \centerline{
       \includegraphics[width=.48\textwidth,height=2in]{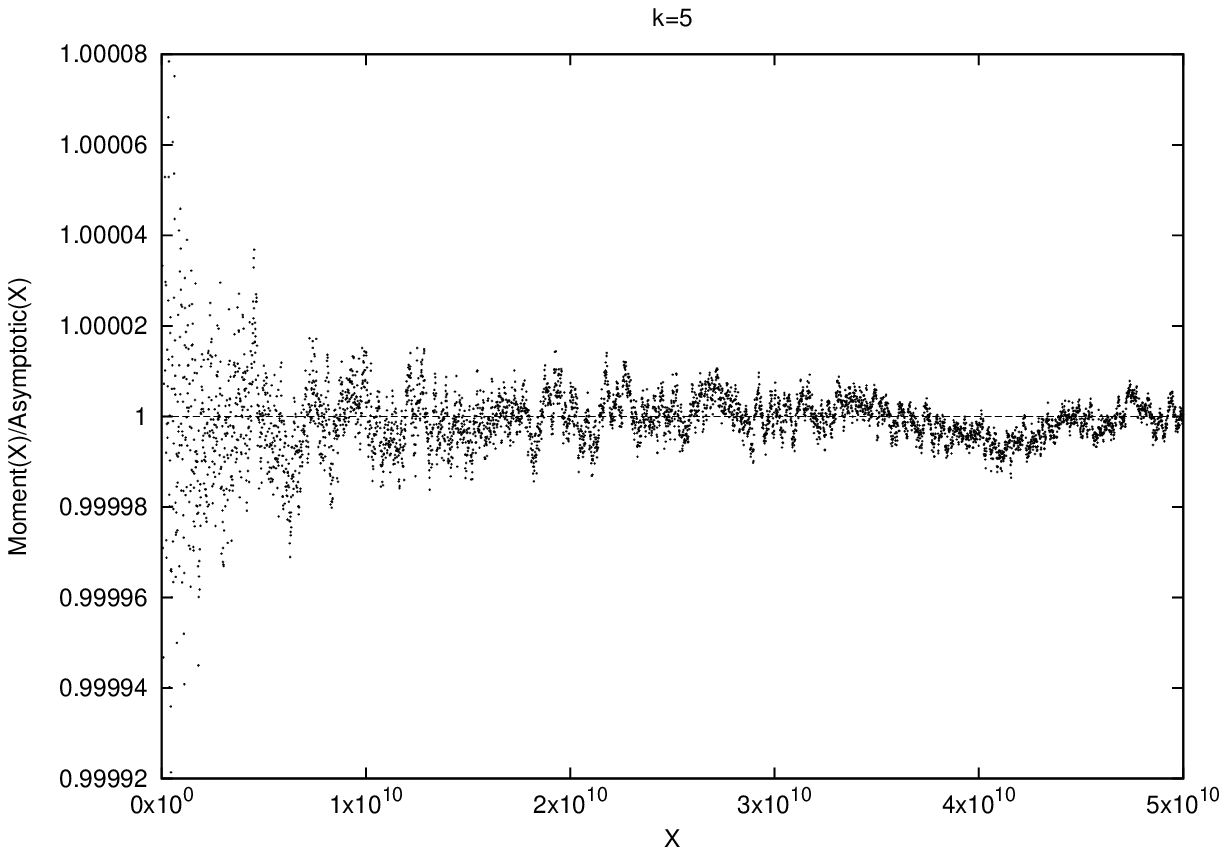}
        \includegraphics[width=.48\textwidth,height=2in]{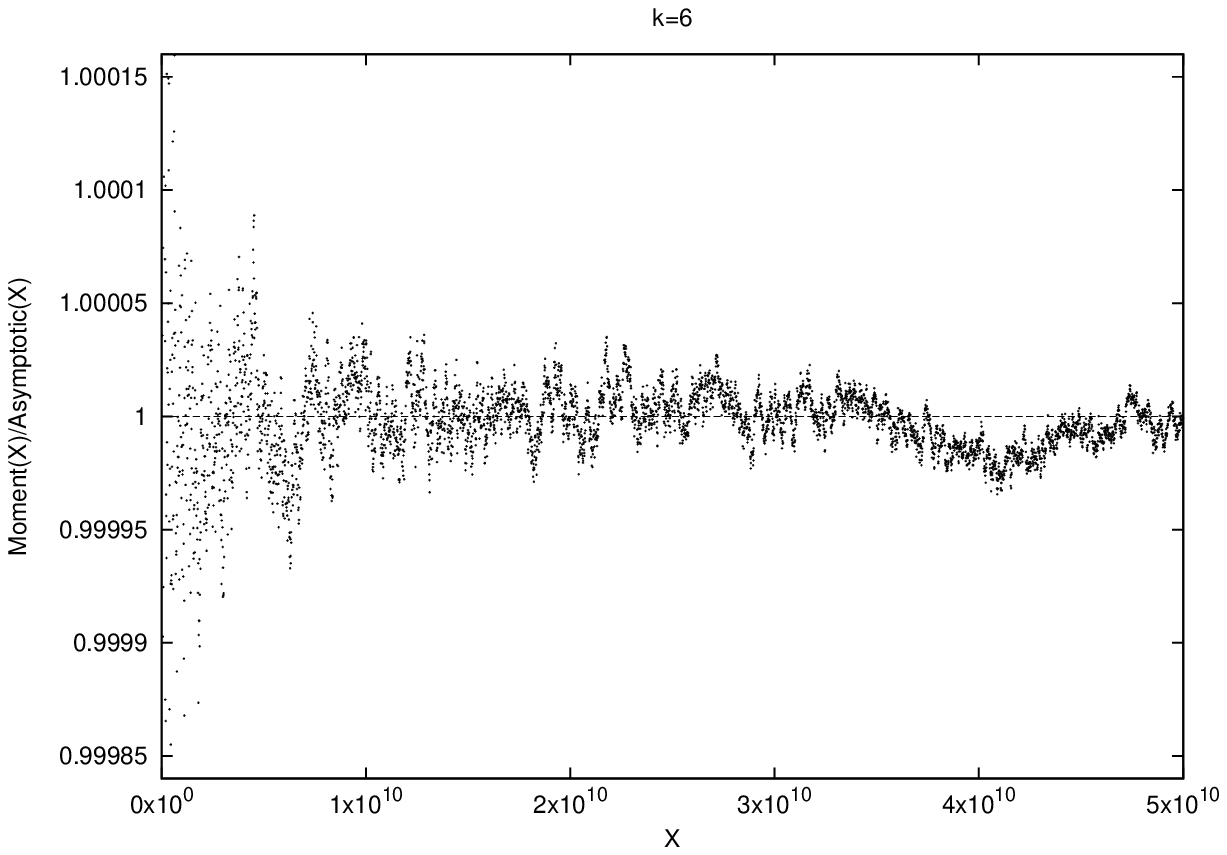}
     }  
     \centerline{
       \includegraphics[width=.48\textwidth,height=2in]{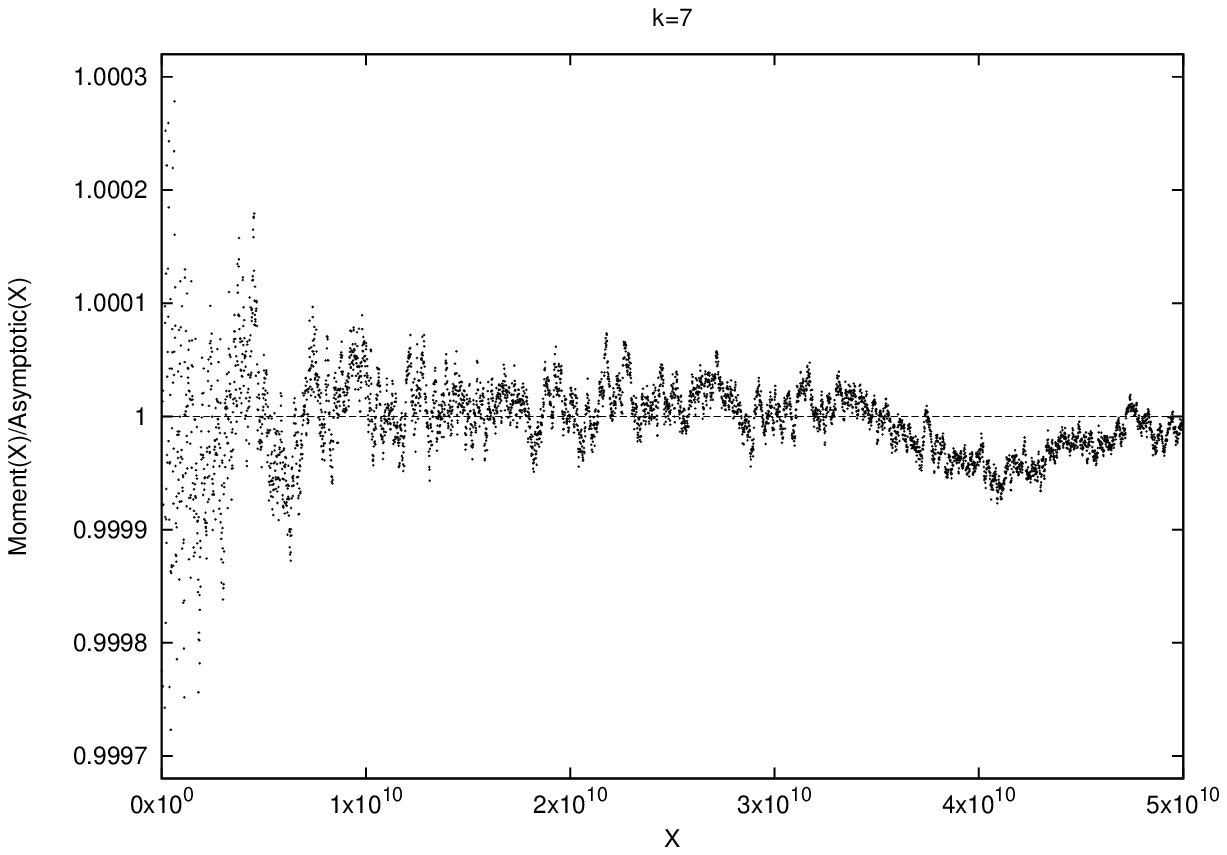}
       \includegraphics[width=.48\textwidth,height=2in]{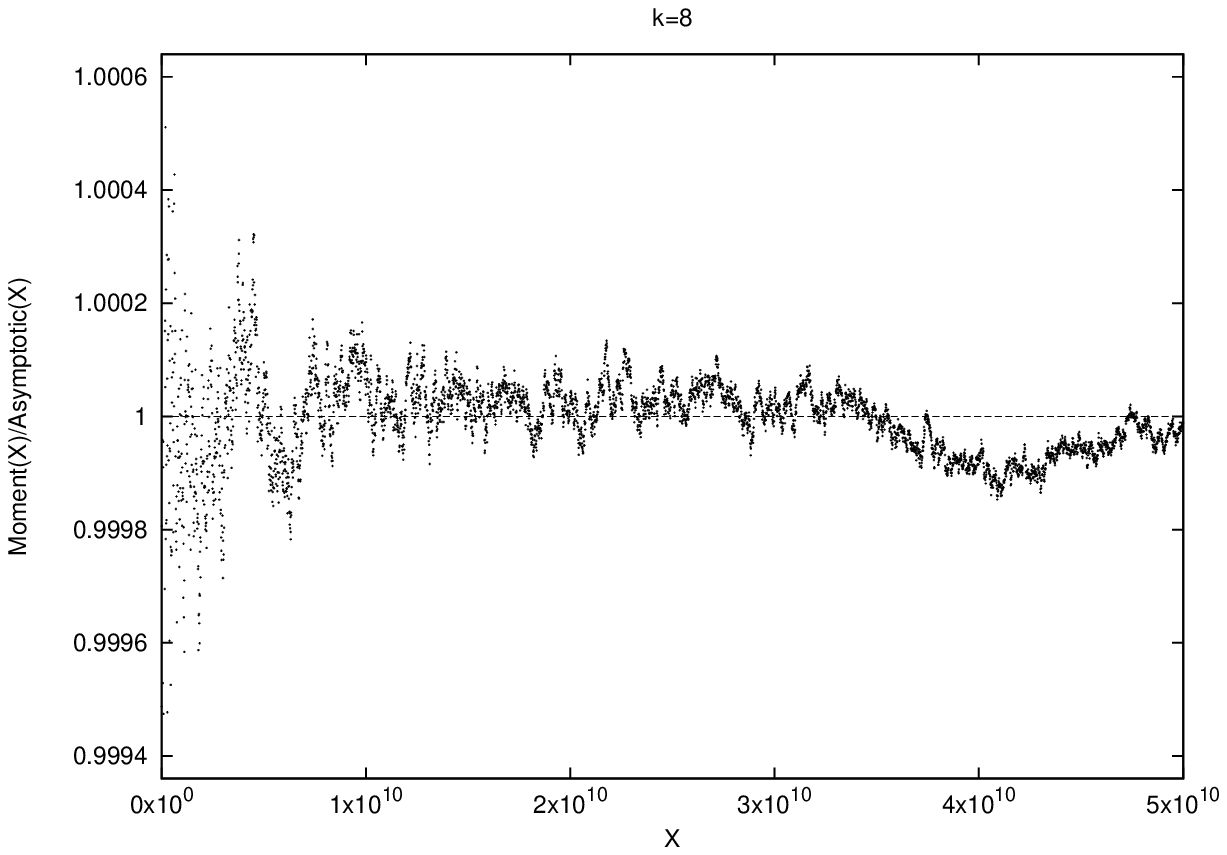}
     }
     \caption[Plot of the ratio $R_-(k,X)$for $k=1,\ldots, 8$]{These
     plots depict the ratio $R_-(k,X)$ of the numerically computed moments compared to
     the CFKRS predictions, for $k=1,\ldots, 8$ and $d<0$, sampled every $10^7$, i.e.
     at $X = 10^7, 2\times 10^7, \ldots, 5\times 10^{10}$. The horizontal axis is
     $X$, the vertical axis is the ratio $R_-(k,X)$.}\label{fig: fig1}
\end{figure}

\newpage
\thispagestyle{empty}
\begin{figure}[H]
 \centerline{
     \includegraphics[width=0.48\textwidth,height=2in]{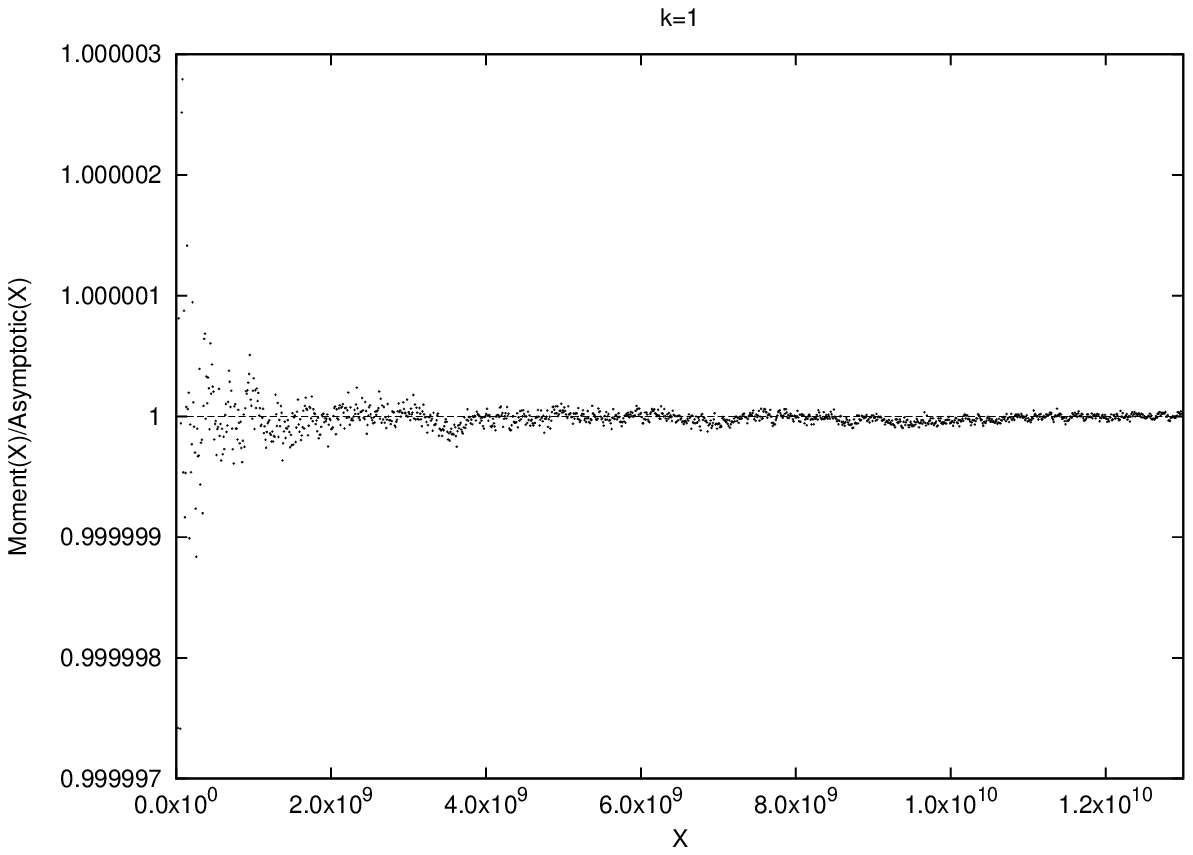}
     \includegraphics[width=0.48\textwidth,height=2in]{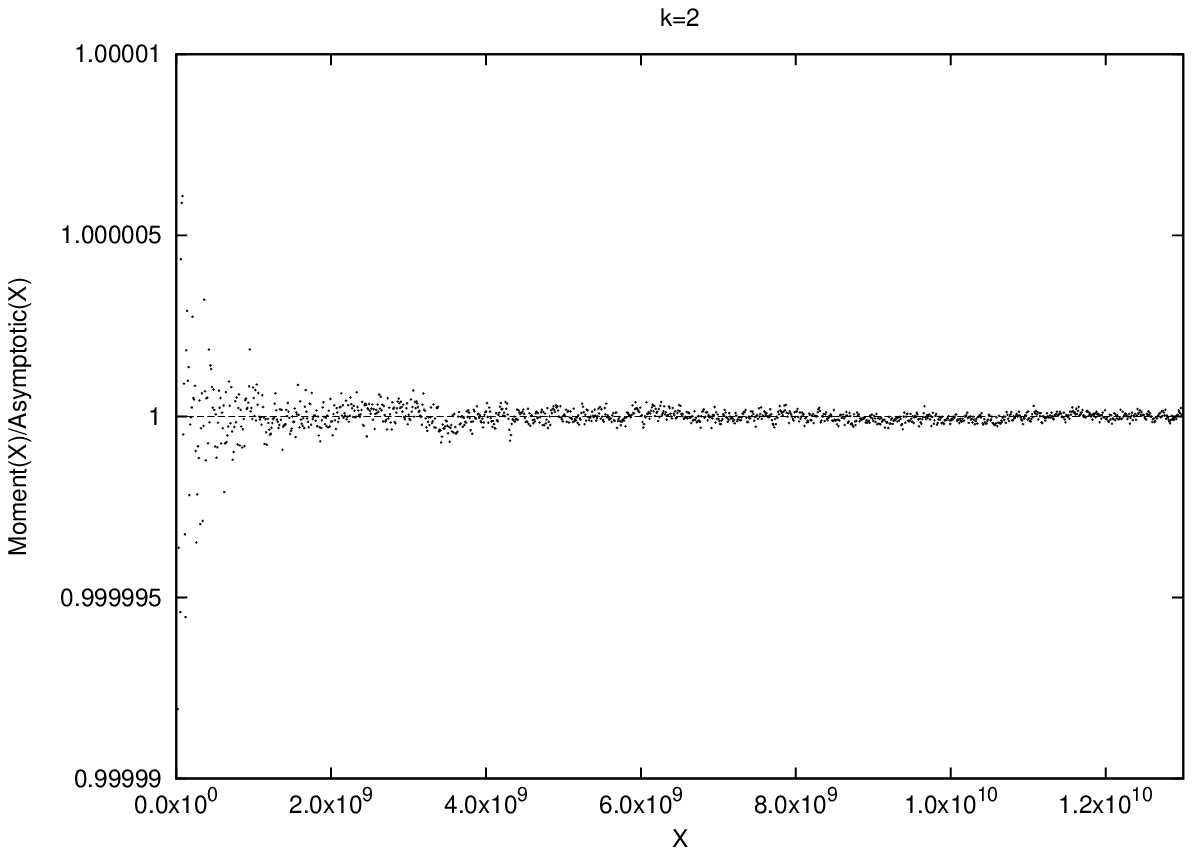}
   }
  \centerline{
     \includegraphics[width=0.48\textwidth,height=2in]{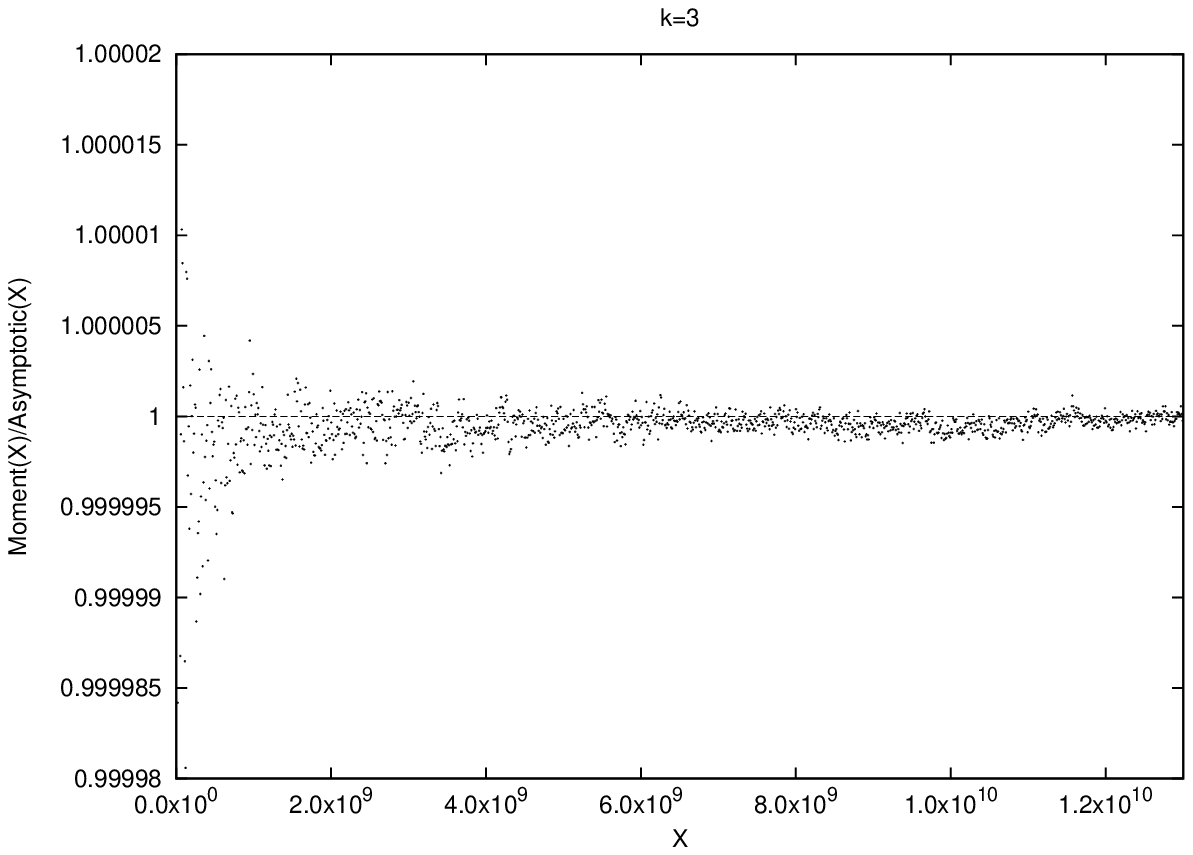}
     \includegraphics[width=0.48\textwidth,height=2in]{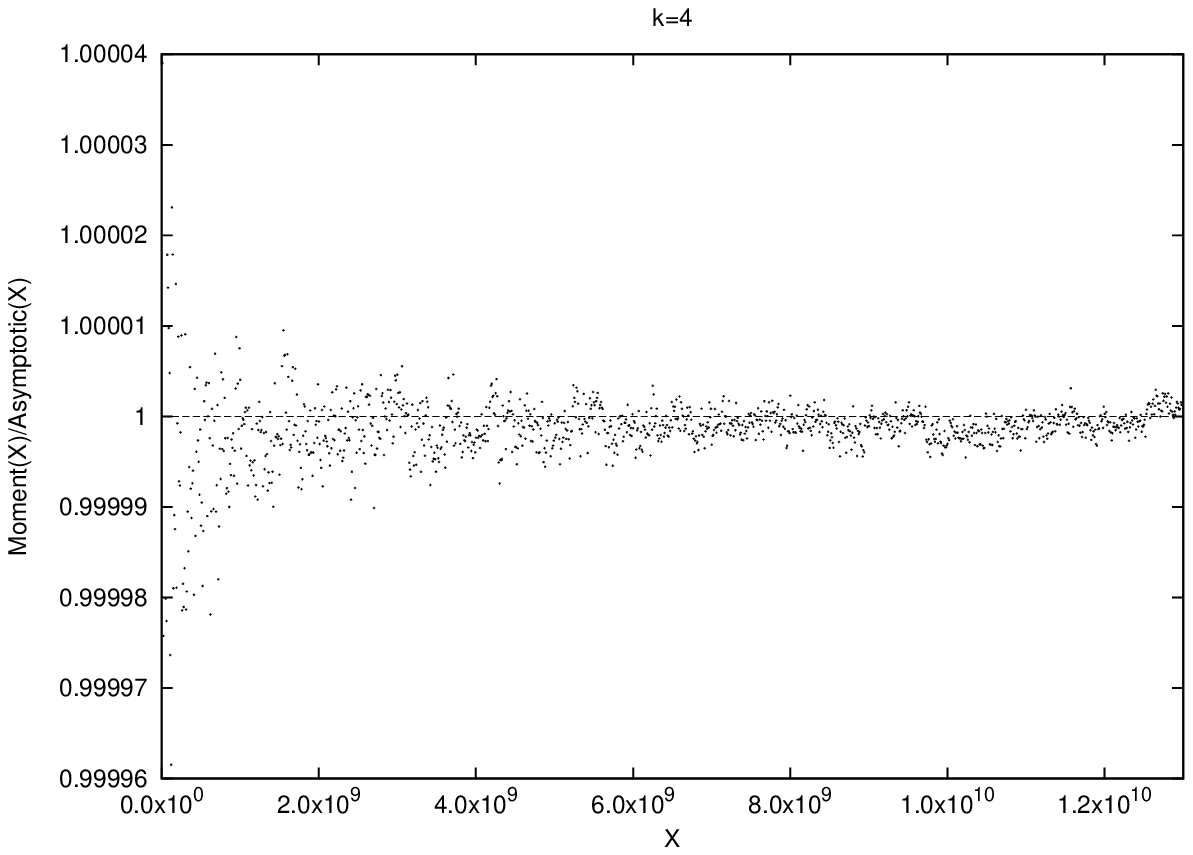}
   }
  \centerline{
    \includegraphics[width=0.48\textwidth,height=2in]{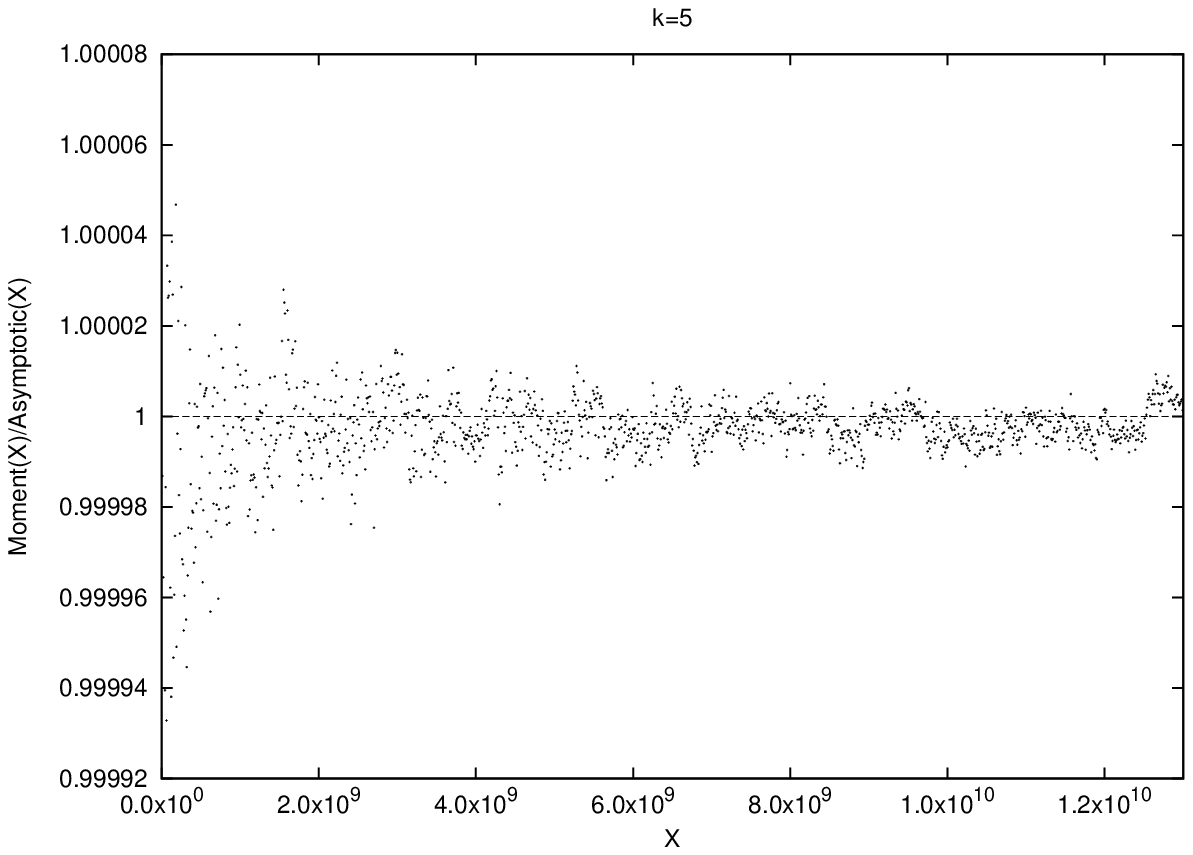}
     \includegraphics[width=0.48\textwidth,height=2in]{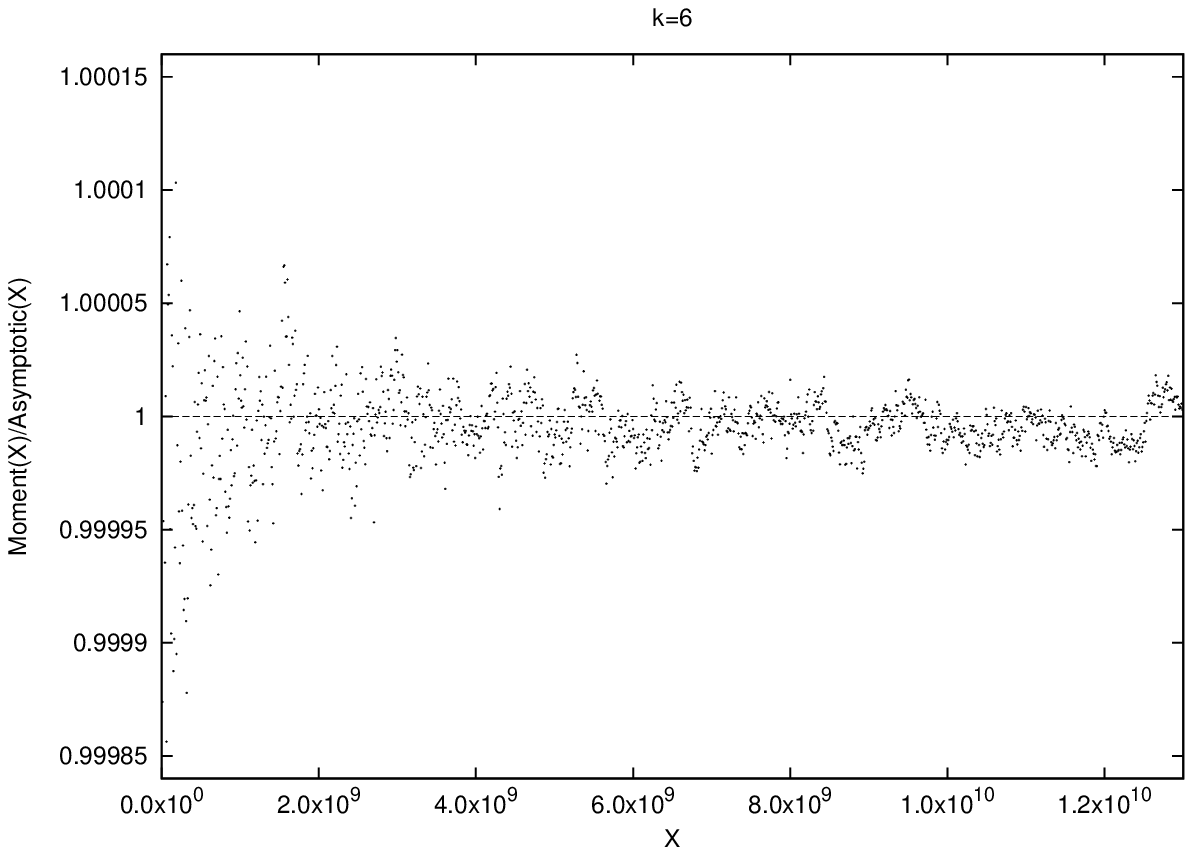}
  }
  \centerline{
    \includegraphics[width=0.48\textwidth,height=2in]{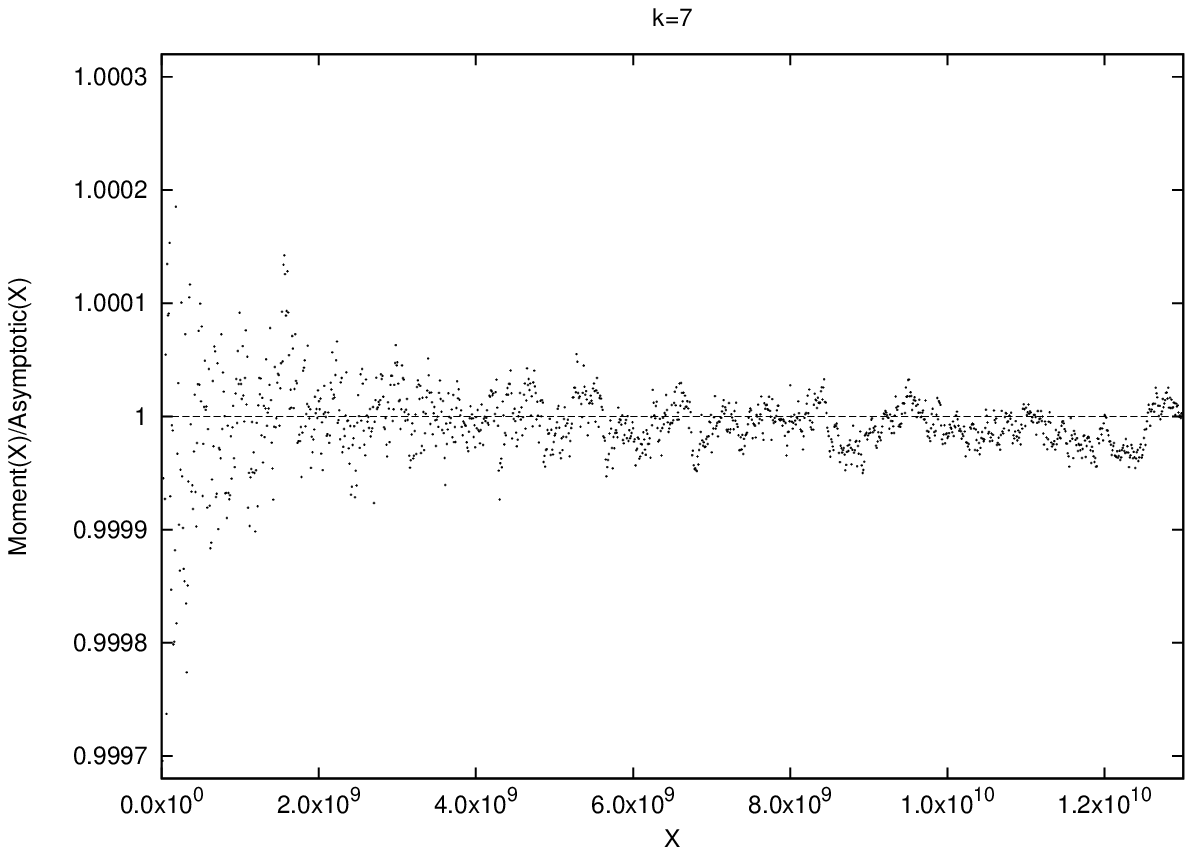}
     \includegraphics[width=0.48\textwidth,height=2in]{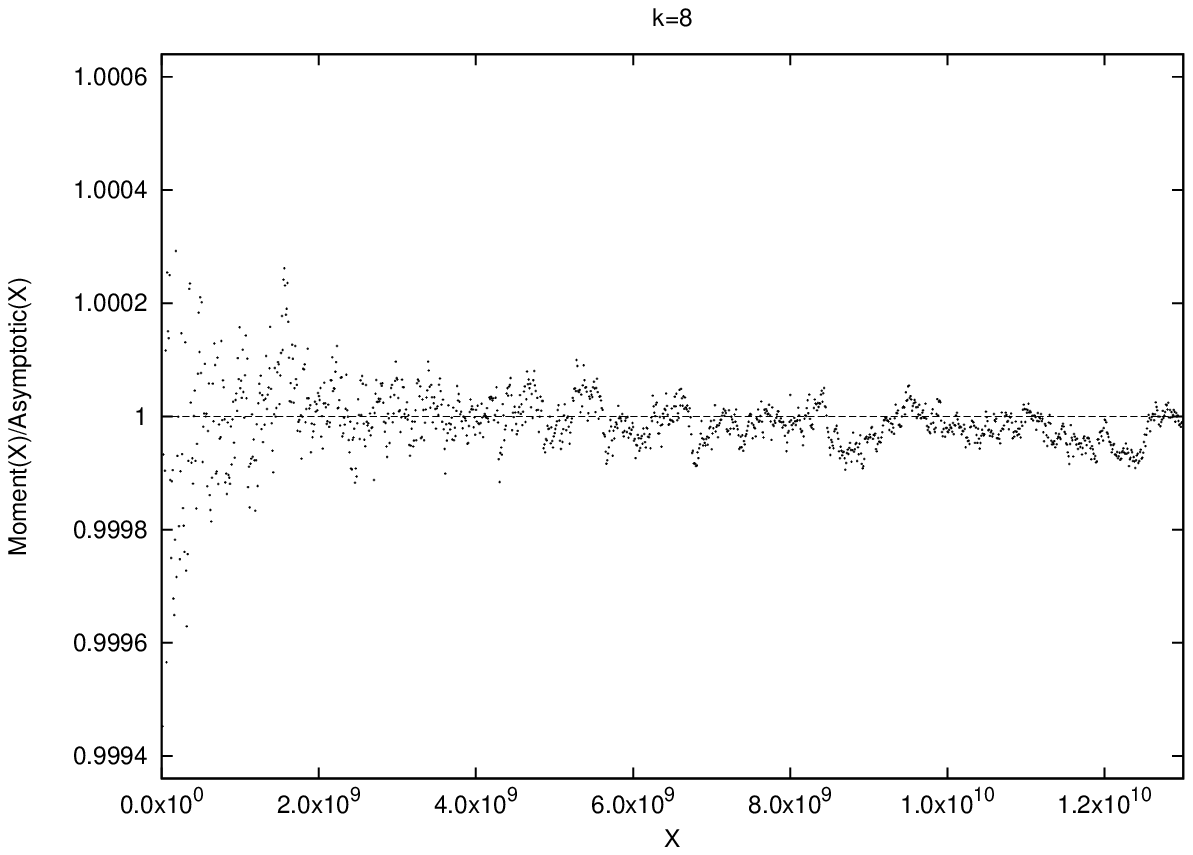}
  }
   \caption[Plot of the ratio $R_+(k,X)$ for $k=1,\ldots, 8$]{These
   plots depict the ratio $R_+(k,X)$ of the numerically computed moments compared to
   the CFKRS predictions, for $k=1,\ldots, 8$ and $d>0$, sampled every $10^7$, i.e.
   at $X = 10^7, 2\times 10^7, \ldots, 1.3\times 10^{10}$. The horizontal axis is
   $X$, the vertical axis is the ratio $R_+(k,X)$.}\label{fig: fig2}
\end{figure}

\newpage
\thispagestyle{empty}
\begin{figure}[H]
 \centerline{
     \includegraphics[width=0.48\textwidth,height=2in]{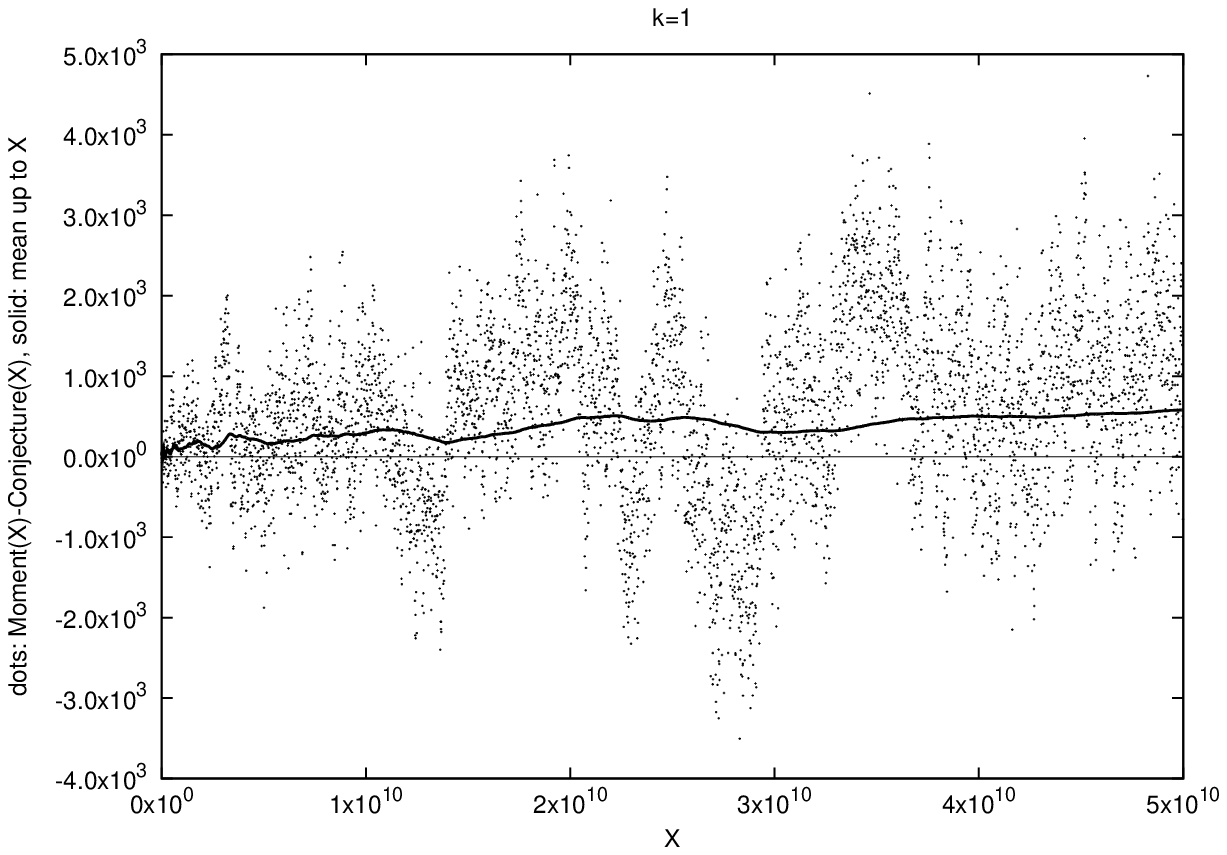}
     \includegraphics[width=0.48\textwidth,height=2in]{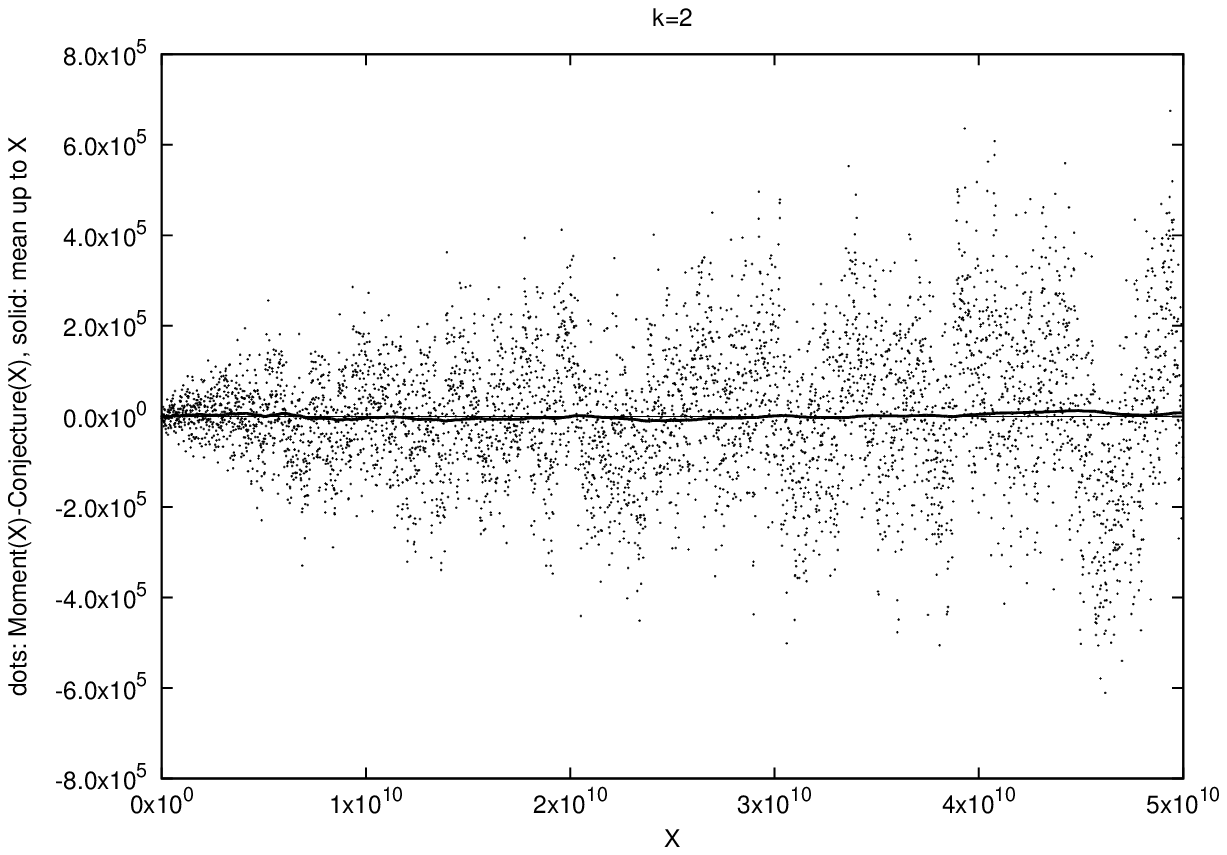}
   }
  \centerline{
    \includegraphics[width=0.48\textwidth,height=2in]{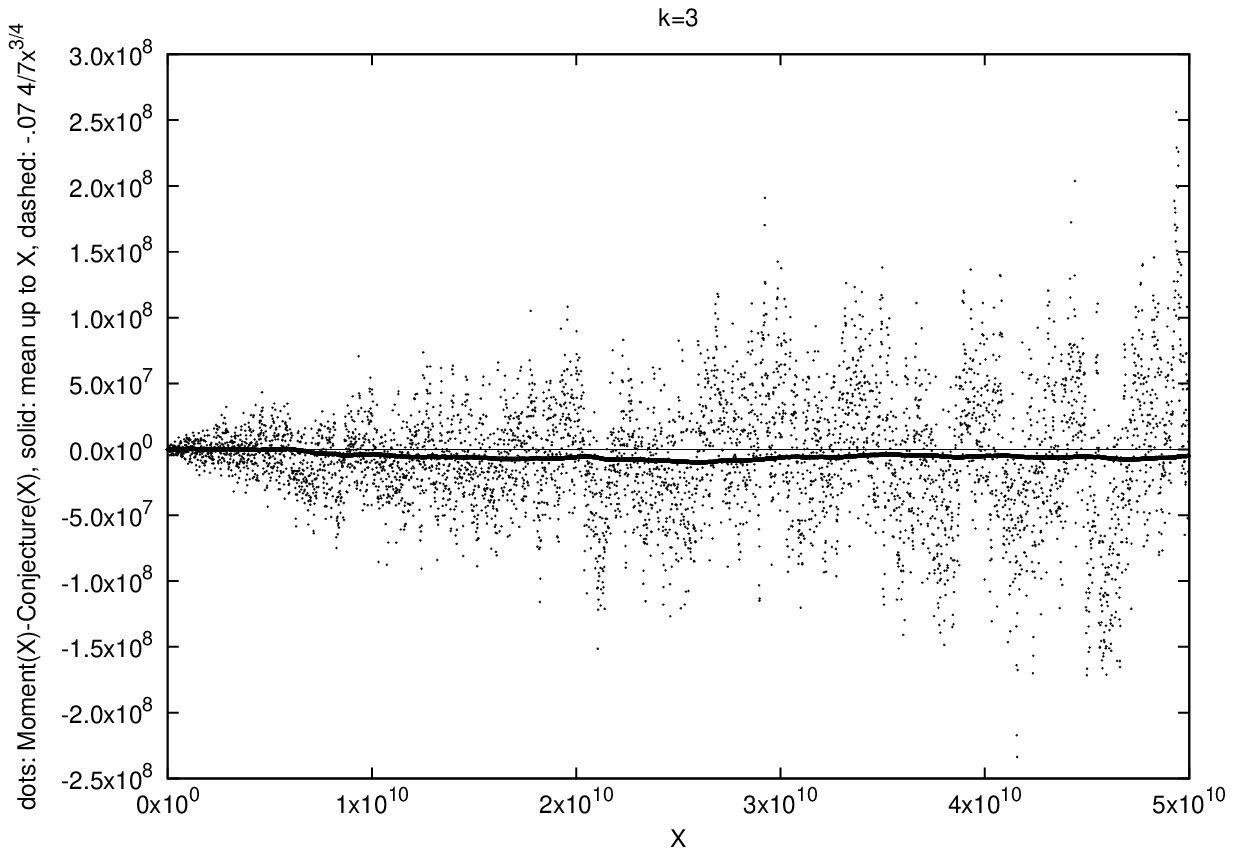}
     \includegraphics[width=0.48\textwidth,height=2in]{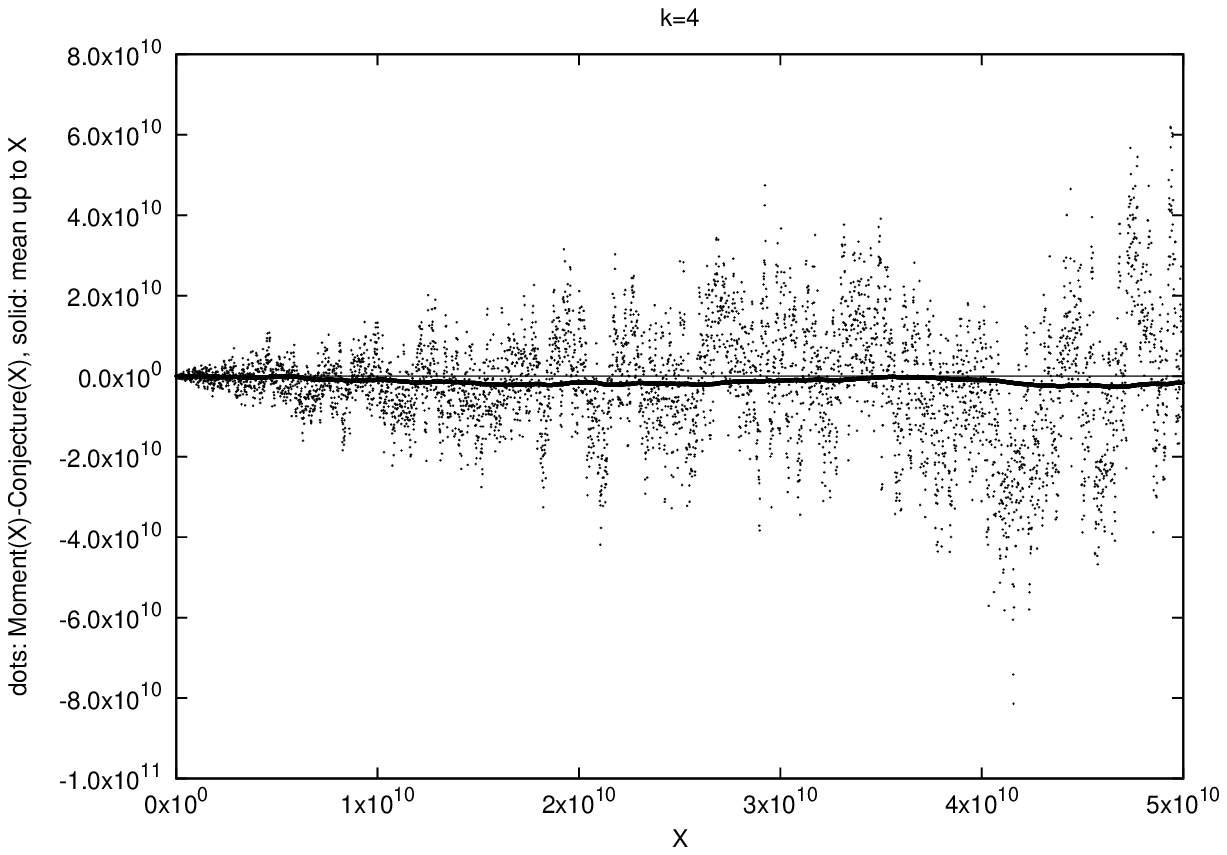}
  }
  \centerline{
    \includegraphics[width=0.48\textwidth,height=2in]{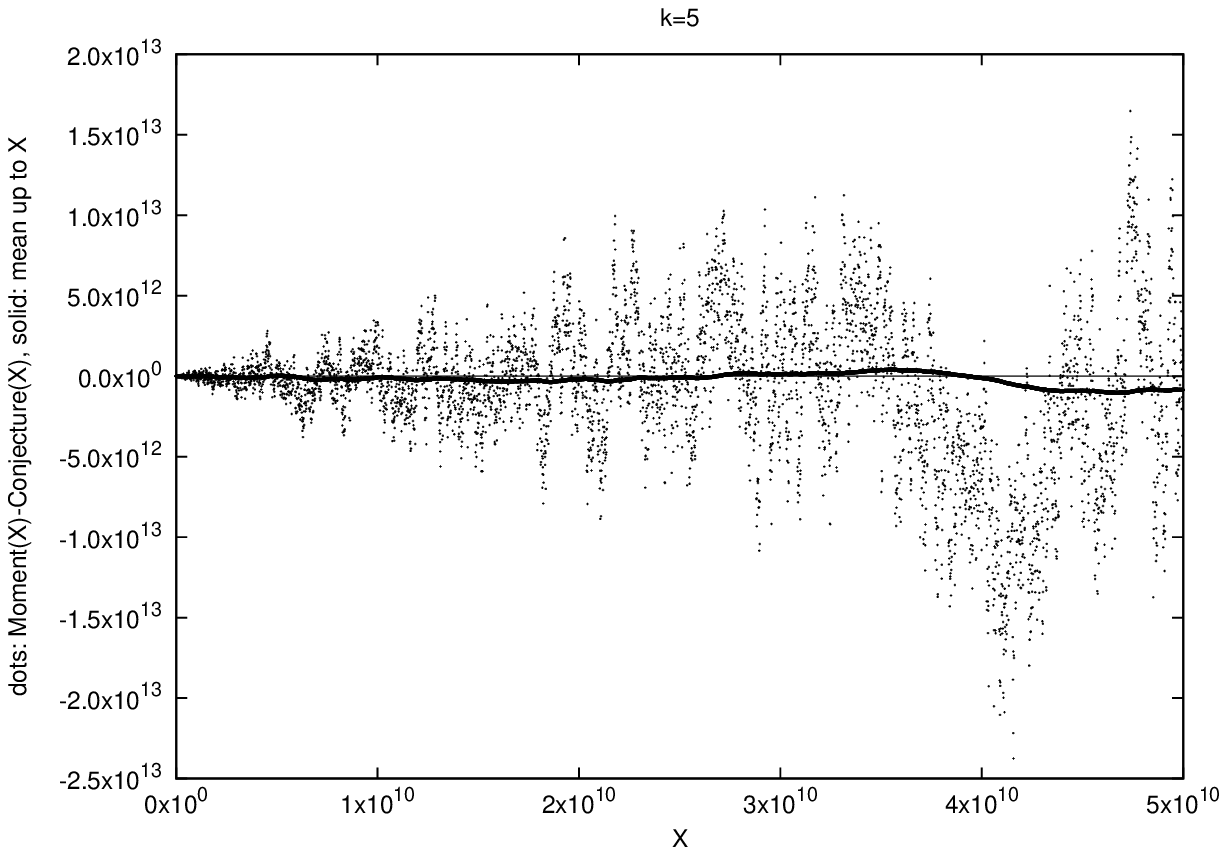}
     \includegraphics[width=0.48\textwidth,height=2in]{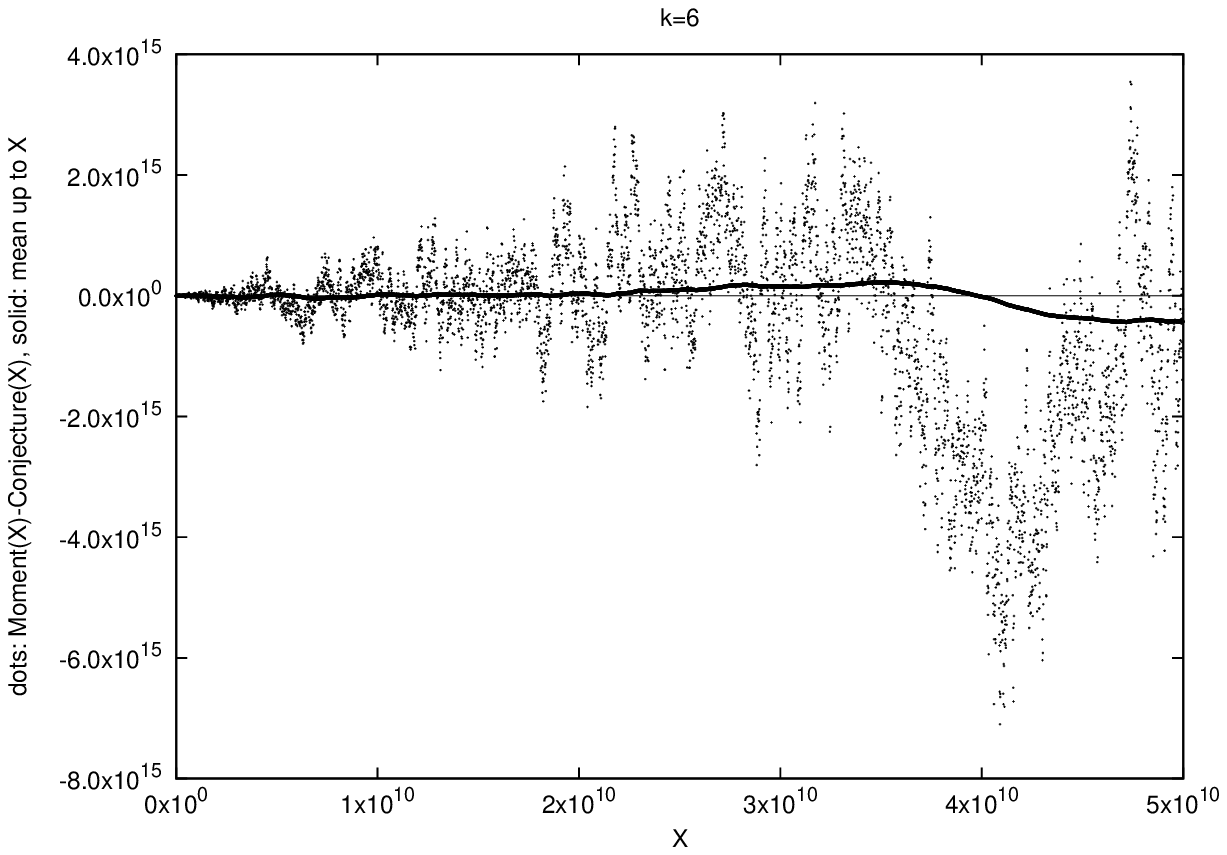}
  }
  \centerline{
    \includegraphics[width=0.48\textwidth,height=2in]{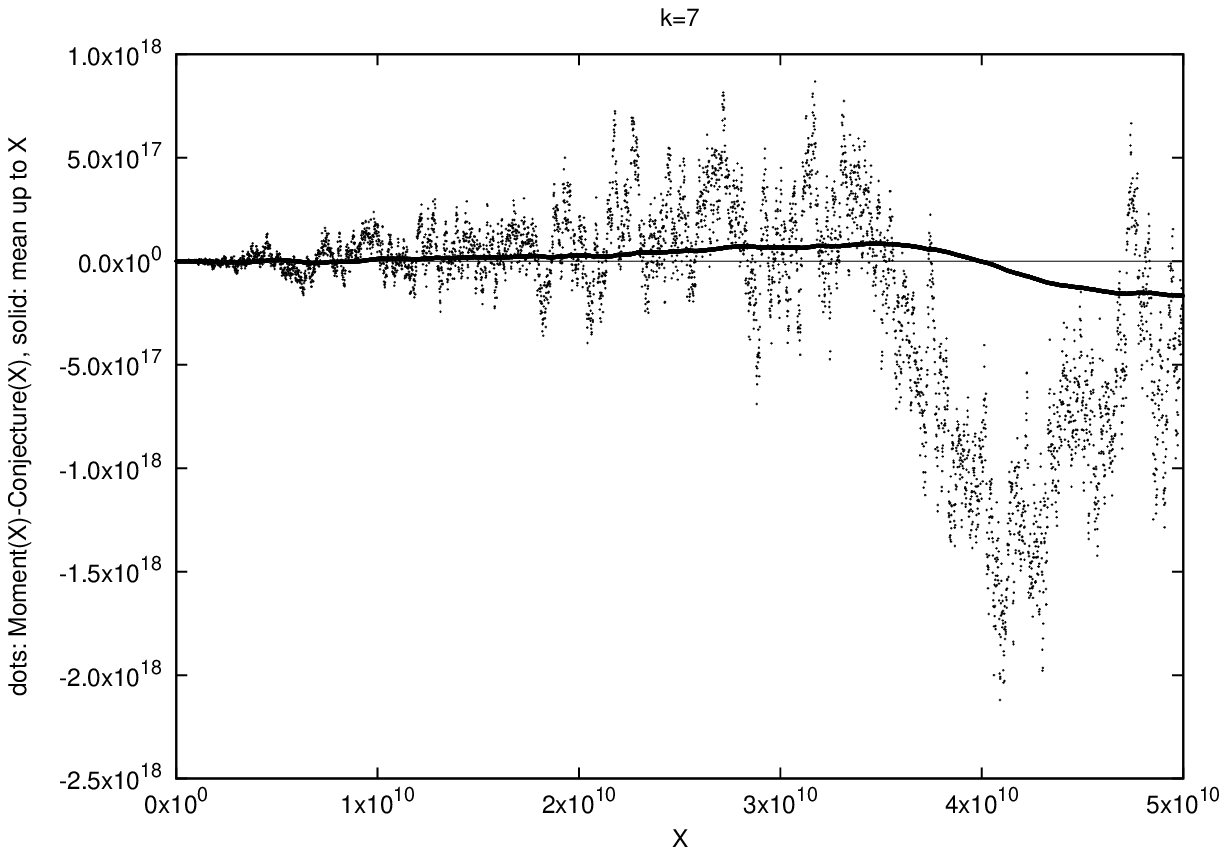}
     \includegraphics[width=0.48\textwidth,height=2in]{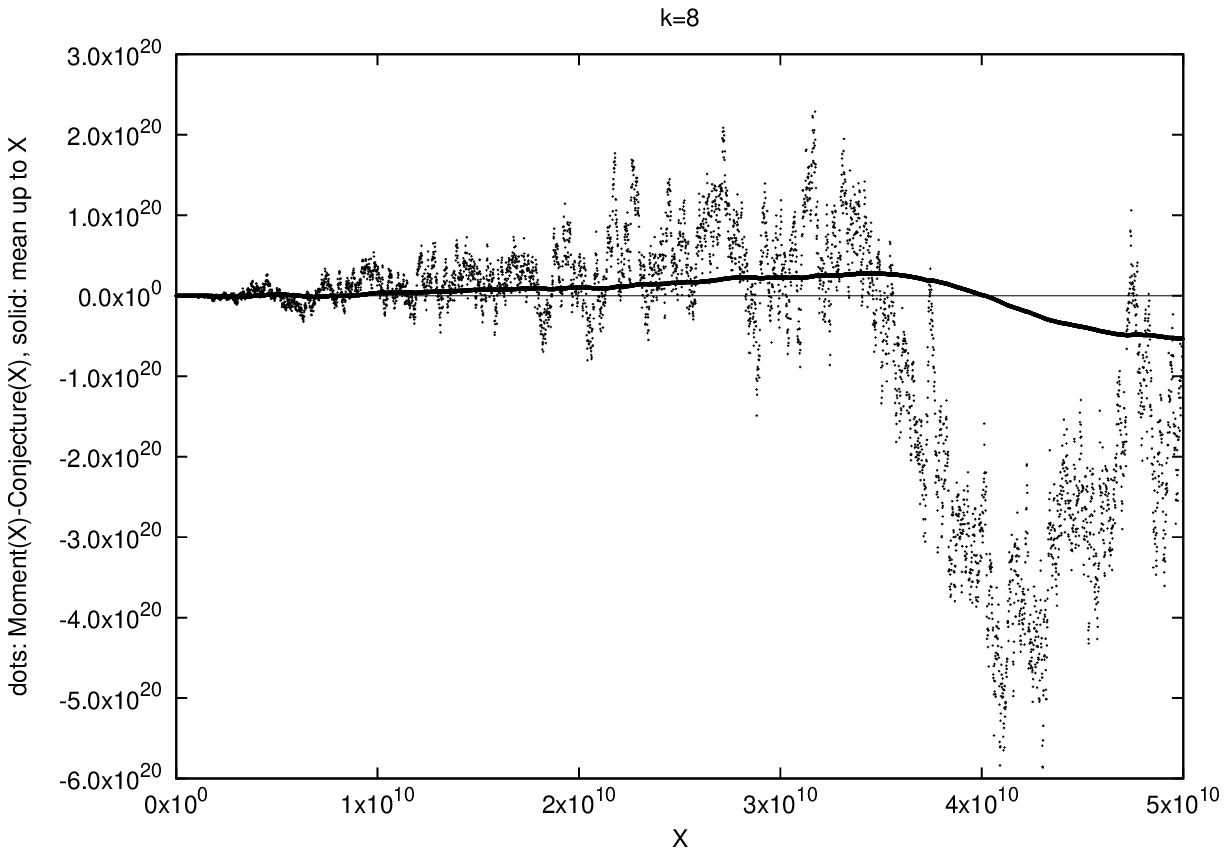}
  }
  \caption[Plot of the difference $\Delta_-(k,X)$ for $k=1,\ldots, 8$]{These
  plots depict the difference $\Delta_-(k,X)$ between the numerically computed moments and
  the CFKRS prediction, for $k=1,\ldots, 8$ and $d<0$,
  sampled at $X =  10^7, 2\times 10^7, \ldots, 5\times 10^{10}$. The horizontal axis
  is $X$, the vertical axis is the difference $\Delta_-(k,X)$, and the solid curve
  is the mean up to $X$ of the plotted differences (see the discussion above).}\label{fig: fig3}
\end{figure}

\newpage
\thispagestyle{empty}
\begin{figure}[H]
  \centerline{
     \includegraphics[width=0.48\textwidth,height=2in]{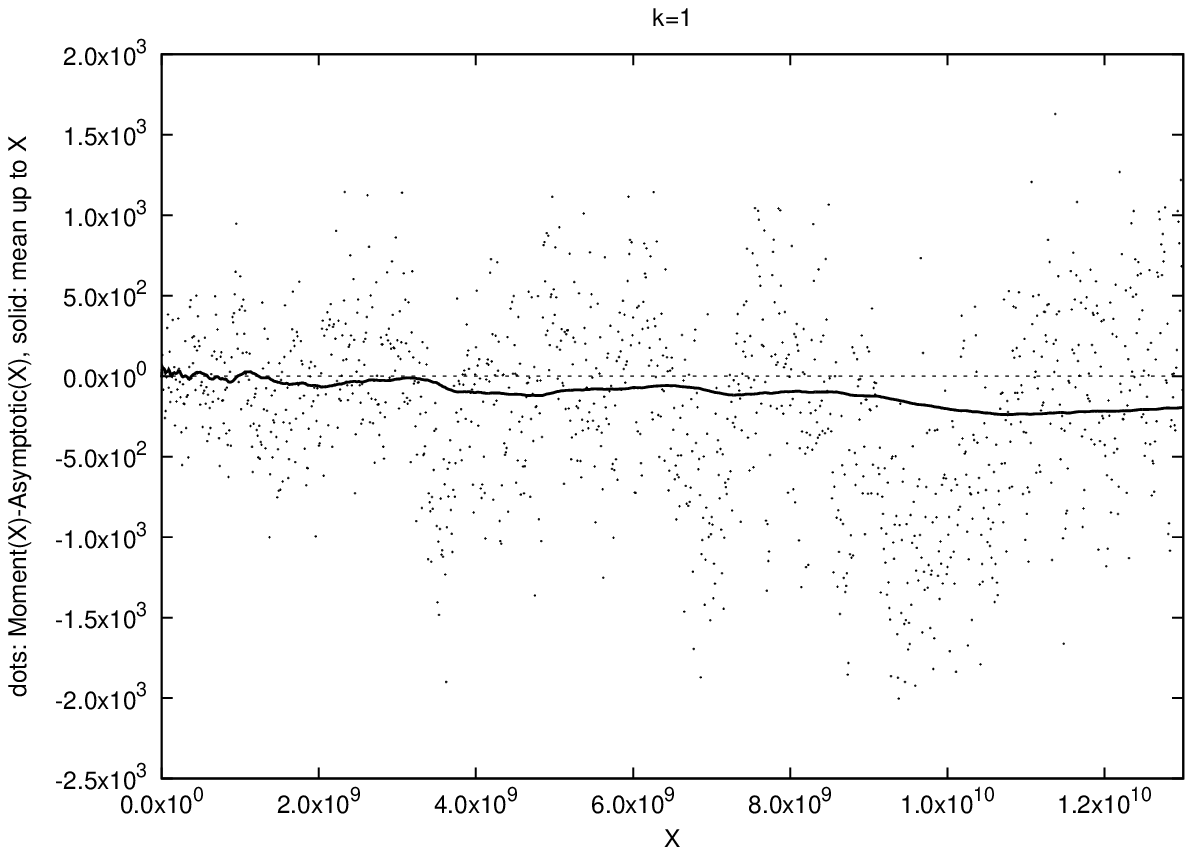}
     \includegraphics[width=0.48\textwidth,height=2in]{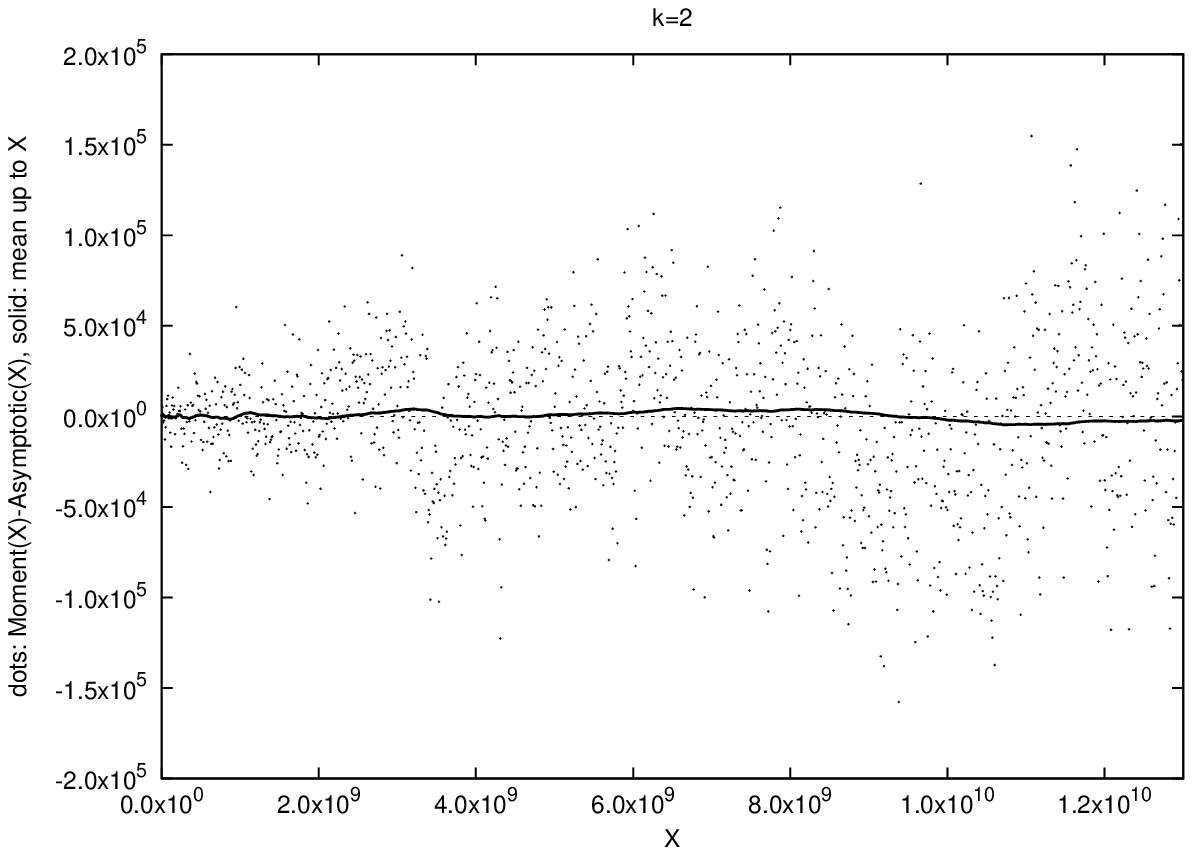}
   }
   \centerline{
     \includegraphics[width=0.48\textwidth,height=2in]{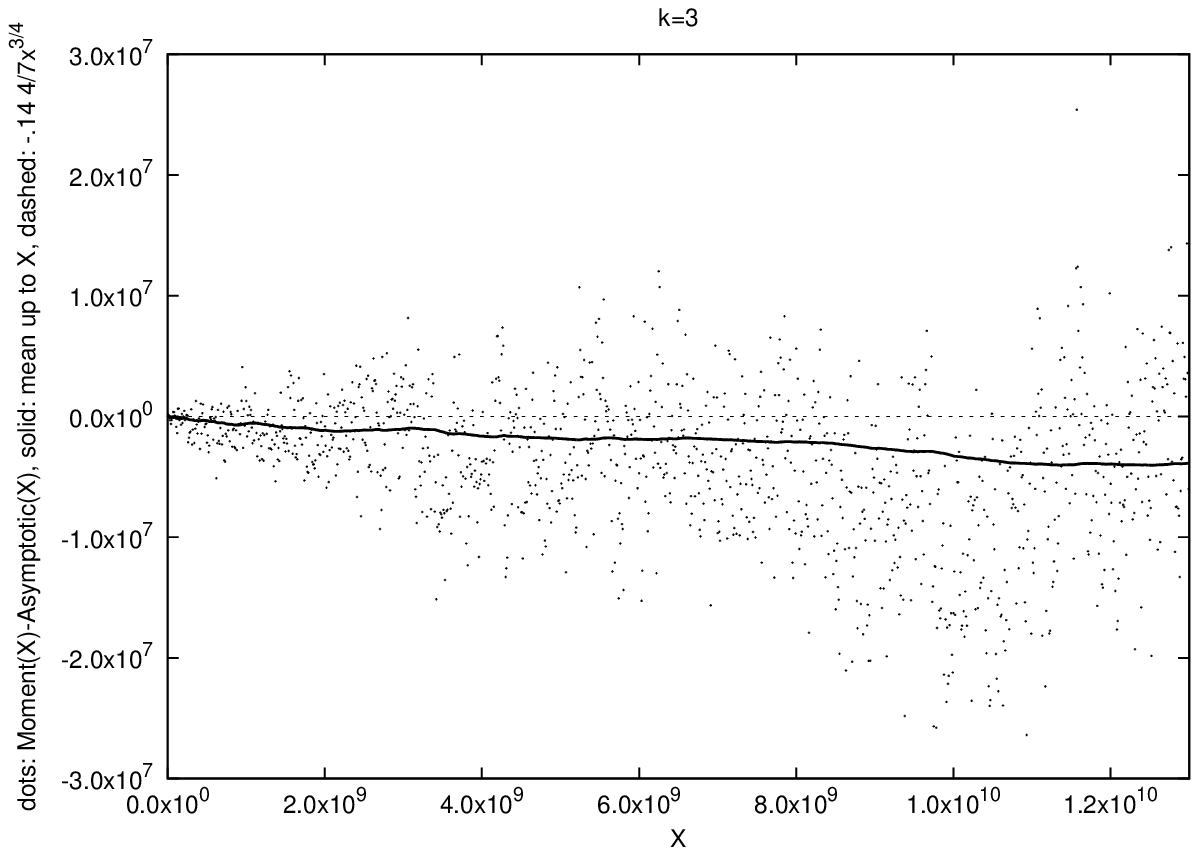}
     \includegraphics[width=0.48\textwidth,height=2in]{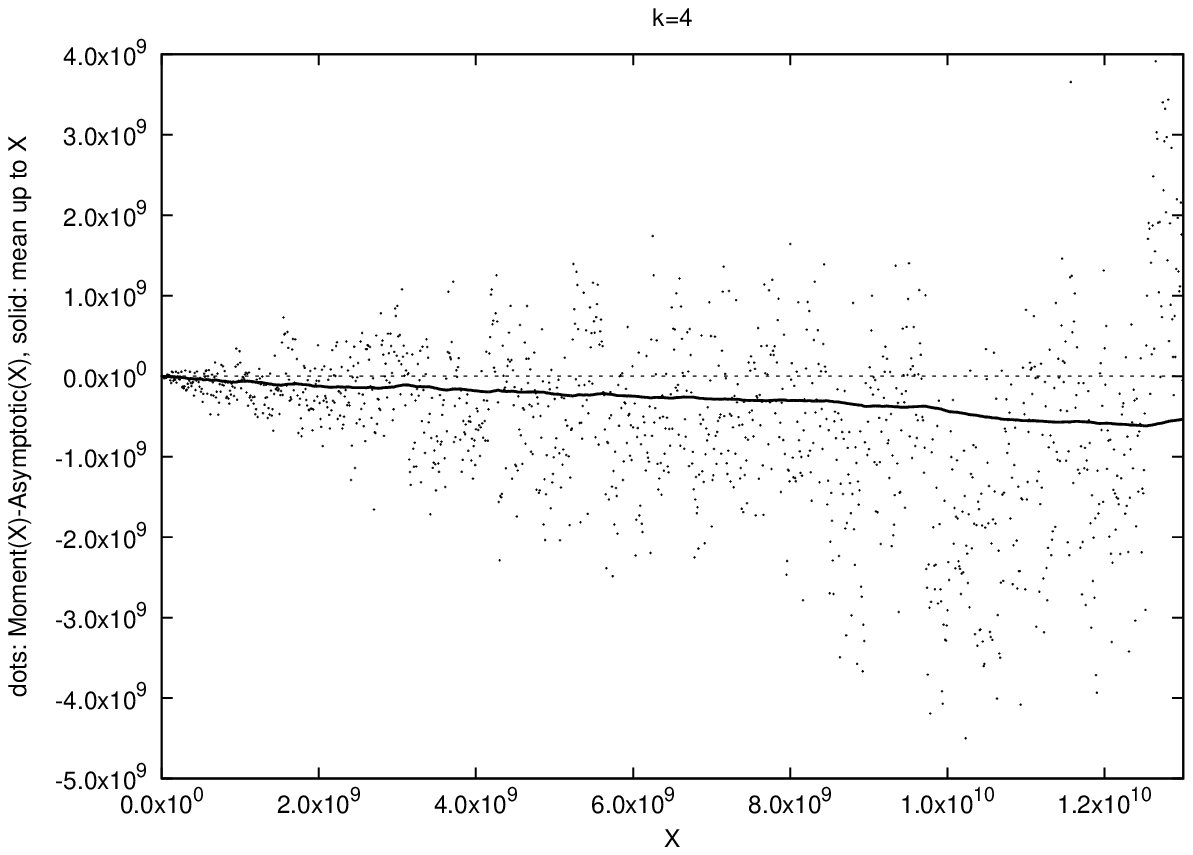}
   }
  \centerline{
   \includegraphics[width=0.48\textwidth,height=2in]{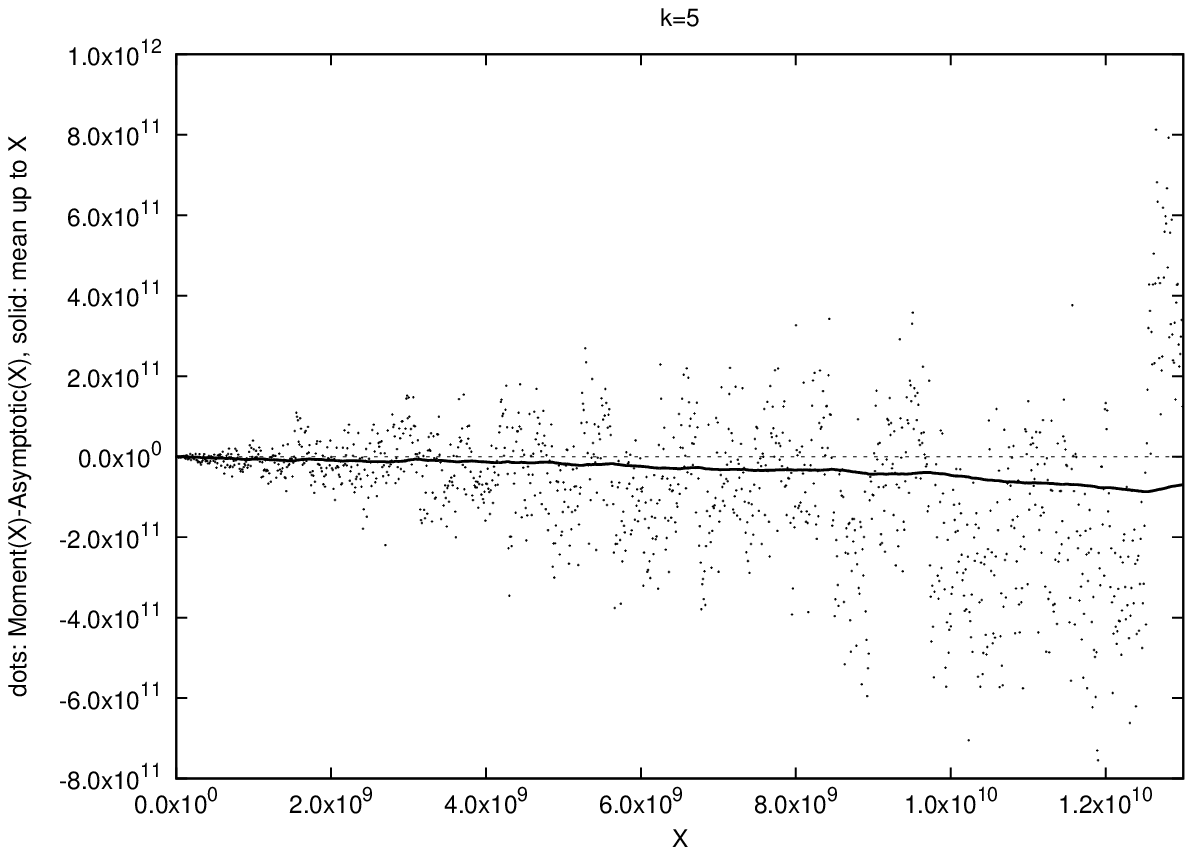}
     \includegraphics[width=0.48\textwidth,height=2in]{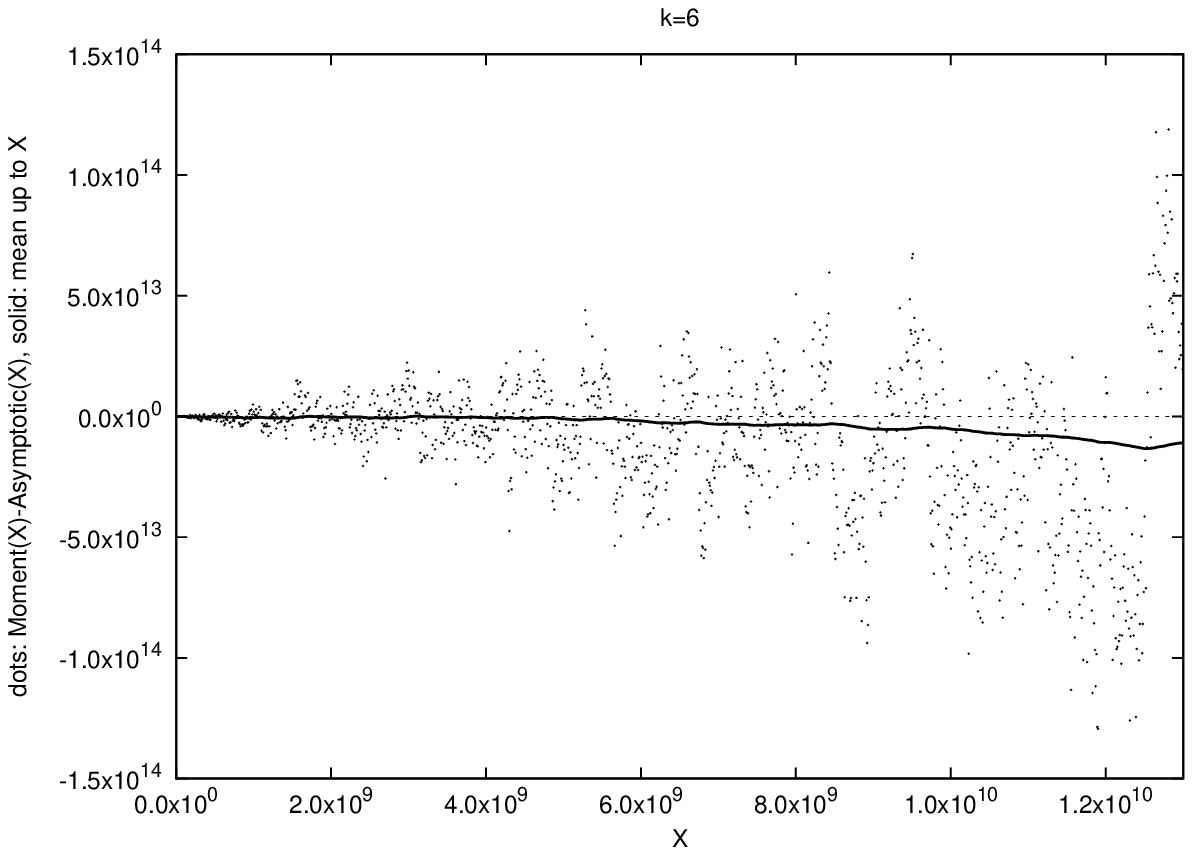}
  }  
  \centerline{
    \includegraphics[width=0.48\textwidth,height=2in]{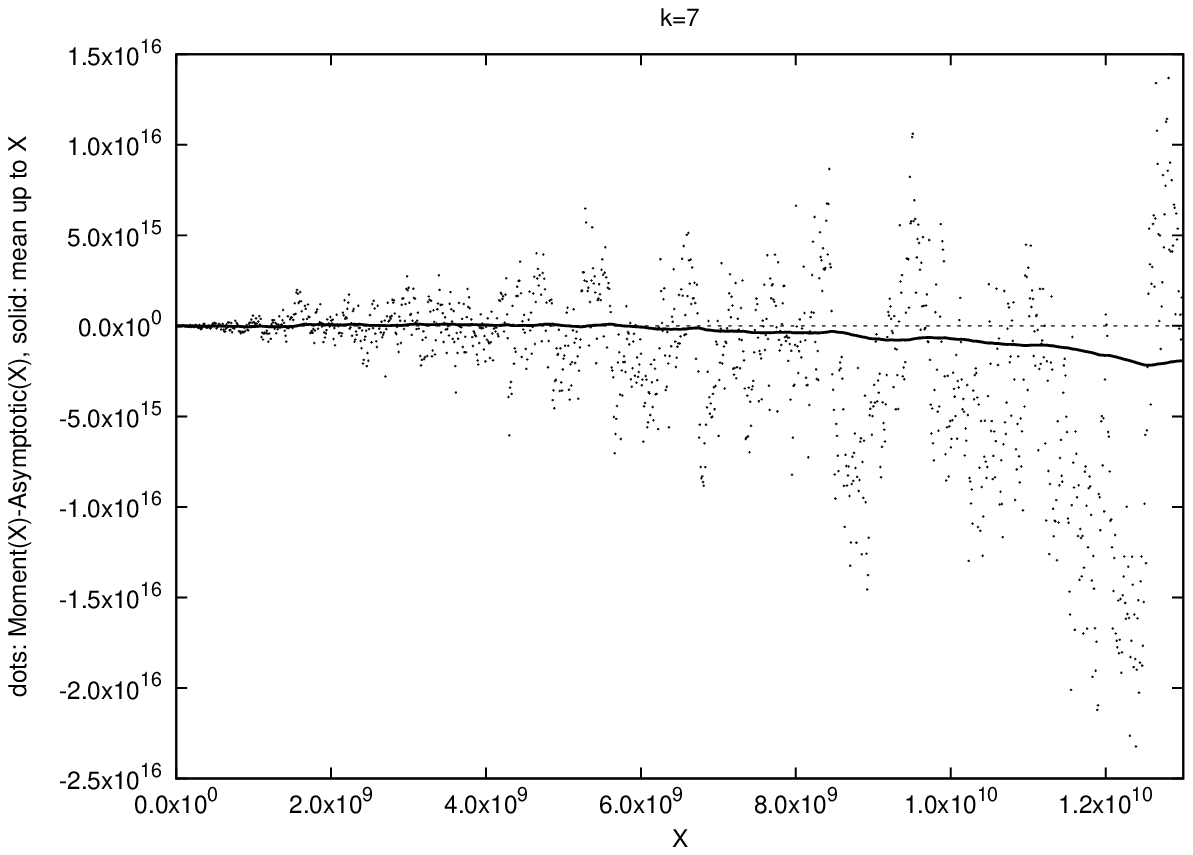}
     \includegraphics[width=0.48\textwidth,height=2in]{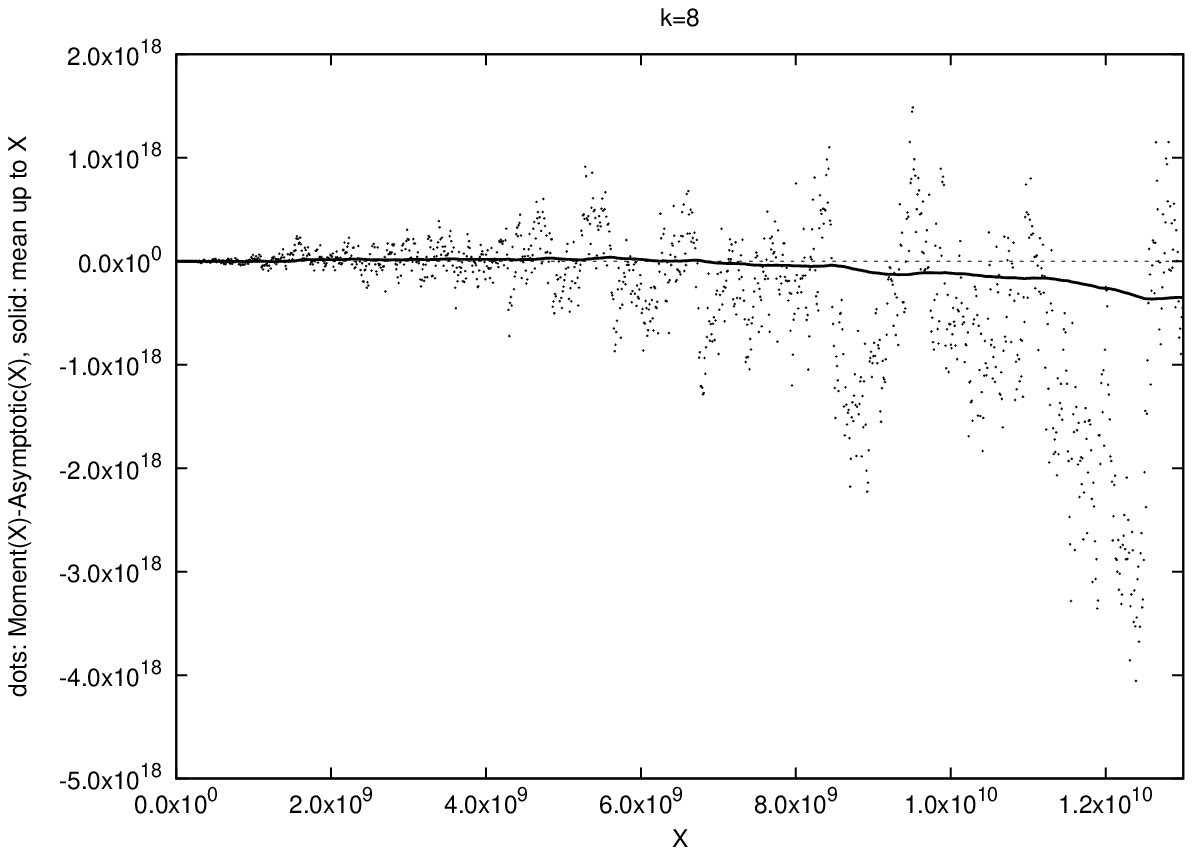}
  }
  \caption[Plot of the difference $\Delta_+(k,X)$ for $k=1,\ldots, 8$]{These
  plots depict the difference $\Delta_+(k,X)$ between the numerically compute moments and
  the CFKRS prediction, for $k=1,\ldots, 8$ and $d>0$,
  sampled at $X =  10^7, 2\times 10^7, \ldots, 1.3\times 10^{10}$. The horizontal axis
  is $X$, the vertical axis is the difference $\Delta_+(k,X)$, and the solid curve
  is the mean up to $X$ of the plotted differences (see the discussion above).}\label{fig: fig4}
\end{figure}

\newpage
\begin{figure}[H]
    \centerline{
        \includegraphics[width =1\textwidth,height=3.7in]{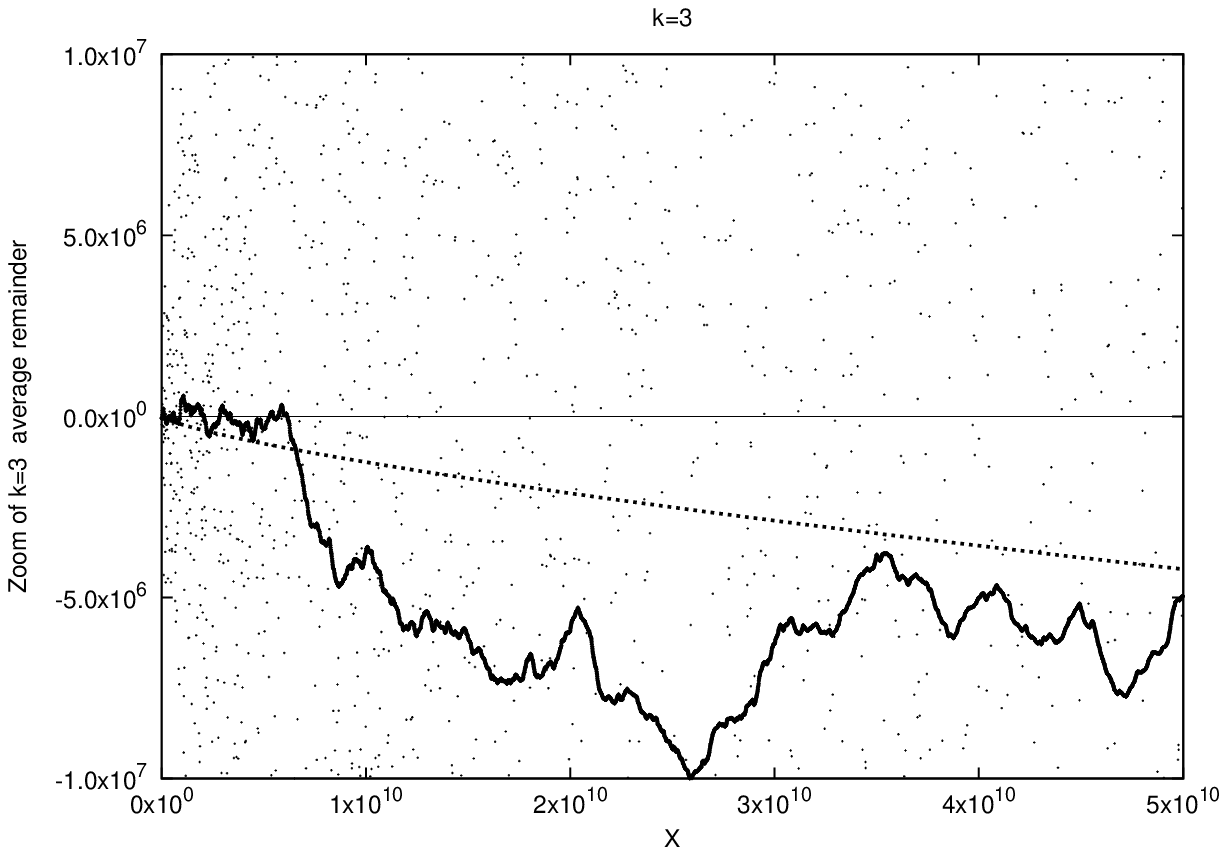}
    }
    \centerline{
       \includegraphics[width=1\textwidth,height=3.7in]{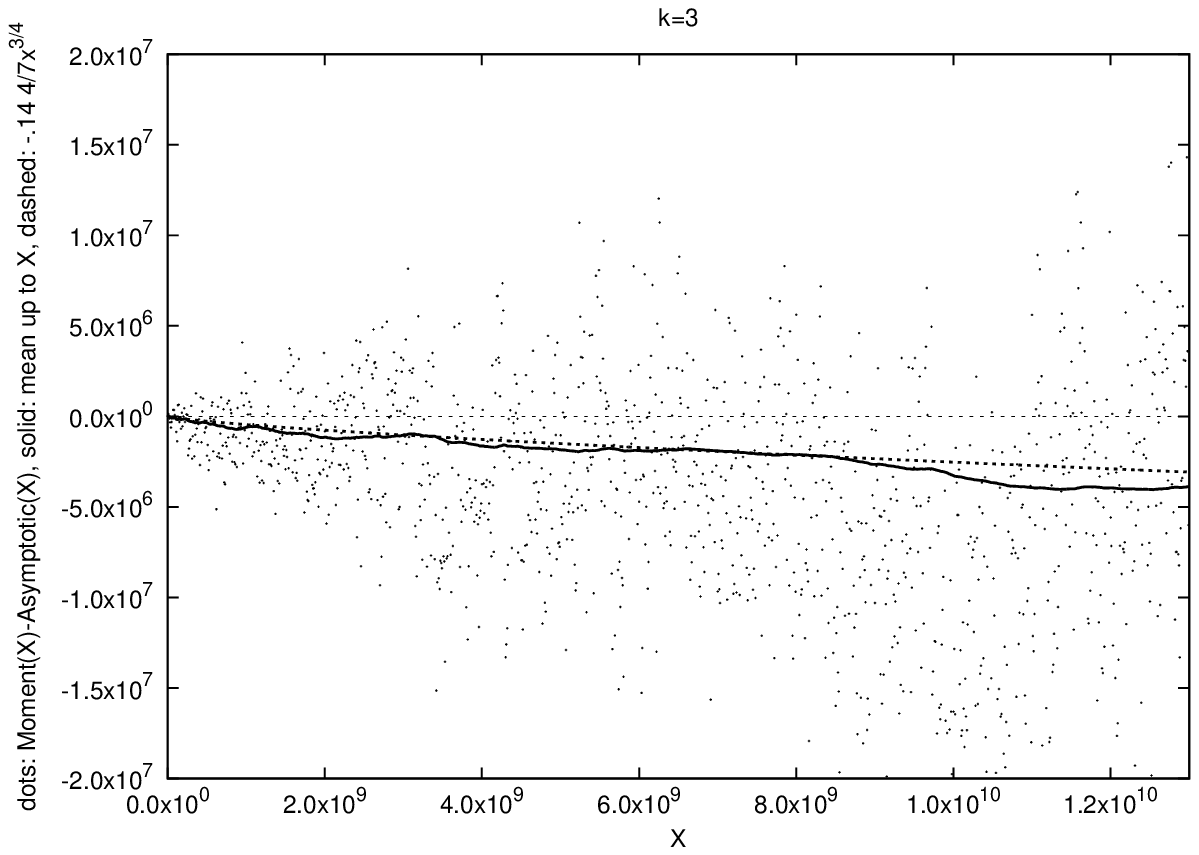}
    }
    \caption[Plot of $\Delta_\pm(3,X)$]{This graph depicts the remainder term for
    the third moment for $d<0$ (top) and $d>0$ (bottom), i.e. $\Delta_-(3,X)$ and $\Delta_+(3,X)$.
    The solid line is the average remainder, and the dashed line is Zhang's prediction.
    See the discussion above concerning the quality of the fit.
    }\label{fig: fig5}
\end{figure}

\newpage
\thispagestyle{empty}
\begin{figure}[H]
 \centerline{
    \includegraphics[width=.48\textwidth,height=2in]{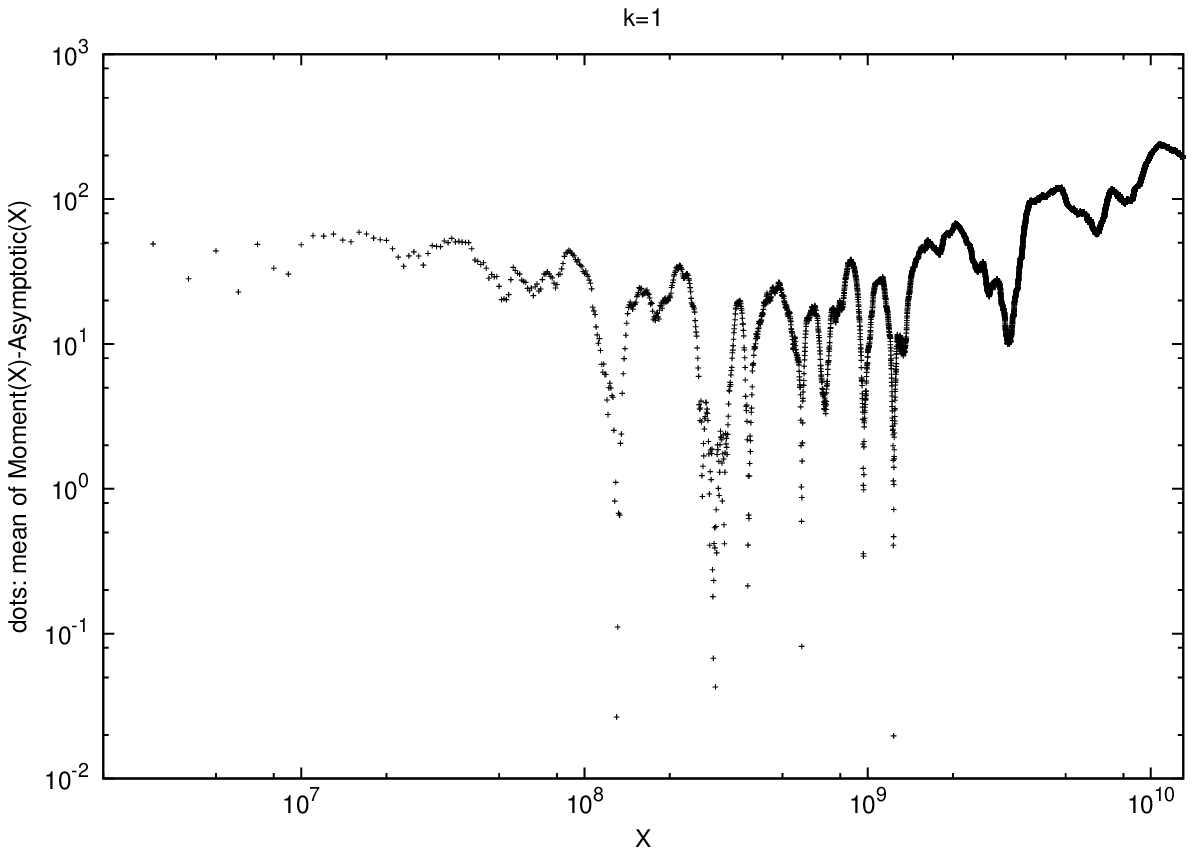}
    \includegraphics[width=.48\textwidth,height=2in]{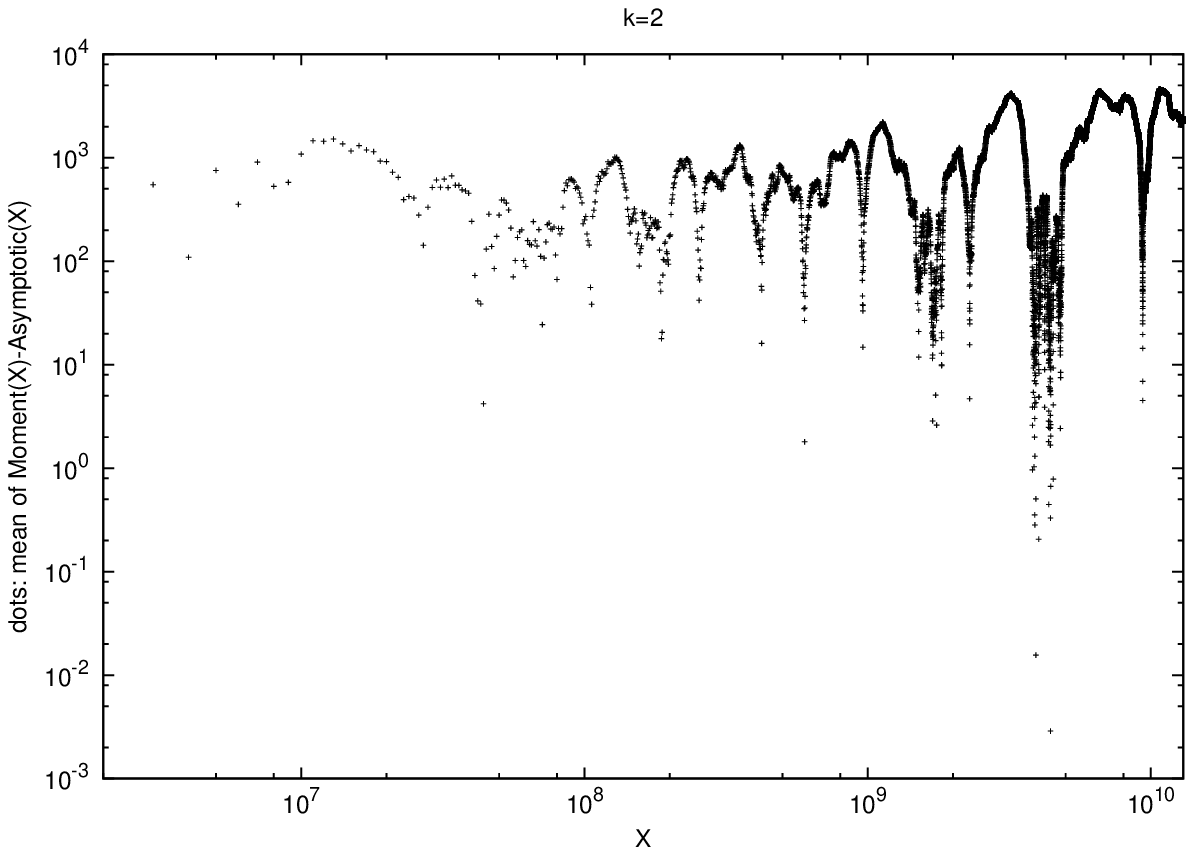}
  }
 \centerline{
    \includegraphics[width=.48\textwidth,height=2in]{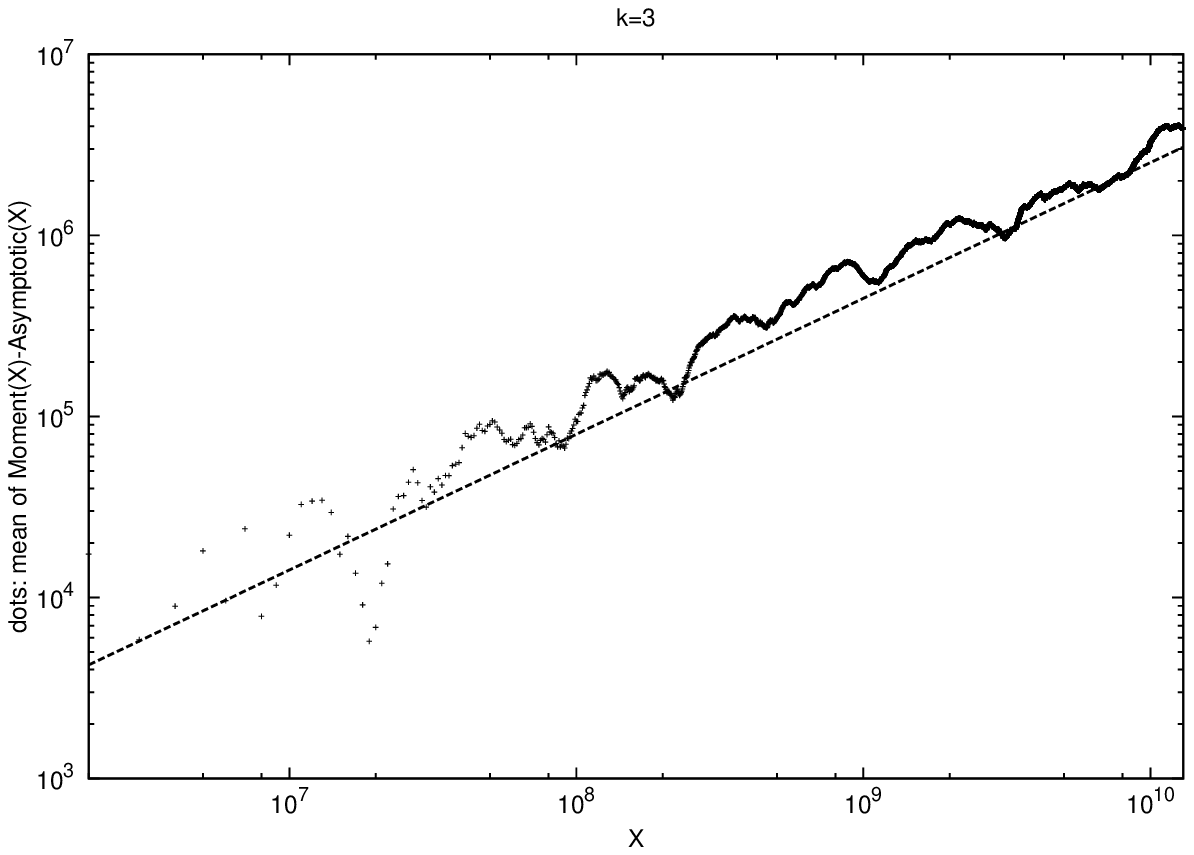}
    \includegraphics[width=.48\textwidth,height=2in]{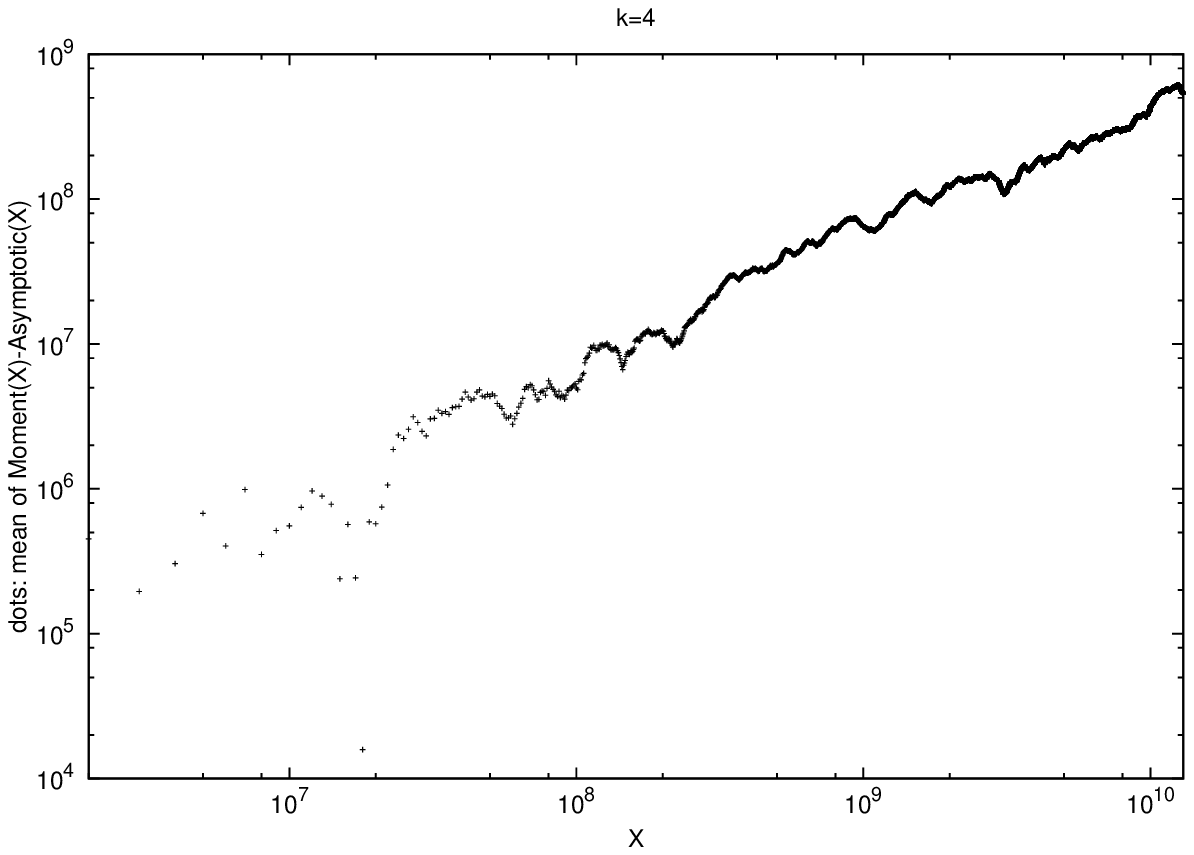}
  }
 \centerline{
    \includegraphics[width=.48\textwidth,height=2in]{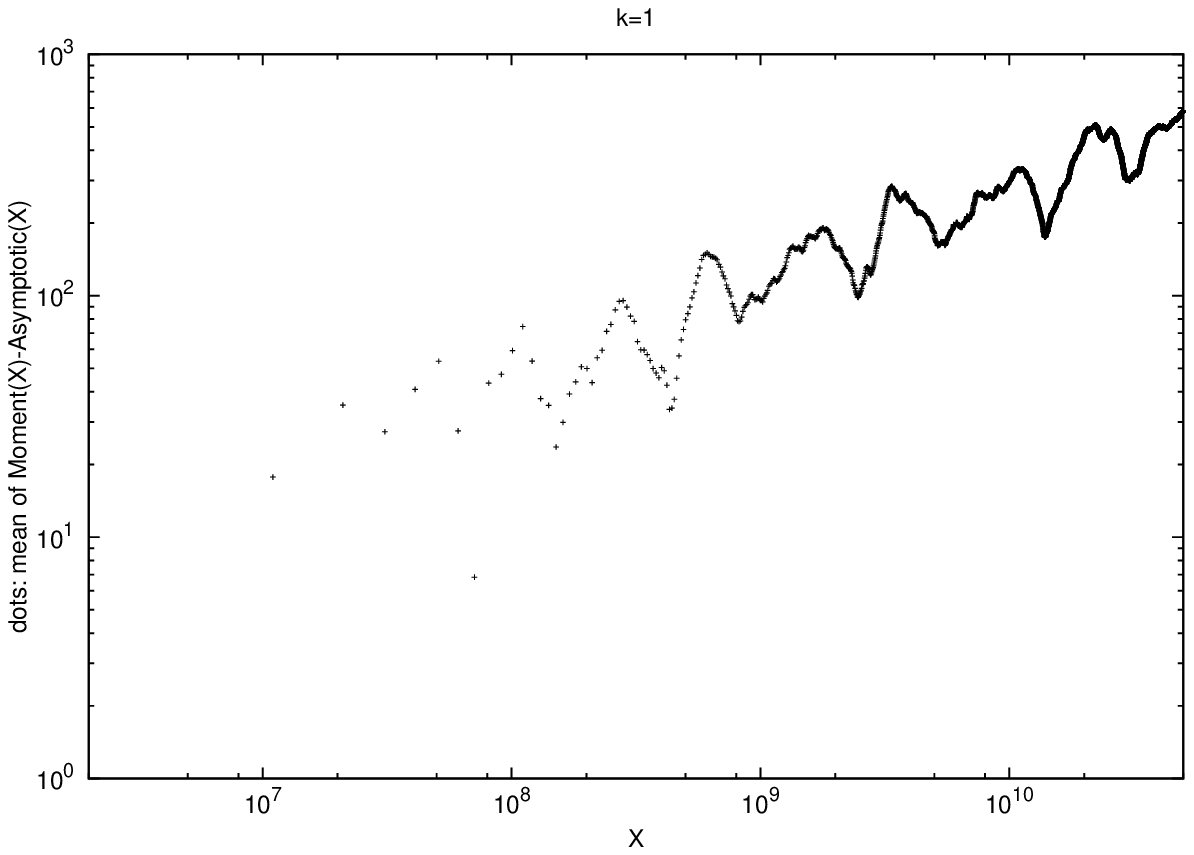}
    \includegraphics[width=.48\textwidth,height=2in]{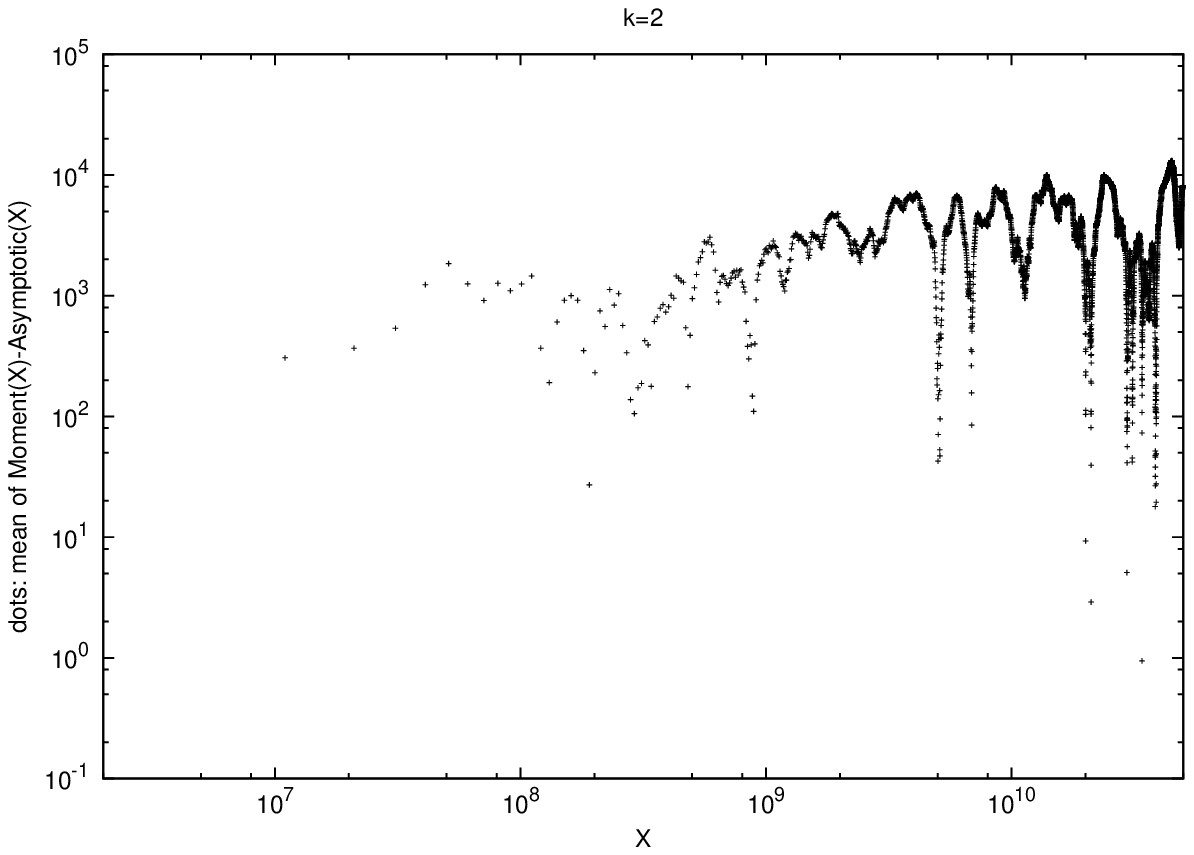}
  }
 \centerline{
    \includegraphics[width=.48\textwidth,height=2in]{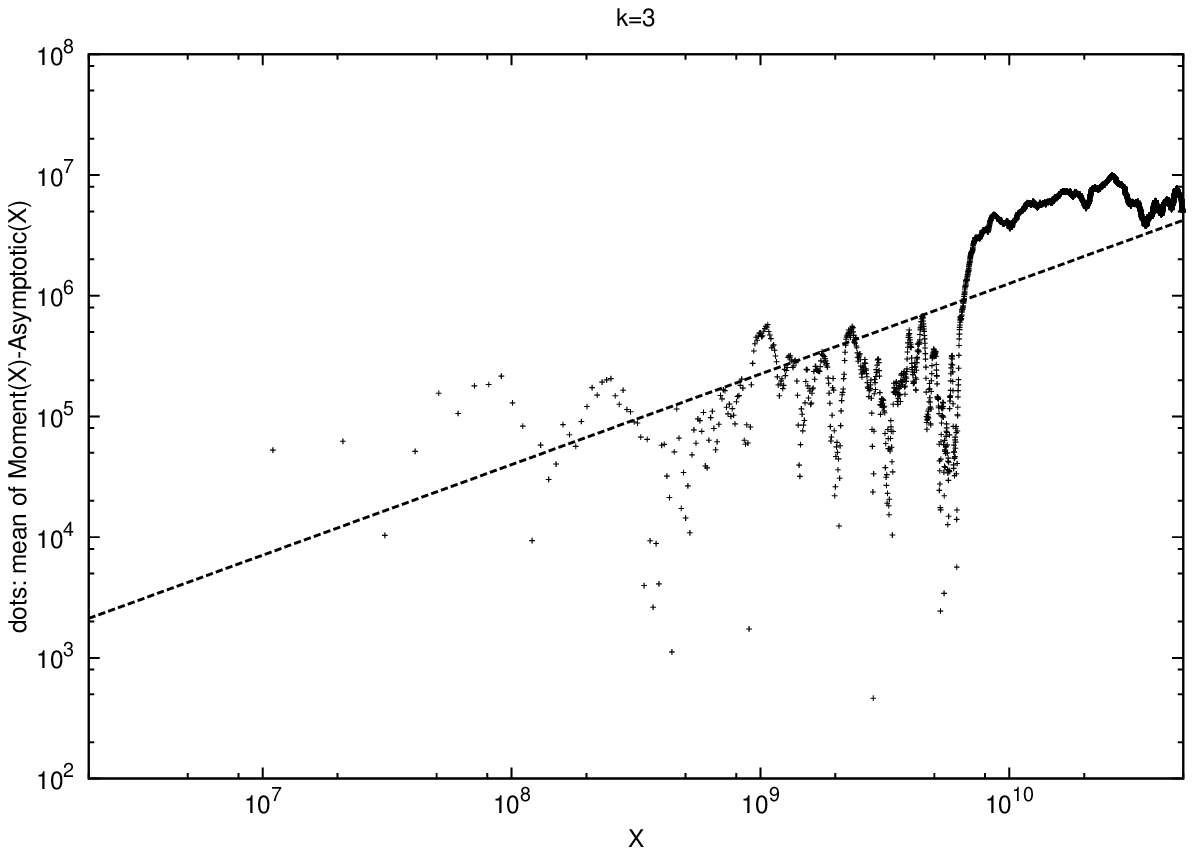}
    \includegraphics[width=.48\textwidth,height=2in]{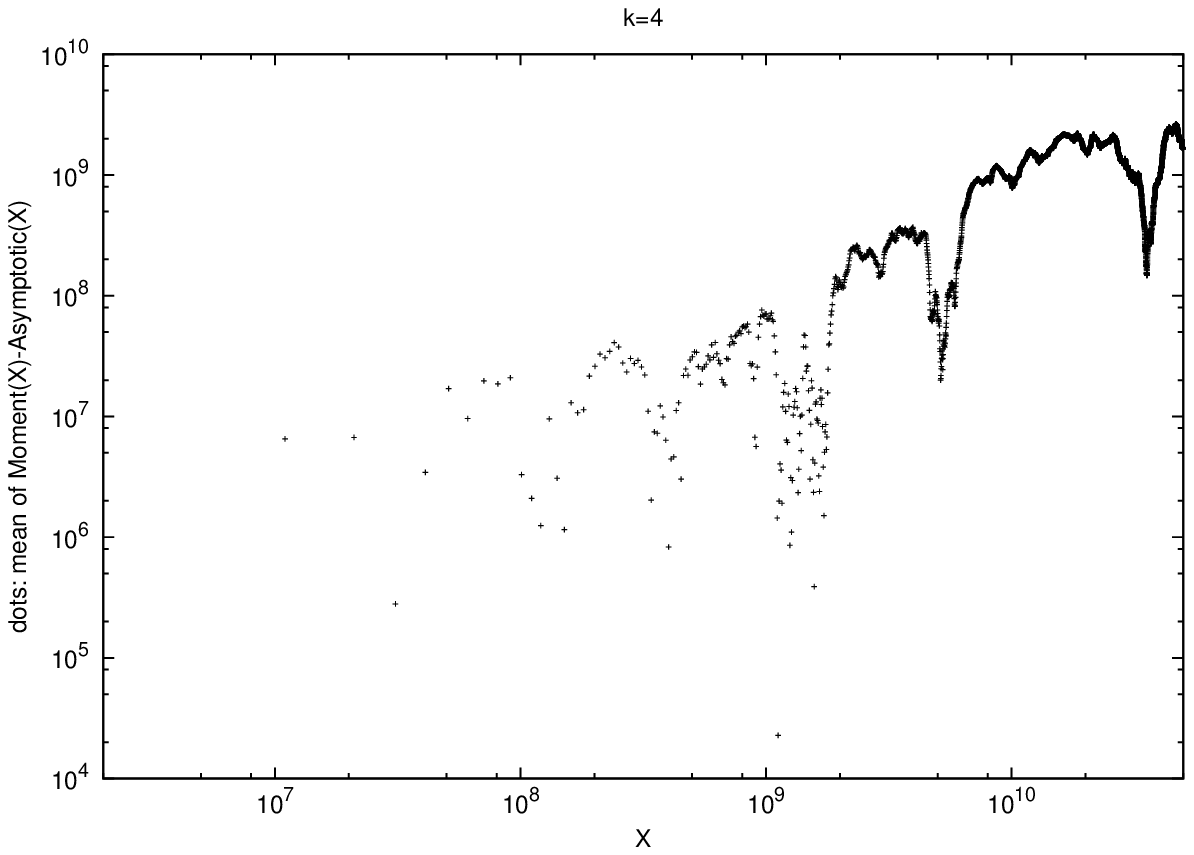}
  }
   \caption[log log plot of the average reaminder]{Plots, on a log log scale, of
   the absolute value of the average remainder term depicted in Figures 3--4,
   for $1 \leq k \leq 4$, $d>0$ (top four plots), and $d<0$ (bottom four plots).
   For the 3rd moment we compare to Zhang's predictions of $.14 \frac{4}{7} x^{3/4}$
   (3rd plot) and $.07 \frac{4}{7} x^{3/4}$ (7th plot).
   For the 1st moment, $d<0$, there seems to be a bias of size roughly $x^{1/4}$.
   }\label{fig: fig6} \end{figure}

\newpage

\section{Our Computational Formulae}

The computations for the moments of $L(1/2,\ch_d)$
hinge on the efficient computation of $L(1/2,\chi_d)$ itself for many discriminants $d$.
This computation is split into two cases according to whether $d$ is positive
or negative. In the former case we calculate $L(1/2,\chi_d)$ using a
smooth approximate functional equation for $L(s,\ch_d)$, which is
representable in terms of the incomplete gamma function. In the latter case,
we consider the Dedekind zeta function for the associated quadratic field, and reduce the
computation of $L(1/2,\chi_d)$ to a sum over binary quadratic forms, and the
$K$-Bessel expansions of their Epstein zeta functions as determined by Chowla and
Selberg \cite{cs: epzeta}.

Testing the conjectures described in the introduction also involves numerical
values for the coefficients of the polynomials $Q_\pm(k,x)$. We reran
the program used in~\cite{cfkrs: intmom} on a faster machine and for
a longer amount of time to get slightly more accurate coefficients
for these polynomials.

\subsection{Computational Formula for $L(1/2,\ch_d)$, $d<0$}

Let
$$
    \zeta_{\Q({\sqrt{D}})}(s)=\zeta(s)\ L(s,\chi_{d}) \label{eq: algDede}
$$
be the Dedekind zeta function of the quadratic number field $\mathbb Q({\sqrt{D}})$,
and $h(d)$ the corresponding class number.

Let $a_j m^2 + b_j mn + c_j n^2$, $j=1,\dots,h(d),$
be representatives for the $h(d)$ equivalence classes of primitive positive
definite binary quadratic forms of discriminant $b_j^2- 4a_j c_j=d<0$
\cite{land: elnumthe}.

Dirichlet proved (see also \cite{hdaven: mnumth}, Chapter 6) that:
\begin{equation}
    \label{eq:Dirichlet}
    \zeta_{\Q(\sqrt{D})}(s)=\frac{1}{\omega}\sum_{j=1}^{h(d)}\ \sum{}^{'} \biggl(a_j{m^2 +
    b_j mn + c_j n^2 \biggr)}^{-s}, \qquad \Re{s} > 1
\end{equation}
where
$$
    \omega=\begin{cases}
        2 & d<-4, \\
        4 & d=-4, \\
        6 & d =-3,
    \end{cases}
$$
\noindent and where $\sum{}^{'}$ denotes the sums over all pairs
$(m,n) \in \mathbb Z^2$, $(m,n) \ne (0,0)$.

Chowla and Selberg \cite{cs: epzeta} obtained the meromorphic
continuation of the Epstein zeta function
$$
Z(s):=\sum{}^{'} (am^2 + bmn + cn^2)^{-s}, \qquad \Re{s}>1,
$$
with $d=b^2-4ac<0$, $a,c>0$, by giving an expansion for $Z(s)$
as a series of $K$-Bessel functions.  Specifically, they proved that
\begin{equation}
    \label{eq:Z}
    Z(s)=2\zeta(2s)a^{-s} + \frac{2a^{s-1}\ \sqrt\pi}{\Gamma(s) \Bigl(|d|^{1/2}/2 \Bigr)^{2s-1}}
    \ \zeta (2s-1)\Gamma (s - 1/2) + B(s)
\end{equation}
where
\begin{align}
    B(s)&= \frac{8 \pi^s 2^{s-1/2}}
    {a^{1/2}\Gamma(s)\ |d|^{\frac{2s-1}{4}}} \sum_{n=1}^\infty
    n^{s-1/2}
    \ \sigma_{1-2s}(n)\cos\biggl(\frac{n \pi b}{a}\biggr)K_{s-1/2}
    \biggl(\frac{\pi n \ |d|^{1/2}}{a}\biggr), \\
    \sigma_\omega(n) &= \sum_{m|n} m^\omega,
\end{align}
and
$$
K_\omega(z)= \frac{1}{2}\int _0^\infty \ {\exp\biggl({-\frac{z}{2}(y+ 1/y
    )}\biggr) {y^{\omega-1}} dy, \qquad \Re{z} >0.}
$$

$K_{s-1/2}(x)$ decreases exponentially fast as $x \to \infty$,
uniformly for $s$ in compact sets, and
the above expansion gives $Z(s)$ as an analytic function throughout $\C$ except
for a simple pole at $s=1$. Note that the poles at $s=1/2$ in \eqref{eq:Z} of
the terms with $\zeta(2s)$ and $\Gamma(s-1/2)$ cancel out.

Specializing the above formula to $s = 1/2$ gives
\begin{align}
   Z(1/2) &= \fr{2}{a^{\fr{1}{2}}}\pr{\gm +
   \log\pr{\fr{\lr{d}^{\fr{1}{2}}}{8\pi a}}} \notag\\ &+
   \fr{8}{a^{\fr{1}{2}}}\sum_{n\geq 1} \s_0(n) \cos\pr{\fr{\pi n b}{a}}
   K_0\pr{\fr{\pi n\lr{d}^{\fr{1}{2}}}{a}}. \label{eq: expansion}
\end{align}
Substituting this into \eqref{eq:Dirichlet}, for a given a set of representative
quadratic forms, one for each equivalence class, yields a formula
that can be used to numerically compute $L(1/2,\chi_d)$.

An explicit bound on $K_0(x)$ can be obtained as follows.
\begin{equation}
    K_0(x)=\int_1^\infty \ \exp\biggl(-\frac{x}2 (y+1/y)\biggr) {\frac{dy}{y}}
\end{equation}
so that
\begin{equation}
    \label{eq:K-Bessel decay}
    |K_0(x)| < {\sqrt{\frac{\pi}{2x}}e^{-x}}.
\end{equation}
The last inequality can be seen by writing $y+1/y = (y^{1/2}-y^{-1/2})^2+2$,
changing variables $u=x^{1/2}(y^{1/2}-y^{-1/2})$, so that
$dy/y =x^{-1/2} 2 du/(y^{1/2}+y^{-1/2})$,
and using, from the AGM inequality,  $y^{1/2}+y^{-1/2} \geq 2 $.

Next, we discuss implementation issues and complexity arising from this formula.

Lagrange proved that each of the equivalence classes of primitive positive definite
forms of discriminant $d$ contains exactly one form $ax^2+bxy+cy^2$ for which
$-a < b \leq a < c$ or $0 \leq b \leq a=c$. Roughly, this is the set
$ 0 \leq |b| \leq a \leq c$, with some exceptions.
Furthermore, $a<(|d|/3)^{1/2}$. Recall that in this context primitive means that
$\gcd(a,b,c)=1$.

Therefore, let $A(X)$ denote the set of triples:
\begin{equation}
    A(X):= \{(a,b,c) \in \Z^3 | 4ac-b^2\leq X, \text{$-a < b \leq a < c$ or $0 \leq b \leq a=c$}   \}
    \label{eq: A}
\end{equation}
and $A'(X)$ the set of primitive triples in $A(X)$:
\begin{equation}
    A'(X):= \{(a,b,c) \in A(X) | \gcd(a,b,c)=1   \}.
    \label{eq: A'}
\end{equation}

Our first step was to distribute the computation across several processors,
each one handling a range of discriminants, in order
to speed up the computation and also to reduce the memory requirements per processor.
Therefore, suppose $0 < -d \leq X$ for some positive integer $X$,
and let $\Delta X$ be a positive integer dividing $X$.
We partitioned the interval into blocks, $X_{i-1} < \lr{d} \leq X_i$, of equal length
$\Delta X = X_i - X_{i-1}$ for $i=1,2,\ldots,m$ (so $X_m := X$).

We then looped through all integers $(a,b,c)\in A'(X)$,
with corresponding $|d|$ lying in $(X_{i-1},X_i]$, i.e. satisfying 
the following properties:
\begin{equation}
    \label{eq: loops}
    0< a \leq \sq{\fr{X_i}{3}}, \qquad 0 \leq |b| \leq a
    \leq c, \qquad \fr{b^2 + X_{i-1}}{4a} < c \leq \fr{b^2 + X_i}{4a},
\end{equation}
(taking care to throw away the terms above at the endpoints that are not in
$A'(X)$, for example, terms with $-b=a$).

We computed $d = b^2-4ac$ and updated
$L(1/2,\chi_d)$, stored in an array, using \eqref{eq: expansion} to calculate
the corresponding contribution to~\eqref{eq:Dirichlet} from the triple $(a,b,c)$.

Combining $a<(|d|/3)^{1/2}$ with the exponential decay
of $K_0(x)$ in \eqref{eq:K-Bessel decay} shows that
few terms are needed to compute~\eqref{eq: expansion}
to given precision.
For example, the terms in~\eqref{eq: expansion}
with $n\geq 7$ contribute, in absolute value, less than $10^{-15}$
to the sum, and can be ignored. Smaller $a$ require even fewer terms.

Because we are evaluating $K_0(x)$ for a limited range of values
and to machine double precision (15-16 digits), we used a precomputed 
table of the first five terms
of several thousand Taylor series expansions:
\begin{equation}
    \label{eq:Ktaylor}
    K_0(x) \approx K_0(x_j) + K_0^{(1)}(x_j) (x-x_j)+ \ldots (K_0^{(4)}(x_j)/4!)  (x-x_j)^4,
\end{equation}
each one centred on a point $x_j$ of the form $x_j=j/200$ with
$5< x_j <37$, $j \in \Z$. The interval $[5,37]$ was used, because, on one end,
$\pi 3^{1/2}>5$, the lhs being a lower bound for the smallest possible $x$
for which we would need to evaluate $K_0(x)$. We chose $37$ for the other end of
the interval because $\exp(-37)<10^{-16}$, i.e. smaller than our desired precision.

While our code was written in {\tt C++},
the precomputation of these Taylor expansions were carried out in Maple, with the coefficients
stored in a file that could be read into our {\tt C++} program. Our {\tt C++} code
was compiled with the GNU compiler {\tt GCC}.

Only five terms were computed and stored because our Taylor expansions were
always applied with $|x-x_j|<1/400$. Combined with the exponential decay
of the Taylor coefficients (as a function of $x$), at most 5 terms,
and often fewer, were needed to evaluate the sum to within $10^{-15}$.

We also make note of a few of the hacks that helped to increase the speed of our
program:

\begin{itemize}

     \item For a given $a,b$, only one cosine needs to be computed. Indeed,
     given $\cos(\pi b/a)$, we can compute $\cos(\pi nb/a)$,
     for $n=1,2,\ldots, 7$, using standard trigonometric identities.
     For instance, the double angle identity computes the expression for $n=2$.

     \item To test for primitivity, we must check whether $\gcd(a,b,c)=1$.
     If one computes $\gcd(a,b)$ outside the $c$ loop as previously
     mentioned, then for a given $\gcd(a,b)$, we can use
     $$\gcd(a,b,c) = \gcd(\gcd(a,b), \, c \mod \gcd(a,b)),$$
     and thus compute, and then store within the $c$ loop,
     at most one gcd per residue class mod $\gcd(a,b)$.

     \item
     When reading an array, the computer loads blocks of consecutive
     bytes of the array from RAM into the CPU's cache where it can be accessed
     quickly by the CPU.
     On profiling, we found, after optimizing and streamlining
     the bulk of our code, that a significant amount of time was being spent
     accessing the array in which we were storing $L(1/2,\chi_d)$ so
     as to increment it by the contribution from a given triple $(a,b,c)$.
     The reason was that, as the inner $c$ loop increments by 1,
     the value of $d=b^2-4ac$ changes by $4a$, i.e. often a large step.
     Non-sequential array accesses are expensive timewise.
     We were able to significantly
     decrease the time spent on array accesses by anticipating the subsequent $d$'s, and
     using {\tt GCC}'s `\_\_builtin\_prefetch' function to fetch the corresponding array
     entry for $L(1/2,\chi_d)$ eight turns in advance - the eight was determined
     experimentally on the hardware that we used and for the range of $d$'s that we considered.

     \item The computation of \eqref{eq: expansion} involves $\log(|d|)$. Therefore,
     we precomputed these and stored them in an array, but were again faced with the same kind of
     expensive memory accesses as for the values of $L(1/2,\chi_d)$. Rather than
     prefetch these separately, we created a {\tt C++} struct to hold both $L(1/2,\chi_d)$ and
     $\log(|d|)$ together. That way a single prefetch would load both at once.

     \tb{Remark.} On combining the last two hacks, the array access portion of our code sped up by a factor of 4,
     and the overall running time of the program sped up by a factor of 2. These two hacks were
     the last two implemented, and the speed up achieved indicates how expensive non-sequential
     memory accesses can be, and how optimized the rest of our code was.

     \item To avoid repeatedly checking whether the quantity $d=b^2-4ac$ is a fundamental
     discriminant, we precomputed, for each block $X_{i-1} < \lr{d} \leq X_i$ whether
     $d$ is a fundamental discriminant and stored that information in an array of boolean
     variables. We essentially sieved for squarefree numbers and this
     can be done, for each block of length $\Delta X$, in $O(\Delta X)$ steps,
     because the sum of the reciprocal of the squares converges. Hence, doing so across
     all blocks up to $X$ costs $O(X)$ arithmetic operations and array accesses.
     No prefetching was used on this array as it did not seem to give a benefit, perhaps
     because the array constists of single bit boolean variables, rather than 64 bit doubles
     of the $L$-value and $\log(|d|)$ arrays, and fits more easily within cache.

     \item Since $\cos(x)$ is an even function and $b$ gets squared in the
     discriminant equation $\lr{d} = 4ac - b^2$, we can group $\pm b$
     together, when possible, and restrict our attention to non-negative $b$ values.
     Only a relatively small subset of triples cannot be paired in this fashion, namely
     when $a=b$, $b=0$, or $c=a$.

     \item Terms such as $$\fr{2}{a^{\fr{1}{2}}}\pr{\gm - \log\pr{8\pi a}}$$
     appearing in the leading term of (\ref{eq: expansion}), depend solely on
     $a$. As such, it is to our advantage to compute this, and all other terms
     depending solely on $a$, outside the $b$ and $c$ loops. Similarly, we
     compute expressions like $\gcd(a,b)$ outside the $c$ loop, and so on.
     While this is standard,
     we took it to a meticulous extreme to save on as many arithmetic operations as possible.

\end{itemize}

\subsection{Complexity for $d<0$}\label{sec-negcomplexity}

First observe that the number of candidate triples $(a,b,c)$ defined by~\eqref{eq: A}
that we must loop over, satisfies
\begin{equation}
    |A(X)| \sim \frac{\pi}{18} X^{3/2}.
    \label{eq: A asympt}
\end{equation}
See \cite{kuhleitner} for a discussion on this counting problem.
Furthermore, the number of triples that survive the condition $\gcd(a,b,c)=1$
is:
\begin{equation}
    |A'(X)| \sim \frac{\pi}{18\zeta(3)} X^{3/2}.
    \label{eq: A' asympt}
\end{equation}
The latter asymptotic formula was essentially stated by Gauss in his Disquitiones
in connection to the sum of class numbers $A'(X) = \sum_{-X<d<0} h(d)$.
A proof, with a lower term and a bound on the remainder,
as well as further references can be found in \cite{ci}.

We can get a lower bound for the amount of computation required by
simply counting the number of triples $(a,b,c)$ that are considered. Furthermore,
because we are pairing together $\pm b$, the number of triples that survive the
gcd condition is roughly half of~\eqref{eq: A' asympt}, i.e.
\begin{equation}
    \sim
    \frac{\pi}{36\zeta(3)} X^{3/2}
    \label{eq: number of triples}
\end{equation}
The relatively small constant of $\pi/36 \zeta(3)$ helps to makes this approach
very practical.

While the above asymptotic gives a lower bound on the number of operations
of our method
in computing all $L(1/2,\chi_d)$, for $0 < -d \leq X$, it ignores the amount of
work needed for each triple.
Most of the work involves: checking bounds on each loop,
testing whether $\gcd(a,b,c)=1$ and whether $d=b^2-4ac$ is a fundamental discriminant,
and carrying out simple arithmetic and array accesses
related to the evaluation of ~\eqref{eq: expansion}.

Each arithmetic operation can be done in polynomial time in the size of the
numbers involved (in fact, $\log(X)^{1+\epsilon}$ by using the FFT).
However, in the range of
discriminants we considered ($|d| < 5 \times 10^{10}$), the bit length is
quite small: 32 bit {\tt C++} ints
for $a,b,c$ sufficed and 64 bit long longs were used for discriminants.
We also used 64 bit machine doubles for the floating point arithmetic.
Therefore all arithmetic was carried out in hardware.

Recall that we are assuming that
$0 < -d \leq X$, and partitioning this interval as:
$$\underbrace{1,\ldots, \Dl{X}}_{\tx{Block 1}},
\underbrace{\Dl{X} + 1, \ldots, 2\Dl{X}}_{\tx{Block 2}}, \ldots,
\underbrace{(m-1)\Dl{X} + 1, \ldots, m\Dl{X}}_{\tx{Block m}}, \ldots,$$ where
$\Dl{X}$ is a positive integer assumed to divide $X$ (in practice one can
take $X$ and $\Delta X$ to be powers of, say, ten). The number of blocks equals
$X/\Dl{X}$, where, for reasons made clear below, we will eventually take $\Dl{X} \gg  X^{1/2} \log(X)^2$.

As mentioned earlier,
we precomputed, via sieving, a table of fundamental discriminants
for each block of length $\Delta X$ using $O(\Delta X)$ arithmetic
operations and array accesses, hence $O(X)$ operations across all blocks,
i.e a relatively small cost compared to $O(X^{3/2})$.

The number of terms needed in the expansion of $K_0$ Bessel function
is proportionate to the number of digits of precision desired, and the number
of Taylor coefficients used in each Taylor series for $K_0(x)$ also depends on
the ouput precision. We worked with machine doubles and hence about 15-16
digits precision, and, as explained in the text near~\eqref{eq:Ktaylor}, we
used at most five terms in each Taylor series.

Next we consider the time used to carry out the $\gcd$ computations.
We described in the hacks listed above that, for each triple $(a,b,c)$, we computed
$\gcd(a,b,c)$ by first computing $\gcd(a,b)$ outside the
$c$-loop and then computing $\gcd(\gcd(a,b),c \mod \gcd(a,b))$ inside the
$c$-loop, with the latter calculation being performed at most once per residue
class modulo $\gcd(a,b)$, and then stored in the $c$-loop for subsequent use.

The number of $\gcd(a,b)$'s that are computed across all blocks is
\begin{equation}
    \label{eq:gcd bound1}
    \ll
    X \msum \asum \bsum 1 \ll X^2/\Dl{X}.
\end{equation}

Next, we compute $\gcd(a,b,c) = \gcd\pr{\gcd(a,b),c \mod \gcd(a,b)}$,
storing these values, inside the $a,b$ loops, for each residue class
$\mod \gcd(a,b)$ so as to only compute one $\gcd$ per residue class.
Therefore, the number of additional $\gcd$'s required is
\begin{equation}
    \ll \msum \asum\bsum \gcd(a,b). \label{eq: initialization}
\end{equation}

It is known that
\begin{equation}
    \pasum \bsum \gcd(a,b) \sim \fr{x^2\log x}{2\z(2)} \label{eq: gcd sum};
\end{equation}
see, for example, \cite{Bo} or \cite{Br}.
So, from (\ref{eq: initialization}) and the above asymptotic, we see that
the number of additional $\gcd$'s computations required across all blocks is
\begin{equation}
    \label{eq: gcd bound 2}
    \ll \msum m\Dl{X} \log\pr{m\Dl{X}} \ll \fr{X^2\log X}{\Dl{X}}.
\end{equation}
Therefore, combining with \eqref{eq:gcd bound1}, the total number of $\gcd$
calls is
\begin{equation}
    \label{eq: gcd total bound}
    \B\pr{\fr{X^2\log X}{\Dl{X}}}.
\end{equation}
Furthermore each $\gcd(a,b)$ can be computed, using the Euclidean algorithm, in
$\B\pr{\log(X)^2}$ bit operations, since the binary length of both $a$ and
$b$ is $O(\log{X})$, and so the total number of bit operations coming
from $\gcd$'s is $\ll X^2 \log(X)^3/\Delta X$.

Hence, we can make the overall time required for $\gcd$
evaluations an insignificant portion of the overall time by choosing
\begin{equation}
    \Delta X \gg X^{1/2} \log(X)^2,
\end{equation}
as the number of bit operations for the remaining work (looping
through $a,b,c$ and $m$, and carrying out the required integer and
floating point arithmetic) is then $$O(X^{3/2} \log(X)^{1+\epsilon}).$$
The $X^{3/2}$ accounts for the overall number of triples $(a,b,c)$
considered, and the $\log(X)^{1+\epsilon}$ for the cost of arithmetic on
numbers of bit length $O(\log{X})$. The implied constant depends on the
number digits of precision desired for the $L$-values.

While the best choice might seem to be to take $\Delta X$ equal to $X$
so as to minimize the number of $\gcd$ calls,
this would come at a substantial price.
First, such a large $\Delta X$ would prevent us from simply distributing the
computation across several processors, each one handling one block at a time.

Second, the memory (RAM) requirements needed would be enormous.
There is also an advantage to having arrays that can fit entirely or significantly
within the CPU's cache, so as to avoid too many expensive memory fetches
from RAM, and, even with smaller $\Delta X$, there is a tradeoff between
minimizing calls to the Euclidean algorithm and memory accesses.

We determined a good choice of $\Delta X$ experimentally, since, in practice,
the big-Oh constants in the above estimates depend on the speed of
individual arithmetic and memory operations on given hardware and context
in which they are called.

Nonetheless, since the Euclidean algorithm is very simple, and the remaining work
associated with looping, computing the $K$-Bessel function, and updating
$L$-values involves a moderate number of arithmetic and memory operations, one
expects that the benefit should be felt sooner rather than later. Indeed, we
found that, in our range of $d$'s, a choice that eliminated the $\gcd$'s as a
bottleneck, while not paying too high of a cache size penalty, was
$\Delta X = 10^6$, i.e. blocksizes of one million.

\subsection{Computational Formula for $L(1/2,\ch_d)$, $d>0$}

The proof of de la Vall\'ee Poussin of the functional equation for $L(s,\chi_d)$ imitates that of Riemann
for his zeta function. It yields the analytic continuation of $L(s,\chi_d)$ and also the following formula,
an example of a `smoothed approximate functional equation', useful for its evaluation:
\begin{equation}
    (d/\pi)^{s/2} \Gamma(s/2) L(s,\chi_d)
    = \sum_{n=1}^{\infty}
    \chi_d(n)( G(s/2, \pi n^2 /d) + G((1-s)/2, \pi n^2 / d ) )
\end{equation}
where $G(z,w)$ denotes the normalized incomplete gamma function
\begin{equation}
    G(z,w) := \int_1^\infty x^{z-1} e^{-wx} dx
    = w^{-z}\int_w^\infty x^{z-1}e^{-x} dx
    = w^{-z}\Gm(z,w), \quad \Re{w} > 0,
    \label{eq:G(z,w)}
\end{equation}
with $\Gm(z,w)$ the incomplete gamma function.
See for instance page 69 of \cite{hdaven: mnumth}, or Section 3.4.1 of \cite{rub: annumthe},
and use Gauss' formula for the Gauss sum, namely $\tau(\chi_d)=d^{1/2}$, when $d>0$.

Therefore, on specializing to $s=1/2$, we have
a smooth approximate functional equation for $L(1/2,\chi_d)$, namely
\begin{align}
    L(1/2,\chi_d) = 2\pr{\fr{\pi}{d}}^{\fr{1}{4}}
    \sum_{n\geq 1} \chi_d(n) \fr{G(1/4, n^2\pi/d)}{\Gm(1/4)} =
    2
    \sum_{n\geq 1} \fr{\chi_d(n)}{\sq{n}} \frac{\Gm(1/4, n^2\pi/d)}{\Gm(1/4)},
    \label{eq: pos-appformula}
\end{align}
valid for positive fundamental discriminants $d$.

To estimate the size of the terms being summed, first notice, from the definition, that
$G(1/4,w)>0$ for real $w$. Furthermore,
integrating by parts, gives an upper bound:
\begin{equation}
     G(1/4,w) = \frac{e^{-w}}{w}- \frac{3}{4w} \int_1^\infty e^{-wx} x^{-7/4} dx < \frac{e^{-w}}{w}
    \label{eq:G bound}
\end{equation}

This inequality tells us that the terms in~(\ref{eq: pos-appformula}) decrease exponentially fast
in the quantity $\pi n^2/d$, so that, roughly speaking, we need to truncate the sum when 
$n$ is of size $d^{1/2}$ to achieve a small tail.

Let us consider this estimate more carefully. Set
\begin{equation}
    f(t) = \fr{2}{\sq{t}} \fr{\Gm\pr{1/4, t^2\pi/d}}{\Gm\pr{1/4}}.
    \label{eq:f(t)}
\end{equation}
In light of bound~\eqref{eq:G bound}, we have
\begin{equation}
    |f(n)| < \frac{2}{\Gamma(1/4)} \left(\frac{d}{\pi}\right)^{3/4} \frac{e^{-\pi n^2/d}}{n^2}.
    \label{eq:f bound}
\end{equation}
Let the number of working digits be labelled as `Digits'. Hence, for
\begin{equation}
    \label{eq:n drop}
    n> \sq{\fr{d}{\pi}\log(10) \cdot \Digits}
\end{equation}
we generously have
\begin{equation}
    f(n) < 10^{-\Digits}.
\end{equation}
Furthermore, notice that
the terms start off, for smaller $n$ and large $d$, with
$$\frac{\Gm(1/4, n^2\pi/d)}{\Gm(1/4)} \sim 1.$$ Therefore, it does not make
sense to sum the terms beyond~\eqref{eq:n drop}, as those terms are
lost to numerical imprecision.

We must thus see to what extent the ignored tail end of the sum can
contribute to the value of $L(1/2,\chi_d)$.
Summation by parts yields
\begin{equation}
  \sum_{n\leq N} \chi_d(n)f(n) = f(N)\sum_{n\leq N} \chi_d(n) - \int_1^N
  \sum_{n\leq t} \chi_d(n) f'(t) dt. \label{eq: partialsum}
\end{equation}
So on letting $N\rt \infty$, we obtain
\begin{equation}
  L\pr{1/2, \ch_d} = -\int_1^\infty \sum_{n\leq t} \chi_d(n) f'(t) dt. \label{eq: lim-partial}
\end{equation}
Moreover, by subtracting (\ref{eq: partialsum}) from (\ref{eq: lim-partial}), we get a formula for the tail:
\begin{equation}
  \sum_{n = N+1}^\infty \chi_d(n)f(n) = -f(N)\sum_{n\leq N} \chi_d(n) -
  \int_N^\infty \sum_{n\leq t} \chi_d(n) f'(t) dt. \label{eq: tail}
\end{equation}

One could use the inequality of Polya-Vinogradov, Burgess, or even the trivial bound $\lr{\chi_d(n)}\leq 1$,
here to get a reasonable, but not optimal, estimate for the size of the tail.
However, something closer to the truth is obtained by using the conjectured
bound 
\begin{equation}
  \sum_{n\leq x} \chi_d(n) = \B\pr{x^{1/2} d^{\,\ep}}.  \label{eq: chi-conjecture}
\end{equation}
Combined with \eqref{eq:f bound} this gives
\begin{equation}
    f(N)\sum_{n\leq N} \chi_d(n) = \B\pr{\frac{d^{3/4+\ep}}{N^{3/2}} e^{-\pi N^2/d}}, \label{eq: firstterm}
\end{equation}
and, similarly,
\begin{equation}
  \int_N^\infty \sum_{n\leq t} \chi_d(n) f'(t)dt \ll
  d^{\,\ep}\int_N^\infty t^{\fr{1}{2}} f'(t) dt 
  \ll \frac{d^{3/4+\ep}}{N^{3/2}} e^{-\pi N^2/d}, \label{eq: secondterm}
\end{equation}
where we have applied integration by parts to get the last bound.
Applying these bounds to~\eqref{eq: tail} and choosing
\begin{equation}
  N = \sq{\fr{d}{\pi}\log(10) \cdot \Digits}, \label{eq: N}
\end{equation}
gives the following bound for~\eqref{eq: tail}:
\begin{equation}
    \B\pr{10^{-\Digits}\frac{d^\ep}{\Digits^{1/2}}}.
    \label{eq:tail bound}
\end{equation}
We therefore conclude that the tail isn't much bigger than an individual term, and, in principle,
we could compensate for the extra $d^\ep$ by taking Digits slightly
larger than our desired output precision, say by an amount equal to $\ep \log{d}/\log{10}$.

We remark that, by using the trivial estimate or Polya-Vinogradov inequality, we
could get rigorous estimates with explicit constants, but larger by a factor of roughly $d^{1/4}$
as compared to~\eqref{eq:tail bound}.  

\subsection{Cancellation and accuracy}

We can also use the above analysis to show that our approach to computing $L(1/2,\chi_d)$
using the smooth approximate functional equation
is well balanced, i.e. that little cancellation and hence loss of
precision takes place in summing~\eqref{eq: pos-appformula}.
We consider the maximum size that the partial sums can attain so as to
give us a sense of how many digits accuracy after the decimal place are attained
when working with Digits decimal places.

Consider the partial sums ~\eqref{eq: partialsum} (for
a general $N'$, not just our specific choice of $N$), apply the conjectured
bound~\eqref{eq: chi-conjecture}, and integrate by parts:
\begin{equation}
    \sum_{n\leq N'} \chi_d(n)f(n) \ll f(N') {N'}^{1/2} d^{\,\ep} + d^{\,\ep} t^{1/2} f(t) \Big|_1^{N'}
    + d^{\,\ep} \int_1^{N'} t^{-1/2} f(t) dt. \label{eq: estimate partial sum}
\end{equation}
Again, we can get a proven, though weaker, upper bound with explicit constants if we use
a proven bound rather than the conjecture~\eqref{eq: chi-conjecture}.

Next, notice that $\Gamma(1/4,x) < \Gamma(1/4)$, because the definition of the lhs here involves
integrating over a smaller portion of the positive real axis as compared to the rhs.
Thus, from~\eqref{eq:f(t)},
\begin{equation}
    f(t) < 2 t^{-1/2}.
\end{equation}
Applying this to~\eqref{eq: estimate partial sum} gives
\begin{equation}
    \sum_{n\leq N'} \chi_d(n)f(n) \ll d^{\,\ep} \log{N'}.
\end{equation}
Therefore, because we take the partial sums with $N'\leq N = O((d\, \Digits)^{1/2})$,
we have, on adjusting $\ep$ to incorporate the $\log{d}$:
\begin{equation}
    \sum_{n\leq N'} \chi_d(n)f(n) \ll d^{\,\ep} \log{\Digits}.
\end{equation}
Therefore, the partial sums do not get large and we thus have nearly as
many digits accuracy beyond the decimal place as our working precision.

Using a similar analysis, the effect of accumulated round off error can be estimated
by replacing $\chi_d(n)$ with random plus and minus ones
multiplied by a factor of size $10^{-\Digits}$ to model the random rounding up or down
of the terms in the sum. With high probability, we then get an error, due to accumulated
round off of size
\begin{equation}
    (d^{\,\ep} \log{\Digits}) 10^{-\Digits}.
\end{equation}

If one desires rigorous, rather than experimental, values of $L(1/2,\chi_d)$,
an interval arithmetic package should be used in practice. Because our goal was to test
conjectures rather than prove a rigorous numerical result, we were satisfied with an
intuitive understanding of the accuracy
of our computation and carried out several checks of the values attained, for example comparing
a similar smooth approximate functional equation for the case $d<0$ against select values attained by
our implementation using the Epstein zeta function, and also using a high precision version
of Rubinstein's lcalc package to test a few several values.

\subsubsection{Hacks}

We list a few hacks which were helpful in the implementation of the smooth
approximate functional equation (\ref{eq: pos-appformula}).
    \begin{itemize}
      \item $\chi_d(n)$ can be efficiently computed by repeatedly extracting
          powers of 2 and applying quadratic reciprocity.
      \item As in the case for $d<0$, it is to our advantage to partition
          $0< d \leq X$ into blocks and farm the work out to many processors.
      \item Due to the presence of $\chi_d(n)$ in the (\ref{eq: pos-appformula}),
          it is more efficient to place the $d$-loop on the inside and the
          $n$-loop on the outside because $\chi_d(n)$ is periodic in $d$ with
          period either $n$ or $8n$, depending on whether the power of two
          dividing $n$ is even or odd. Furthermore, $n$ is comparatively small
          compared to $d$, by~\eqref{eq: N}. Thus, for each $n$ we
          precomputed a table of $\chi_d(n)$, so as to only compute this values
          once per residue class $d$ mod $n$ or $8n$. This pays off so long as
          each residue class gets hit, on average, more than once (perhaps
          slightly more because of the overhead involved in storing the values
          and looking up the array.) In our implementation, with blocks of
          length $10^6$, $0< d < 1.3 \times 10^{10}$, and $16$ digits working
          precision, it was conducive to do so.

      \item We computed the normalized incomplete gamma function $G\pr{z,w}$,
          evaluated at $z=1/4$ and $w = n^2\pi/d$, as follows. For
          $w>37$, return $0$ (since $\exp\pr{-37} < 10^{-16}$). For $1 < w <
          37$, use a precomputed table of Taylor series, centering each Taylor
          series at multiples of $.01$ (so nearly $4000$ Taylor series) and
          taking terms up to degree $7$ (less for larger $w$ because of the
          exponential decay). Otherwise, for $w<1$, employ the complimentary
          incomplete gamma function $$\gm(z,w) := \Gm(z) - \Gm(z,w) = \int_0^w
          e^{-x} x^{z-1} dx, \qquad \tx{$\Re(z) > 0$, $\lr{\tx{arg} w} <
          \pi$}.$$ Specifically, set $$g(z,w) = w^{-z} \gm\pr{z,w} = \int_0^1
          e^{-wt}t^{z-1} dt,$$ so $G(z,w) = w^{-z}\Gm(z) - g(z,w)$, and
          integrate by parts to get $$g(z,w)  = e^{-w} \sum_{j=0}^\infty
          \fr{w^j}{\pr{z}_{j+1}},$$ where $$\pr{z}_j = \begin{cases}
          z\pr{z+1}\cdots \pr{z+j-1} \quad &\tx{if $j>0$};\\ 1 \quad
          &\tx{if $j=0$}. \end{cases}$$ We stored the value of $\Gm(1/4)$ and
          calculated the above series for $g(1/4,w)$
          by truncating the sum when the tail was less
          than $10^{-16}$.
    \end{itemize}

\subsection{Complexity for $d>0$}

Recall that, as for the case of negative discriminants, we are partitioning the interval
$0 < d < X$ into blocks of length $\Delta X$.

The overall cost for sieving for fundamental discriminants, summed over all
blocks, is a meager $O(X)$ arithmetic operations and array accesses on numbers
of bit length $O(\log{X})$, as for the case of $d<0$.

Next we estimate the overall time, summed over blocks, required to create a
precomputed table of characters $\chi_d(n)$ for all residue classes 
mod $n$ or $8n$.

Summing over blocks $m$, and taking the maximum truncation point~\eqref{eq: N}
that occurs for a given block, the time required is
\begin{equation}
    \ll \log(X)^2 \sum_{m\leq \fr{X}{\Dl{X}}}
    \sum_{n\leq M} n
\end{equation}
where
\begin{equation}
    M = \sq{\fr{m\Dl{X}}{\pi}\log(10)\cdot \Digits}.
    \label{eq: M}
\end{equation}
Here we have used the fact that
each character $\kr{d}{n}$ can be calculated in time
$\B\pr{\tx{size}(d)\tx{size(n)}}$, where size means binary length (see, for example,
\cite{hcohen: numth1}), and that both $d$ and $n$ are of size
$O(\log X)$ in this case.
Summing, the time needed here is therefore
\begin{equation}
    \ll \fr{X^2\log(X)^2\Digits}{\Dl{X}}.
\end{equation}

So, by choosing $\Dl{X}\gg X^{1/2+\ep}$, we can make the overall
time spent on computing the Kronecker symbol $o(X^{3/2})$.
As $\Dl{X}$ increases, there is a tradeoff between spending less time
on the character computation and
having larger arrays, similar to our computation of gcd's in the $d<0$ case.
There is a definite advantage, depending on the particular hardware, to having
smaller arrays, i.e. smaller
$\Delta X$, to reduce the number of calls to move data from RAM into cache.
On our hardware, and in our range $0<d<1.3 \times 10^{10}$, we found
that a value of $\Delta X = 10^6$ worked well.

Thus, the bulk of the work is spent on looping, for each block, through $n$ and
$d$, looking up the precomputed character values, computing the normalized
incomplete gamma function $G(1/4,\pi n^2/d)$ to given precision, and updating
the corresponding value of $L(1/2,\chi_d)$ by the amount $\chi_d(n) f(n)$.

The kind of work and operations required is thus very similar to our
approach for the $d<0$ case, with the handling of characters similar
to our handling of gcd's, and the approach to computing the incomplete gamma function
similar to that of the $K$-Bessel function.

However, there is one significant difference in the two methods.
For $d<0$, equation~\eqref{eq: number of triples} tells us that
our Epstein zeta function method
loops through $\frac{\pi}{36\zeta(3)} X^{3/2} =\approx 0.0726 X^{3/2}$
triples $a,b,c$. Not only is the constant $.0726$ small,
but the desired precision does not affect the number
of triples required. Precision
becomes a factor only in regards to computing the particular contribution from
each triple, for example the number of terms needed for the various $K$-Bessel Taylor
series expansions.

But, in the present case of the smooth approximate functional equation, both
the {\it length of the sum} and the amount of work needed to compute the individual
terms of the sum depends on the desired precision. So, the main difference in
these two approaches is the length of the sum.

In the case of $d>0$, the length of the main $d,n$
loops, summed over all blocks of length $\Dl{X}$, is quantified by
\begin{equation}
    L_{\tx{pos}} =
    \sum_{m \leq \fr{X}{\Dl{X}}}
    \sum_{n\leq M}
    \sum_{(m-1)\Dl{X} < d \leq m\Dl{X}} 1,
    \label{eq: Lpos defn}
\end{equation}
with $M$ given by (\ref{eq: M}). Simplifying the two inner sums,
this quantity is easily estimated to asymptotically be 
\begin{equation}
    (\Delta X)^{3/2} \sqrt{\fr{\log(10)\cdot \Digits}{\pi}} \sum_{m \leq \fr{X}{\Dl{X}}} \sqrt{m} 
    \;\sim\; \frac{2}{3} \sqrt{\fr{\log(10)\cdot \Digits}{\pi}} X^{3/2}
    \label{eq: Lpositive}
\end{equation}

So, if $\Digits=16$, then $L_{\tx{pos}} \approx 2.28 X^{\fr{3}{2}}$,
which is more than twenty times larger than the number of triples, $0.0726 X^{3/2}$,
considered for $d<0$.

It is impossible to precisely pin down, theoretically, the constant factor
savings in the runtime of our method for $d<0$ compared to the approach used
for $d>0$ as it depends on the speed of the various arithmetic and memory
operations on particular hardware. Furthermore, these are not easily
quantifiable as they change according to how the various resources of the
machine are being used at a given moment. Another obstacle to a precise
comparison is that one would need to take into account implementation choices
made by the programmer and also by the compiler at the minutest of levels.

Nonetheless, the rough comparison between the lengths of the main loops
involved, i.e. $2.28 X^{\fr{3}{2}}$ for $d>0$ and $0.0726 X^{3/2}$ for $d<0$
(see equation~\eqref{eq: number of triples}), does reflect the different
runtimes when compared experimentally.

We ran our computation for $d<0$ on {\tt mod.math.washington.edu}, which is a
Sun Fire X4450, dating from 2008, with 24 Intel Xeon X7640 2.66 GHz CPUs (we
used 12 of them), and 128 GB RAM. Our computation for $d>0$ was carried out
on {\tt pilatus.uwaterloo.ca} which is an older SGI Altix 3700
machine, dating from around 2003, with 64 Intel Madison Itanium CPUs (we used
55 of these) running at 1.3 GHz, and 192 GB of RAM.

Our computation, for $d>0$, took roughly $18.9$ CPU years, and
about $3.9$ CPU years for $d<0$. Recall that we went up to $0 < - d < 5\times10^{10}$,
whereas for $d>0$, we managed to get to $1.3\times10^{10}$. So not only
did our computation for $d<0$ take much less CPU time, but we went significantly further.
To make this more meaningful, we should compare intervals of similar length, i.e.
the subset of $0 < - d < 1.3\times 10^{10}$. This interval required
$.4$ CPU years, i.e. about 47 times faster than our computation for the interval
$0 <  d < 1.3\times 10^{10}$. However, because different machines were used for
$d<0$ and $d>0$,
we should compensate by dividing the time used for $d>0$ by a factor of $2.5$ to
account for the fact that these $d$ were handled on an older and slower machine.
The value of $2.5$ was decided by rerunning select blocks of $d>0$ on both
machines, using the same {\tt C++} code, and comparing their runtimes, which
were about 2-3 times faster on the newer machine. Therefore, dividing $47$ by
$2.5$, our code ran around 20 times faster for $d<0$ than
it did for $d>0$, consistent with our rough expectations based on the
lengths in both methods of the main loops.

\addcontentsline{toc}{chapter}{

\end{document}